\documentclass[11pt,leqno]{article}
\usepackage[T1]{fontenc}
\usepackage[utf8]{inputenc}
\usepackage{lmodern,amsmath,amssymb,amsthm,mathtools,microtype,geometry,booktabs,array,xcolor}
\usepackage{comment}
\definecolor{darkgreen}{RGB}{0,90,48}
\usepackage{tikz}
\usetikzlibrary{arrows.meta}
\usetikzlibrary{calc}
\usepackage[authoryear,round]{natbib}
\usepackage[colorlinks=true,linkcolor=blue,citecolor=blue,urlcolor=blue]{hyperref}
\mathtoolsset{showonlyrefs=true,showmanualtags=false}
\allowdisplaybreaks[2]
\newtheorem{theorem}{Theorem}[section]
\newtheorem{proposition}[theorem]{Proposition}
\newtheorem{lemma}[theorem]{Lemma}
\newtheorem{corollary}[theorem]{Corollary}

\theoremstyle{remark}
\newtheorem{remark}[theorem]{Remark}
\theoremstyle{remark}
\newtheorem{example}[theorem]{Example}
\newcommand{\R}{\mathbb R}
\newcommand{\E}{\mathbb E}
\newcommand{\Prb}{\mathbb P}
\newcommand{\Lip}{\operatorname{Lip}}

\newcommand{\clip}{\operatorname{clip}}

\newcommand{\GS}{\operatorname{GS}}
\newcommand{\dd}{\,\mathrm d}
\newcommand{\ind}{\mathbf 1}

\hypersetup{pdfauthor={},pdfsubject={GroupSort approximation, covering numbers, nonparametric regression and flow matching}}
\hypersetup{pdftitle={Approximation and estimation with Lipschitz-constrained GroupSort networks}}
\begin{document}
\title{Approximation and estimation\\with Lipschitz-constrained GroupSort2 networks}
\author{Alberto Gonz\'alez-Sanz, Lea Kunkel and Sophie Langer}
\date{}
\maketitle
\vspace{-2em}
\begin{abstract}
Approximating smooth targets by neural networks does not automatically control the Lipschitz constant of the approximating network. In flow matching, the effect of velocity estimation error on the generated distribution can grow exponentially with the time integral of the Lipschitz constant of the fitted field. We study fully connected GroupSort2 networks, whose activation sorts coordinates in pairs and is \(1\)-Lipschitz in the maximum norm. By requiring every hidden weight matrix to have absolute row sums at most one, we bound the Lipschitz constant of the network by the norm of its output row, independently of its width and depth. For $\alpha$-H\"older functions on \([0,1]^q\), we derive approximation rates for corresponding GroupSort2 network classes with these parameter constraints. Here, the concrete architecture depends on the smoothness of the target function. A corresponding covering number bound shows the optimality of the derived rates in case of $\alpha \geq 1$. 

We apply our results to flow matching under a Gaussian density-ratio model with smoothness \(\alpha>1\), and to nonparametric least-squares regression. For flow matching, we construct network classes with a common spatial Lipschitz bound independent of the sample size and derive the minimax optimal Wasserstein-\(2\) rate (up to logarithmic factors) for the corresponding flow estimator. For regression with bounded responses or conditionally subgaussian errors, the constrained least-squares estimator attains the classical Hölder minimax rate in squared prediction error up to logarithmic factors.
\end{abstract}

\tableofcontents

\section{Introduction}\label{sec:intro}
Uniform approximation of a Lipschitz function does not by itself control the
Lipschitz constant of the approximating network. In flow matching, this
distinction has a direct consequence for distribution estimation. To generate
a sample, one draws an initial point from a reference distribution and solves
an ordinary differential equation driven by an estimated velocity field. An
error in the velocity can displace the trajectory and thereby change the
points at which the fitted field is evaluated later. If $L_t$ denotes the
spatial Lipschitz constant of the fitted field at time $t$, a standard bound
on the Wasserstein distance to the target contains a factor of the form
$\exp\bigl(\int_0^1 L_t\,dt\bigr)$ multiplying the velocity estimation error
\citep{AlbergoVandenEijnden2023, BentonEtAl2024}. In turn, the dependence of this stability factor on the sample size must be
retained when deriving statistical rates.

A natural way to enforce a prescribed Lipschitz bound is to constrain the product of the operator norms of the weight matrices in the network.
The restriction of such a constraint on approximation is already visible in one
dimension. If the product of the induced $\ell_\infty$ operator norms is allowed to be at
most one, ReLU networks cannot approximate $x\mapsto |x|$ arbitrarily well
on $[-1,1]$ \citep[Section~4.2]{HusterChiangChadha2018}. A network that
sorts coordinates in pairs represents the same function exactly under the
same constraint. It suffices to sort $(x,-x)$ and select the larger
coordinate. This pairwise sorting operation is the GroupSort2 activation.
\citet{AnilLucasGrosse2019} used the ability to take maxima and minima to
show that, under suitable layerwise norm constraints, GroupSort2 networks
can uniformly approximate every real-valued $1$-Lipschitz function on a
compact domain. For smooth targets, this raises a quantitative question.
{\it How do the network width and depth required for a given accuracy depend on the
smoothness of the target?}

\subsection{Main contributions}\label{sec:contributions}
We analyze the approximation and statistical properties of GroupSort2 networks under explicit constraints on their weight matrices. The derived approximation and complexity bounds are applied to flow matching and nonparametric regression. Our main contributions are as follows.
\paragraph{\it Approximation of smooth functions under Lipschitz constraints.}
We derive a new approximation result for target functions in a bounded H\"older class of order $\alpha\geq 1$ on $[0,1]^q$ using GroupSort2 networks with prescribed hidden-layer
matrix constraints. In particular, for every sufficiently small
$\epsilon>0$, we construct a fully connected GroupSort2
network $g_\epsilon$ such that
\[
 \|f-g_\epsilon\|_{\infty,[0,1]^q}\lesssim\epsilon,
 \qquad
 \operatorname{Lip}(g_\epsilon)\le C,
 \qquad
 S_\epsilon\lesssim\epsilon^{-q/\alpha},
\]
where $S_\epsilon\asymp w^2L$ denotes the number of weights and biases.
The constant $C$ is independent of $\epsilon$ and uniform
over the target class. For $\alpha>1$, we establish a range of widths
$w$ and depths $L$ attaining error
$O((w^2L)^{-\alpha/q})$ while preserving a common Lipschitz
bound. When $\alpha>1$ and $2\alpha\ge q+1$, a depth of
order $[\log(ew)]^2$ suffices. At $\alpha=1$, the same
approximation rate with a bounded Lipschitz constant is
achieved using fixed width and depth
$L\asymp\epsilon^{-q}$, see Theorem~\ref{th:holder-approximation}.
\paragraph{\it Entropy bound of constrained GroupSort2 networks.} We show that, when the parameter magnitudes grow at most polynomially in $\epsilon^{-1}$, the corresponding GroupSort2 class $\mathcal F_\epsilon$ satisfies an entropy bound
\[
 \log \mathcal{N}(\epsilon,\mathcal F_\epsilon,\|\cdot\|_{\infty,[0,1]^q})
 \lesssim S_{\epsilon}\log\!\left(\epsilon^{-1}\right)
 \lesssim \epsilon^{-\frac{q}{\alpha}}\log\!\left(\epsilon^{-1}\right), 
\]
see Lemma~\ref{lem:network-entropy}. Comparison with the metric entropy of the $\alpha$-H\"older ball on $[0,1]^q$ shows that our derived approximation rate is optimal up to logarithmic factors among network classes satisfying the same covering estimate.
\paragraph{\it Minimax rates for distribution estimation by flow matching.}  We apply our results to flow matching and estimate the underlying velocity field by least squares over a class $\mathcal G_n$ of GroupSort2 networks. Combining the matrix constraints with a dyadic partition of the time interval, we construct $\mathcal G_n$ so that
\[
  \sup_{g\in\mathcal G_n}
  \int_0^1 \operatorname{Lip}_x(g_t)\,dt \leq C,
\]
where $\operatorname{Lip}_x$ denotes the spatial Lipschitz constant and $C$ is independent of $n$. Consequently, the stability factor relating velocity estimation error to Wasserstein error remains bounded for the network selected from the data, even as $\mathcal G_n$ grows.

Under the Gaussian density-ratio model, see \eqref{eq:model} for a precise definition, with smoothness \(\alpha>1\) and dimension $q$, the resulting distribution estimator $\widehat P_n$ satisfies
\[
\sup_{P\in\mathcal P_\alpha}
\mathbb E_P W_2(\widehat P_n,P)
=\widetilde O\!\left(
n^{-\min\{\frac12,\frac{\alpha+1}{2\alpha+q}\}}
\right).
\]
This rate is minimax optimal up to logarithmic factors (\citet{NilesWeedBerthet2022}) and is attained with at most
\(\widetilde O(n^{q/(2\alpha+\max\{q,2\})})\) parameters.
\paragraph{\it Minimax rates for nonparametric regression.} We analyze GroupSort2 networks as estimators in nonparametric regression and show that for a regression function $f_0$ in a bounded H\"older class of order $\alpha>0$ on $[0,1]^q$, the least-squares estimator over the constrained GroupSort class satisfies
\[
  \mathbb E\|\widehat f_n-f_0\|_{L^2(P_X)}^2
  = \widetilde O\!\left(n^{-\frac{2\alpha}{2\alpha+q}}\right), 
\]
see Theorem~\ref{thm:regression}.  The estimator thus attains the classical minimax convergence rate of \citet{Stone1982} under our smoothness assumptions. The upper bound holds for every covariate distribution $P_X$ supported on $[0,1]^q$. Corollary~\ref{cor:regression-subgaussian} extends the result to subgaussian responses.

\subsection{Related work}\label{sec:related-work}

\textit{GroupSort networks and approximation under norm constraints.}
Early work on GroupSort \citep{AnilLucasGrosse2019} and FullSort
\citep{CohenEtAl2019} showed that sorting activations permit universal
approximation of $1$-Lipschitz functions on compact domains under
suitable layerwise norm constraints. \citet{TanielianSangnierBiau2021} derived corresponding quantitative results on the network size.  GroupSort networks are applied in the context of robust classification by \citet{serrurier2021robustness}, who trained
$1$-Lipschitz networks using a hinge-regularized optimal transport
loss, and \citet{bethune2022loss}, who studied how the choice of loss
affects the trade-off between accuracy and robustness.

For ReLU networks, much of the approximation theory analyzes how network accuracy depends on the regularity and dimension of the target function \citep{Yarotsky2017,GuhringEtAl2020} and how network width can
be traded for depth \citep{NaglerLanger2026,Nagler2026}. Those results only focus on the approximation accuracy of the networks without accounting for the Lipschitz constant of the approximating network. Simultaneous control of the Lipschitz constant was derived by \cite{HongKratsios2024}, who showed that sparse ReLU
networks can also inherit the Lipschitz bound of their targets. \citet{JiaoWangYang2023} further state that when enforcing Lipschitz control through the network weights, the available norm budget enters the approximation rate. Our analysis reveals how accuracy improves with smoothness
at a given parameter count when each row of every hidden matrix has
$\ell_1$-norm at most one. 
\medskip
\\
\textit{Flow matching and distribution estimation.}
Flow matching fits a velocity field by least squares and integrates the
fitted field to obtain a pushforward map that transforms a given latent distribution into a distribution estimator \citep{LiuGongLiu2023,LipmanEtAl2023, AlbergoEtAl2025}.

To evaluate the performance of the estimated  distribution $\widehat P_n$ in Wasserstein distance, a common approach (cf.~\cite[Lemma~1]{Stephanovitch2026}) is  to use the inequality 
\begin{equation}\label{eq:stability-intro}
 W_2(\widehat P_n,P)
 \le\int_0^1
 \exp\!\left(\int_t^1\Lip_x(\widehat v_{n,s})\dd s\right)
 \|\widehat v_{n,t}-v_t^P\|_{L^2(P_t)}\dd t,  
\end{equation}
with $v_t^P$ the population vector field, which transfers the approximation rate of the fitted velocity $\widehat v_{n,s}$ to the estimated  distribution $\widehat P_n$.
To bound this term, one needs, besides an accurate estimator, a bound on the Lipschitz constant of the population vector field and one for its estimator, which is commonly a neural network. Without this control, a gap in the analysis occurs as, e.g., in the approach of \cite{FukumizuEtAl2025}, pointed out by \cite{KunkelTrabs2026}. \cite{KunkelTrabs2026} solve this problem by considering highly overparameterized neural networks and connect their flow matching analysis to a kernel density estimation problem. However, this analysis does not explain the outperformance of flow matching over kernel density estimators in empirical results. A bounded Lipschitz constant of the population velocity field can be derived by imposing regularity properties on the target distribution  \citep{GaoEtAl2024, Stephanovitch2026, Kunkel2026, GaoEtAl2026}. One setting is the Gaussian density-ratio model, which we consider in the following for simplicity.

A bound on the Lipschitz constant of the network estimator can be derived by, e.g., taking ReLU networks that simultaneously approximate functions and their Lipschitz constants, see \cite{GuhringEtAl2020}. However, this result requires comparably large networks, which leads to suboptimal rates for the Wasserstein distance of the distribution estimators \citep{ZhouLiu2025, Kunkel2026}. Our GroupSort2 class combines efficient approximation of the
population vector field with a uniform spatial Lipschitz bound
enforced by explicit layerwise matrix-norm constraints.
This yields the minimax $W_2$ rate up to logarithmic factors
on our Gaussian density-ratio class.
Parallel to our work, the recent revision of \citet{GaoEtAl2026}
constructs a finite ReLU sieve with polynomial-size networks and
uniform spatial Lipschitz control.
Their admissible class is defined through constraints on the
realized functions, without an efficient procedure for finding these networks. Their $W_2$ upper bound is not matched by their minimax lower bound.

\paragraph{\it Nonparametric regression.}
In nonparametric regression, the minimax optimal convergence rate over an $\alpha$-H\"older ball in $q$
variables is $n^{-2\alpha/(2\alpha+q)}$ for squared prediction error
\citep{Stone1982}. For networks, these optimal rates were shown to hold for sparsely connected \citep{SchmidtHieber2020} and fully connected ReLU networks \citep{KohlerLanger2021}, including those of flexible width and depth \citep{NaglerLanger2026, JiaoEtAl2023}. We show that the classical H\"older rate survives the layerwise constraints of our GroupSort network class. In parallel to our work, the very recent work of \citet{ShiEtAl2026}\footnote{The comparison refers to arXiv:2611.23979v1.
This version contains a proof sketch, but the detailed
proofs are deferred to supplementary appendices that
are absent from both the public PDF and the deposited
source archive. Consequently, the available arguments
do not permit a complete verification that the stated
results hold under the claimed assumptions.} studies penalized least-squares estimators
with ReLU networks, using a split spectral norm product as the penalty. For sufficiently smooth target functions with $\alpha > q/2+1$, they obtain the
optimal rate up to logarithmic
factors for networks with bounded expected Lipschitz constant. In comparison, we analyze an unpenalized empirical risk minimization over GroupSort2
classes with prescribed layerwise norm constraints and flexible width-depth architecture for target functions with arbitrary smoothness.

\paragraph{\it Outline.} In Section~\ref{sec:groupsortnn} we introduce the class of GroupSort2 networks. Section~\ref{sec:approximation} provides their approximation rates for $\alpha$-H\"older functions and the corresponding covering number bounds. Sections~\ref{sec:flow-method} and~\ref{sec:regression-method} apply them in the context of flow matching and nonparametric regression, respectively. The proofs are deferred to the appendix.

\paragraph{\it Notation.} For a vector $x \in \mathbb{R}^q$, write $|x|$ and $\|x\|_\infty$ for the Euclidean and maximum norms, respectively. For a function $f$ on $\Omega$ its sup-norm and $L_2$-norm are defined as $\|f\|_{\infty,\Omega}=\sup_{x\in\Omega}|f(x)|$ and $\|f\|_{L^2(P)}^2=\int|f|^2\dd P$, respectively. We omit $\Omega$ when it is the full domain. For a time-dependent field $g$ we write $g_t(x)=g(t,x)$ and denote its Euclidean spatial Lipschitz constant by $\Lip_x(g_t)$. We use $\Lip$ for the Euclidean Lipschitz constant and $\Lip_\infty$ when both the domain and codomain carry the maximum norm. We write $a_+=\max\{a,0\}$, $\ind_A$ for the indicator of $A$, and $\clip_D(x)=(\min\{\max\{x_i,a_i\}, b_i\})_{i=1}^q$ for $x \in \mathbb{R}^q$ coordinatewise clipping to a rectangle $D=\prod_{i=1}^q [a_i,b_i]$. The implicit constants in $\lesssim$, $\gtrsim$ and $\asymp$ may depend on fixed dimensions, smoothness and model bounds. 

For a multi-index $\beta$, let $|\beta|=\sum_i\beta_i$, $\beta!=\prod_i\beta_i!$, and $D^\beta=\partial_1^{\beta_1}\cdots\partial_q^{\beta_q}$. On a closed rectangle, derivatives extend continuously to its boundary. For $0<\rho<1$, write $[f]_{\rho,\Omega}=\sup_{x\ne y\in\Omega}\frac{|f(x)-f(y)|}{|x-y|^\rho}$. For noninteger $\alpha>0$ and $s=\lfloor\alpha\rfloor$, set
\begin{equation}\label{eq:holder-norm-definition}
 \|f\|_{C^\alpha(\Omega)}=
 \max\left\{\max_{|\beta|\le s}\|D^\beta f\|_{\infty,\Omega},
             \max_{|\beta|=s}[D^\beta f]_{\alpha-s,\Omega}\right\}.
\end{equation}
At integer order $s\ge0$, use
\begin{equation}\label{eq:integer-norm-definition}
 \|f\|_{C^s(\Omega)}=\max_{|\beta|\le s}\|D^\beta f\|_{\infty,\Omega}.
\end{equation}
For vector-valued functions we take the maximum of the componentwise $C^\alpha$ norms. All unqualified logarithms are natural.
We write $a_n=\widetilde O(b_n)$ if
$a_n=O(b_n[\log(en)]^c)$ for some fixed $c\ge0$.

\section{GroupSort neural networks}\label{sec:groupsortnn}\label{sec:networks}
For integers $r\ge2$ and $m$ divisible by $r$, the GroupSort
activation $\GS_r:\mathbb R^m\to\mathbb R^m$ partitions its input
into $m/r$ consecutive blocks
$\bigl(x_{(i-1)r+1},\ldots,x_{ir}\bigr),$ for $ i=1,\ldots,m/r,$ 
and sorts the coordinates within each block in nonincreasing order
\citep{AnilLucasGrosse2019}. In the following, we suppress the dimension $m$ in the
notation. 

Let $L \in \mathbb{N}_0$ and $k=(k_0,\ldots,k_{L+1})\in\mathbb N^{L+2}$.
A fully connected GroupSort network with architecture $(L,k)$
has input dimension $k_0$, output dimension $k_{L+1}$, and $L$
hidden layers of widths $k_1,\ldots,k_L$, each divisible by $r$.
We call $L$ the depth and $w=\max\{1,k_1,\ldots,k_L\}$ the
maximum hidden width. Its affine maps are
\[
  a_\ell(z)=W_\ell z+b_\ell,
  \qquad
  W_\ell\in\mathbb R^{k_{\ell+1}\times k_\ell},
  \quad
  b_\ell\in\mathbb R^{k_{\ell+1}},
  \qquad \ell=0,\ldots,L.
\]
Writing $\theta=((W_\ell,b_\ell))_{\ell=0}^{L}$ for the collection
of network parameters, the network is defined as 
\begin{equation}\label{GSnet}
  g_\theta(x)
  =
  \bigl(
    a_L\circ\GS_r\circ a_{L-1}
    \circ\cdots\circ\GS_r\circ a_0
  \bigr)(x).
\end{equation}
For $L=0$, \eqref{GSnet} becomes an affine map
$g_\theta(x)=a_0(x)=W_0x+b_0$.

\begin{figure}[htbp]
\centering
\begin{tikzpicture}[
  x=1cm,y=1cm,
  stage/.style={font=\small\bfseries,text=black,inner sep=0pt},
  mathlabel/.style={font=\small,text=black,inner sep=1pt},
  neuron/.style={circle,draw=black,fill=white,line width=.6pt,
                 minimum size=5.8mm,inner sep=1pt,text=black},
  weight/.style={draw=black,line width=.55pt},
  selected/.style={draw=black,line width=1.35pt},
  feed/.style={draw=black,line width=.55pt,-{Stealth[length=1.25mm,width=1.05mm]}},
  sortbox/.style={draw=black,fill=white,line width=.65pt,rounded corners=2pt}
]

\foreach \left/\right/\middle/\title in {
  -.50/.56/.03/Input, 1.37/5.60/3.485/Hidden layer 1,
  6.68/10.73/8.705/Hidden layer 2, 11.70/13.23/12.465/Output}{
  \node[stage] at (\middle,2.22) {\title};
  \draw[black,line width=.4pt] (\left,1.96) -- (\right,1.96);
}

\node[neuron] (x1) at (.03,.60) {$x_1$};
\node[neuron] (x2) at (.03,-.60) {$x_2$};
\foreach \i/\y in {1/1.37,2/.69,3/-.69,4/-1.37}
  \node[neuron] (z\i) at (1.82,\y) {$z_{\i}$};
\foreach \i in {1,2}
  \foreach \j in {1,2,3,4}
    \draw[weight] (x\i) -- (z\j);

\foreach \i/\j/\ya/\yb in {1/2/1.37/.69,3/4/-.69/-1.37}{
  \path[sortbox] (2.78,{\yb-.28}) rectangle (5.58,{\ya+.28});
  \draw[black,line width=.35pt] (2.89,{(\ya+\yb)/2}) -- (5.47,{(\ya+\yb)/2});
  \draw[feed] (z\i.east) -- (2.78,\ya);
  \draw[feed] (z\j.east) -- (2.78,\yb);
  \node[mathlabel] at (4.18,\ya) {$h_{\i}=\max(z_{\i},z_{\j})$};
  \node[mathlabel] at (4.18,\yb) {$h_{\j}=\min(z_{\i},z_{\j})$};
  \coordinate (h\i) at (5.58,\ya);
  \coordinate (h\j) at (5.58,\yb);
}

\node[neuron] (u1) at (7.10,.58) {$u_1$};
\node[neuron] (u2) at (7.10,-.58) {$u_2$};
\foreach \i in {1,2,3,4}
  \foreach \j in {1,2}
    \draw[weight] (h\i) -- (u\j);
\draw[selected] (h1) -- (u1);
\draw[selected] (h3) -- (u2);

\path[sortbox] (7.90,-.86) rectangle (10.70,.86);
\draw[black,line width=.35pt] (8.01,0)--(10.59,0);
\draw[feed] (u1.east) -- (7.90,.58);
\draw[feed] (u2.east) -- (7.90,-.58);
\node[mathlabel] at (9.30,.58) {$v_1=\max(u_1,u_2)$};
\node[mathlabel] at (9.30,-.58) {$v_2=\min(u_1,u_2)$};
\node[neuron,minimum size=9mm] (output) at (12.465,0) {$g(x)$};
\draw[weight] (10.70,-.58) -- (output);
\draw[selected] (10.70,.58) -- (output);

\node[mathlabel] at (.03,-2.05) {$x\in\R^2$};
\node[mathlabel] at (3.485,-2.05) {$z=W_0x+b_0,\quad h=\GS(z)$};
\node[mathlabel] at (8.705,-2.05) {$u=W_1h+b_1,\quad v=\GS(u)$};
\node[mathlabel] at (12.465,-2.05) {$g(x)=W_2v+b_2$};

\draw[black,line width=.4pt] (-.50,-2.46)--(13.23,-2.46);
\draw[selected] (-.43,-2.86)--(.13,-2.86);
\node[font=\small\bfseries,text=black,anchor=west,inner sep=0pt]
  at (.31,-2.86) {Example};
\node[mathlabel,anchor=west] at (2.0,-2.86)
  {$u_1=h_1,\quad u_2=h_3,\quad g(x)=v_1$
   \quad$\Longrightarrow$\quad
   $g(x)=\max\{z_1,z_2,z_3,z_4\}$.};

\end{tikzpicture}
\caption{GroupSort2 network with $L=2$ hidden layers and architecture $k=(2,4,2,1)$.}
\label{fig:groupsort}
\end{figure}
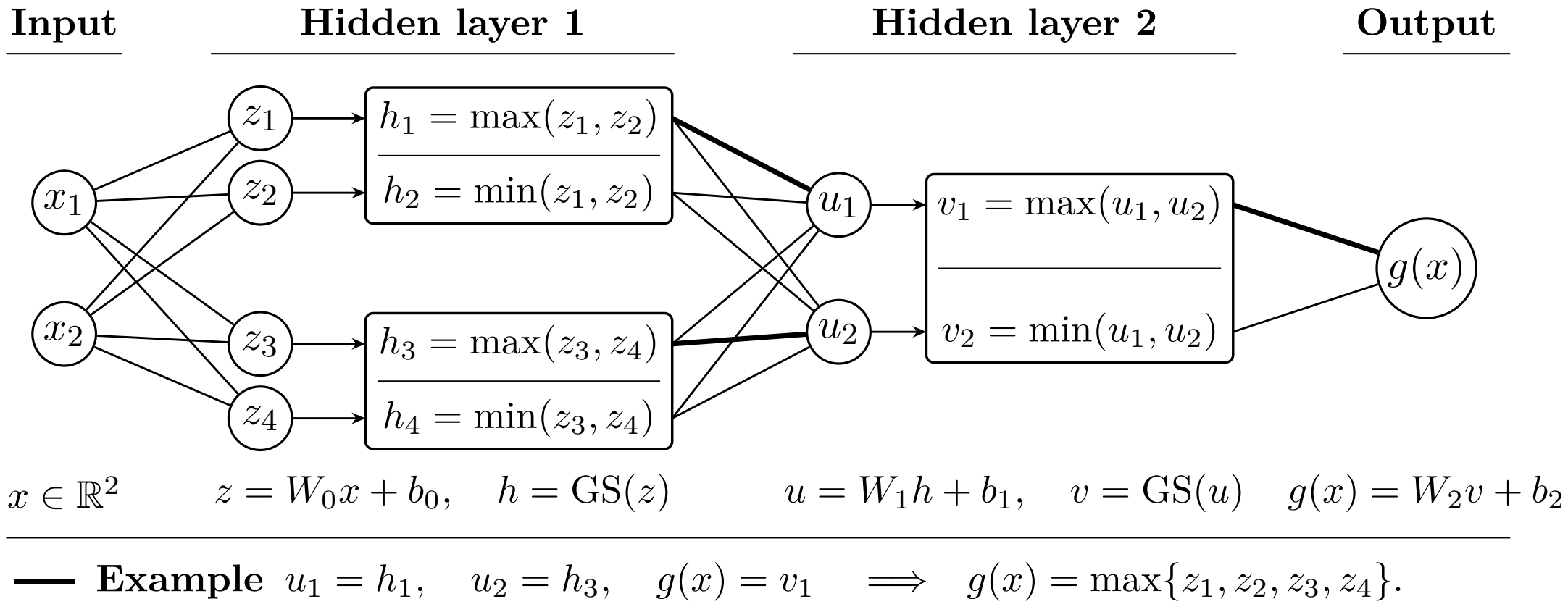
Throughout this paper, we consider GroupSort2 networks, i.e., $r=2$ and we write $\GS=\GS_2$.
Thus all hidden-layer widths are even, and the activation acts
on each adjacent pair according to
\[
  \GS(u,v)
  =
  \bigl(\max\{u,v\},\,\min\{u,v\}\bigr),
\]
see Figure~\ref{fig:groupsort} for a graphical description of these networks. 
The number of weights and biases of an architecture $(L,k)$ is
\begin{equation}\label{eq:dense-parameter-count}
 S(L,k)=\sum_{\ell=0}^{L}k_{\ell+1}(k_\ell+1).
\end{equation}

For scalar-output networks with input dimension $q$ and all $L\ge1$ hidden layers of width $w$, the number of parameters is
\begin{equation}\label{eq:dense-tradeoff}
 S(L,w)=(q+1)w+(L-1)w(w+1)+w+1 \asymp w^2(L+1)
\end{equation}
for $L\ge2$. For hidden widths at most $w$, we have $S(L,k)\le S(L,w)$.

To distinguish entrywise bounds from operator-norm bounds, we write
\[
  \|W\|_{\max}:=\max_{i,j}|W_{ij}|,
  \qquad
  \|\theta\|_{\max}
  :=
  \max_{0\le\ell\le L}
  \max\bigl\{\|W_\ell\|_{\max},\,\|b_\ell\|_\infty\bigr\}.
\]
For $b,B>0$, we impose the parameter bound
$\|\theta\|_{\max}\le b$ and the layerwise constraints
\begin{equation}\label{eq:network-constraints}
  \|W_\ell\|_{\infty\to\infty} =\max_i\sum_j|(W_\ell)_{ij}|
  \le
  \begin{cases}
    1, & \ell=0,\ldots,L-1,\\
    B, & \ell=L.
  \end{cases}
\end{equation}
As $\GS_2$ is nonexpansive
with respect to the maximum norm, i.e.,
\[
  \|\GS_2(x)-\GS_2(y)\|_\infty
  \le \|x-y\|_\infty,
  \qquad x,y\in\mathbb R^m,
\]
we then have
\begin{equation}
\label{eq:layer-norm-lipschitz}
\Lip_\infty(g_\theta)\le\prod_{\ell=0}^L\|W_\ell\|_{\infty\to\infty} \leq B.
\end{equation}
For $L=0$, only the output-layer constraint applies.
Write $\widetilde{\mathcal F}(L,w,b,B)$ for scalar-valued networks of the form \eqref{GSnet} of depth $L$, hidden widths at most $w$, $\|\theta\|_{\max} \leq b$ and parameters satisfying \eqref{eq:network-constraints}. The input dimension is specified by the domain. For $H>0$ define
\begin{equation}\label{eq:regression-class-definition}
 \mathcal F(L,w,b,B,H)
 =\{\clip_{[-H,H]}\circ g:g\in\widetilde{\mathcal F}(L,w,b,B)\},
\end{equation}
which cuts the output of each network in absolute value by $H$. As the clipping operator is 1-Lipschitz, networks $h \in \mathcal F(L,w,b,B,H)$ preserve the bound $\Lip_\infty(h)\le B$.

\section{Approximation and covering numbers}\label{sec:approximation}\label{sec:methodology}

 We first state the approximation result of GroupSort2 networks for H\"older functions and then derive their corresponding uniform covering number bound.

\subsection{Approximation of H\"older functions}\label{sec:holder-approximation}
Depending on the smoothness $\alpha$ of the H\"older function, the network architecture has to be chosen differently. For $0<\alpha<1$, we keep the width of the network fixed and improve approximation power by solely increasing the depth. For $\alpha\ge1$, the network class is more flexible and we can
also use increasing width to achieve the desired accuracy. We track
the norm of the output row explicitly because, under our hidden-layer
constraints, it bounds the Lipschitz constant of the network (see \eqref{eq:network-constraints}).

\begin{theorem}[Uniform approximation of H\"older functions]
\label{th:holder-approximation}
Let $q\ge1$ be an integer, $\alpha>0$, and $H\ge1$.
Suppose that $\|f\|_{C^\alpha([0,1]^q)}\le H$.
Let $0<\epsilon\le1$, let $w$ be a sufficiently large even
integer, and let $L\ge2$ be an integer. For sufficiently large constants $C_B$ and $C_b$, set
\begin{equation}\label{eq:holder-network-budgets}
 B_{\epsilon,w}
 = C_B
 \begin{cases}
   \epsilon^{-(1-\alpha)/\alpha}, & 0<\alpha<1,\\[3pt]
   \log(e/\epsilon)+w^{1-1/q},    & \alpha=1,\\[3pt]
   1,                             & \alpha>1,
 \end{cases}
 \qquad
 b_{\epsilon,w}=C_b B_{\epsilon,w}.
\end{equation}
Suppose that one of the following holds.
\begin{enumerate}
\item[(i)] $0<\alpha<1$, the width $w$ is fixed, and
  $L\gtrsim\epsilon^{-q/\alpha}$.
\item[(ii)] $\alpha\ge1$ and
  \begin{equation}\label{eq:holder-joint-range}
     L\gtrsim[\log(ew)]^2,
     \qquad
     w^2L\gtrsim\epsilon^{-q/\alpha}.
  \end{equation}
  If $1<\alpha<(q+1)/2$, assume in addition that
  \[
     L\ge w^{(q+1-2\alpha)/(\alpha-1)}.
  \]
\end{enumerate}
Then there exists a GroupSort2 network
\begin{equation}\label{eq:holder-approximation}
 g\in\mathcal F(L,w,b_{\epsilon,w},B_{\epsilon,w},H)
 \qquad\text{such that}\qquad
 \|f-g\|_{\infty,[0,1]^q}\lesssim\epsilon.
\end{equation}
For $\alpha=1$ and $L\gtrsim\epsilon^{-q}$, the conclusion
also holds with the smaller parameter bounds
$B_{\epsilon,w}=C_B$ and $b_{\epsilon,w}=C_bC_B$,
independently of $\epsilon$.
\end{theorem}

\begin{remark}[Admissible architecture choices]
\label{rem:holder-architecture-choices}
For $\alpha\ge1$ and sufficiently small $\epsilon>0$,
an admissible choice is
\[
 L_\epsilon\asymp
 \begin{cases}
   \epsilon^{-q},
       & \alpha=1,\\[3pt]
   \epsilon^{-\frac{q(q+1-2\alpha)}{\alpha(q-1)}},
       & 1<\alpha<\frac{q+1}{2},\\[5pt]
   [\log(e/\epsilon)]^2,
       & \alpha>1,\quad 2\alpha\ge q+1,
 \end{cases}
 \qquad
 w_\epsilon\asymp
 \frac{\epsilon^{-q/(2\alpha)}}{\sqrt{L_\epsilon}}.
\]
The multiplicative constants can be chosen so that the conditions of the 
theorem hold. The intermediate regime is empty
when $q=1$. These choices give a total parameter count
\[
 S_\epsilon
 \asymp w_\epsilon^2(L_\epsilon+1)
 \asymp \epsilon^{-q/\alpha},
\]
with parameter and Lipschitz bounds independent of $\epsilon$.
At $\alpha=1$, the width remains fixed and we use the
constant bounds in the final assertion of the theorem.
\end{remark}

At fixed width, the parameter count is proportional to the
depth. Theorem~\ref{th:holder-approximation} therefore gives
uniform approximation error $O(S^{-\alpha/q})$ for every
$\alpha>0$. For $0<\alpha<1$, this results in a growing output-row
bound in \eqref{eq:holder-network-budgets}. For $\alpha>1$,
the same rate holds throughout the stated joint width--depth
range with fixed parameter bounds.

For the choices in Remark~\ref{rem:holder-architecture-choices},
the Lipschitz bound holds uniformly over the corresponding
network class. Uniform bounds on individual parameters,
as considered, e.g., in \citet{Nagler2026}, do not by themselves
provide a Lipschitz bound independent of network size. Section~\ref{sec:entropy} shows that the total parameter
count is optimal up to a logarithmic factor under the
specified constraints and parameter bounds. 

The rate in Theorem~\ref{th:holder-approximation}(ii) is expressed in terms of the total parameter
count $S\asymp w^2L$, yet the theorem also requires the depth to grow
with the width. For ReLU networks with uniformly bounded individual weights, \citet[Theorem~3.1]{Nagler2026} shows that also for fixed depth and increasing width an arbitrarily small approximation error can be achieved. The following example shows why some
growth in depth is unavoidable under our hidden-layer constraints and
the output-row bound in \eqref{eq:holder-network-budgets}. Already in
one dimension, the $C^1$ ball contains functions whose slopes remain
bounded by one while their oscillations become increasingly frequent.

\begin{example}
Under our hidden-layer norm constraints, approximating a fixed
$C^1([0,1])$ ball to uniform accuracy $\epsilon$ requires
\[
    2^L B \gtrsim \epsilon^{-1},
\]
where $B$ bounds the $\ell_1$-norm of the output row. In particular,
$L\gtrsim\log(1/\epsilon)$ both for bounded $B$ and for the
logarithmic budget in \eqref{eq:holder-network-budgets} at
$q=\alpha=1$.

To see this, take $m=\lfloor(8\pi\epsilon)^{-1}\rfloor$ for
sufficiently small $\epsilon$ and set
$f_m(x)=(2\pi m)^{-1}\sin(2\pi m x)$. These functions lie in a
fixed $C^1$ ball since $\|f_m'\|_\infty=1$. Their successive extrema
are $1/(2m)$ apart and have alternating values $\pm(2\pi m)^{-1}$.
Thus, if $\|g-f_m\|_\infty\le\epsilon$, the secant slopes $s_k$ of
$g$ between successive extrema alternate in sign and have magnitude
at least $3/(2\pi)$. Hence $\sum_k|s_{k+1}-s_k|\gtrsim m$.

To bound the same sum from above, write each network coordinate as
a difference of two convex functions, starting with $x=x-0$.
A hidden affine row preserves a common Lipschitz bound on the two
terms by exchanging them for negative coefficients. Sorting a pair
at most doubles this bound, since
$\max(a-b,c-d)=\max(a+d,c+b)-(b+d)$ and minima follow by negation.
Thus $g$ is a difference of two convex $2^LB$-Lipschitz functions.
The secant slopes of a convex $K$-Lipschitz function are nondecreasing
and lie in $[-K,K]$, so their total variation is at most $2K$.
Applying this to both terms gives
\[
    \epsilon^{-1}\asymp m
    \lesssim \sum_k|s_{k+1}-s_k|
    \le 4\,2^L B,
\]
as claimed.
\end{example}

\par

\textit{Proof sketch of Theorem \ref{th:holder-approximation}.} A common approach to approximating smooth functions by ReLU networks
is to use local Taylor polynomials centered at the points of a regular grid that decomposes the input domain
\citep{SchmidtHieber2020,KohlerLanger2021,JiaoWangYang2023,NaglerLanger2026}.
Each polynomial approximates the target near its expansion point.
To obtain a uniform approximation, the network must evaluate these
polynomials and localize their contributions. In our GroupSort2
construction, both steps must respect the layerwise matrix norm bounds.

The approximation of polynomials with ReLU networks is grounded in an approximation result of \citet[Propositions~2--3]{Yarotsky2017} for the square function. Here $t^2$ is approximated on $[0,1]$ by its piecewise affine interpolant on a dyadic grid. The construction
iterates a hat map with slopes $\pm2$. In our class, each hidden
weight matrix has absolute row sums at most one, such that each hidden layer
is nonexpansive in the maximum norm. Lemma~\ref{lemma:squared}
shows that the same interpolant can be derived under this constraint by taking minima
of affine functions with slopes $\pm1$ and rescaling the partial sum
at each refinement. For $u,v\in[-1,1]$, the identity
\[
 uv=\left(\frac{|u+v|}{2}\right)^2
    -\left(\frac{|u-v|}{2}\right)^2
\]
then reduces multiplication to two evaluations of the squaring
network on $[0,1]$ (Lemma~\ref{lemma:multiplication}). Repeated
multiplication gives approximations of the monomials
(Lemma~\ref{lemma:monomial}), which we combine with the Taylor
coefficients to approximate each local polynomial.

For localization, let $h=1/K$ and consider the regular grid
$\mathcal G_h=\{0,h,\ldots,1\}^q$. To each grid point $a$, we
associate the tensor-product hat function
\[
 \eta_{a,h}(x)
 =\prod_{i=1}^q\left(1-\frac{|x_i-a_i|}{h}\right)_+.
\]
These functions are nonnegative, sum to one on $[0,1]^q$, and
vanish whenever $\|x-a\|_\infty\ge h$. Consequently, the sum
\[
 P_h(x)=\sum_{a\in\mathcal G_h}T_a(x)\eta_{a,h}(x)
\]
combines the local Taylor polynomials into a uniform approximation.
At each $x$, only polynomials centered within distance $h$
contribute, so local Taylor errors of order $h^\alpha$ give
$\|P_h-f\|_\infty\lesssim h^\alpha$.

To implement this construction, we need products of monomials
and hat functions. Writing
$T_a(x)=\sum_\beta c_{a,\beta}(x-a)^\beta$, we obtain
\[
 T_a(x)\eta_{a,h}(x)
 =\sum_\beta c_{a,\beta}h^{|\beta|}
   \psi_\beta\!\left(\frac{x-a}{h}\right),
 \qquad
 \psi_\beta(z)=z^\beta\prod_{i=1}^q(1-|z_i|)_+.
\]
However, directly implementing the change of variables
$z=(x-a)/h$ would require input weights of size $h^{-1}$,
violating the hidden matrix norm bounds. We therefore approximate
the scaled localized monomials
\[
 x\longmapsto \frac{h}{C_0}\,
 \psi_\beta\!\left(\frac{x-a}{h}\right),
\]
where $C_0$ is independent of $h$. The prefactor $h$ allows us
to use the same weight matrices at every scale. Indeed, if
$\Phi_\beta$ approximates $\psi_\beta/C_0$ at unit scale,
positive homogeneity of GroupSort implies that
$h\Phi_\beta((x-a)/h)$ can be realized by replacing the input
$x$ with $x-a$ and multiplying all biases by $h$.
The hat functions are represented exactly, and the multiplication
network returns zero whenever either input is zero. Thus the
localized approximant also vanishes whenever
$\|x-a\|_\infty\ge h$
(Lemma~\ref{lem:localized-taylor-sums}). To control the output-row norm, we combine approximations of Taylor polynomials across different grid scales. Set \(h_j=2^{-j}\), with \(h_J^\alpha\asymp\epsilon\), and write
\[
P_{h_J}=P_1+\sum_{j=1}^{J}\bigl(P_{h_j}-P_{h_{j-1}}\bigr).
\]
Comparing Taylor polynomials on consecutive grids shows that each difference has coefficients of order \(h_j^\alpha\) in the localized monomial basis above. Since our localized building blocks carry a prefactor \(h_j\), the construction uses normalization factors \(A_j\gtrsim h_j^{\alpha-1}\). For \(\alpha>1\), the sum
$\sum_{j=0}^{J}h_j^{\alpha-1}$ 
is bounded independently of \(J\). We choose the normalization factors so that \(\sum_{j=0}^{J}A_j\lesssim1\), while satisfying the stated width and depth bounds. The normalized subnetworks are combined through successive weighted averages, whose coefficients are nonnegative and sum to one. These averages preserve the hidden-layer row-norm constraints. The final affine map restores the overall scale, giving an output-row norm bounded independently of \(\epsilon\).

At \(\alpha=1\), this multiscale construction gives the larger output-row bound in the theorem. The separate fixed-width construction achieves a bounded output-row norm.
A complete proof of Theorem~\ref{th:holder-approximation}
is given in Section~\ref{sec:holder-proof}.

\subsection{Covering numbers and optimality}\label{sec:entropy}
Covering numbers serve two purposes in our analysis. They control the complexity entering the statistical
risk bounds and, through comparison with the metric entropy of the H\"older ball, prove optimality of our approximation rates in Theorem \ref{th:holder-approximation}. 

 For a class
$\mathcal G$ of functions on $\Omega$, the covering number
$\mathcal N(\epsilon,\mathcal G,\|\cdot\|_{\infty,\Omega})$ denotes
the smallest $N$ for which there exist $g_1,\ldots,g_N\in\mathcal G$
such that
\[
    \sup_{g\in\mathcal G}\min_{1\le j\le N}
    \|g-g_j\|_{\infty,\Omega}\le\epsilon.
\]
Fix $q\ge1$, $\alpha>0$, and $H\ge1$. On $D=[0,1]^q$, the
H\"older ball $\mathcal H=\{f\in C^\alpha(D):\|f\|_{C^\alpha(D)}\le H\}$
has metric entropy
\[
    \log\mathcal N(\epsilon,\mathcal H,\|\cdot\|_{\infty,D})
    \asymp \epsilon^{-q/\alpha},
    \qquad \epsilon\downarrow0,
\]
see \citet[Theorem~2.7.1]{VanDerVaartWellner1996} for the upper
bound and \citet{KolmogorovTikhomirov1959} for the lower bound.

If a class $\mathcal F$ approximates $\mathcal H$ to uniform error
$O(\epsilon)$, enlarging the balls of an $\epsilon$-cover of
$\mathcal F$ by a constant factor covers $\mathcal H$. Each enlarged
ball that meets $\mathcal H$ can be recentered at a point of
$\mathcal H$ by doubling its radius. Consequently,
\begin{equation}\label{eq:holder-entropy-transfer}
    \sup_{f\in\mathcal H}\inf_{g\in\mathcal F}
    \|f-g\|_{\infty,D}\lesssim\epsilon
    \quad\Longrightarrow\quad
    \log\mathcal N(\epsilon,\mathcal F,\|\cdot\|_{\infty,D})
    \gtrsim\epsilon^{-q/\alpha}.
\end{equation}

\begin{lemma}\label{lem:network-entropy}
Fix an integer $q\ge1$. Let $L\ge1$ be an integer, let $w\ge2$
be even, and let $b,B,R\ge1$. For every $H>0$ and $0<\epsilon<1$,
\begin{equation}\label{eq:network-entropy}
\begin{aligned}
    \log\mathcal N\!\left(
        \epsilon,\mathcal F(L,w,b,B,H),
        \|\cdot\|_{\infty,[-R,R]^q}\right)
    &\le \log\mathcal N\!\left(
        \epsilon,\widetilde{\mathcal F}(L,w,b,B),
        \|\cdot\|_{\infty,[-R,R]^q}\right)\\
    &\lesssim w^2L
        \log\!\left(\frac{e(q+w)(L+1)bBR}{\epsilon}\right).
\end{aligned}
\end{equation}
The implicit constant may depend on $q$. In particular, for fixed
$\alpha>0$, if $w^2L\lesssim\epsilon^{-q/\alpha}$ and $b,B,R$
are bounded by a fixed power of $\epsilon^{-1}$, then
\[
    \log\mathcal N\!\left(
        \epsilon,\mathcal F(L,w,b,B,H),
        \|\cdot\|_{\infty,[-R,R]^q}\right)
    \lesssim \epsilon^{-q/\alpha}\log(e/\epsilon).
\]
\end{lemma}

Clipping the centers of a cover gives a cover of the clipped class
at the same radius. For the second inequality in \eqref{eq:network-entropy},
we discretize the weights and biases as in the proof of
\citet[Lemma~5]{SchmidtHieber2020}. The unit row-norm bounds
ensure that each hidden layer propagates errors from the preceding
layers without amplifying them.

For $L\ge2$, write $S \asymp w^2L$ for the number of
weights and biases in the fully connected architecture with $L$
hidden layers of width $w$. If the network class with parameter bounds \eqref{eq:holder-network-budgets} approximates $\mathcal H$ to
uniform error $O(\epsilon)$, then \eqref{eq:holder-entropy-transfer}
and Lemma~\ref{lem:network-entropy} with $R=1$ give
\[
    \epsilon^{-q/\alpha}
    \lesssim S\log\!\left(
        \frac{e(q+w)(L+1)b_{\epsilon,w}B_{\epsilon,w}}{\epsilon}
    \right).
\]
When $S\le\epsilon^{-q/\alpha}$, both $w$ and $L$ are bounded by
fixed powers of $\epsilon^{-1}$. The same holds for the norm budgets
for every $\alpha>0$. In particular, at $\alpha=1$ we have
$w^{1-1/q}\lesssim\epsilon^{-(q-1)/2}$. Hence
\begin{equation}\label{eq:dense-necessary}
    S\log(e/\epsilon)\gtrsim\epsilon^{-q/\alpha}.
\end{equation}
For $S>\epsilon^{-q/\alpha}$, this inequality is immediate.
Equivalently, any such achievable accuracy satisfies
$\epsilon\gtrsim [S\log(eS)]^{-\alpha/q}$.

Theorem~\ref{th:holder-approximation} achieves uniform
error $O(\epsilon)$ with
$ S\lesssim\epsilon^{-q/\alpha}$ 
at fixed width for every $\alpha>0$, and throughout
its stated joint width--depth range for $\alpha\ge1$.
Conversely, \eqref{eq:dense-necessary} gives
\[
 S\gtrsim
 \frac{\epsilon^{-q/\alpha}}{\log(e/\epsilon)}.
\]
Thus, under the layerwise constraints and parameter
budgets considered above, the approximation rate in Theorem~\ref{th:holder-approximation} is optimal, i.e., the same rate cannot be achieved with a smaller number of network parameters (up to logarithmic factors).

At $\alpha=1$, the joint width--depth construction uses
\[
 B_{\epsilon,w}\asymp
 \log(e/\epsilon)+w^{1-1/q},
\]
while the fixed-width construction permits bounded $B$.
For $q\ge3$ and $L\asymp[\log(ew)]^2$,
Proposition~\ref{prop:bounded-endpoint} shows that
attaining the parameter rate $S^{-1/q}$ uniformly
over the $C^1$ ball requires $B$ to grow with $w$.
The proposition does not establish that the particular
growth prescribed in \eqref{eq:holder-network-budgets}
is necessary.

\section{Flow matching}\label{sec:main-results}\label{sec:flow-method}
In this section, we apply our class of GroupSort2 networks to flow matching. 

Let $Y_1,\ldots,Y_n$ be independent observations from an unknown
distribution $P$ on $\R^q$ with finite second moment. The goal of
flow matching is to estimate a velocity field whose flow transforms
a known (e.g.~$Z\sim N(0,I_q)$) random vector into a random vector with law
close to $P$.  
Let $Y\sim P$ be independent of $Z$. Following \citet{LiuGongLiu2023}, we consider the linear interpolation
\begin{equation}\label{eq:interpolation-definition}
  U_t=(1-t)Z+tY,\qquad P_t=\operatorname{Law}(U_t).  
\end{equation}
Here $U_0 \sim \mathcal{N}(0,I_q)$ and $U_1 \sim P$. Naturally, this specific endpoint coupling is not unique. For the target class introduced below, the marginal laws $P_t$ in \eqref{eq:interpolation-definition} are also obtained by starting from $\bar U_0\sim\mathcal N(0,I_q)$ and solving the ODE (see Lemma~\ref{lem:regularity})
\begin{equation}\label{eq:vectorfield}
    \dot{\bar U}_t= v^P_{t}(\bar U_t),\qquad 0\le t\le1, \qquad \text{where} \quad v_t^P(x)=\E[Y-Z\mid U_t=x].
\end{equation}

Since $P$ is unknown, the population velocity $v^P$ cannot be computed directly. Instead, we want to fit an estimator for $v^P$ chosen from a class of GroupSort2 networks. To construct a suitable objective, we first consider the population loss
\begin{equation*}
     \mathcal L(g)=\int_0^1\E|g_t(U_t)-(Y-Z)|^2\dd t.
\end{equation*}
Analogously to \cite[Theorem 2]{LipmanEtAl2023}, it can be shown that $\mathcal L(g)-\mathcal L(v^P)=\int_0^1\|g_t-v_t^P\|_{L^2(P_t)}^2\dd t$ whenever $\mathcal L(g)<\infty$. Hence optimizing $\mathcal L$ with respect to $g$ over a set $\mathcal{G}$ is equivalent to finding the best approximation for $v^P$ within $\mathcal{G}$ in the $L^2(P_t)$ sense. 
Since $P$ is unknown, we cannot compute $\mathcal L$ and instead consider the empirical counterpart based on the observations $Y_1,...,Y_n \sim P$
\begin{equation}\label{eq:loss}
 \widehat{\mathcal L}_n(g)=\frac{1}{n}\sum_{i=1}^n\int_0^1
 \E_Z\big|g_t((1-t)Z+tY_i)-(Y_i-Z)\big|^2\dd t.
\end{equation}
The expectation is over a standard Gaussian vector independent of the observations. In practice, this integration can be approximated up to arbitrary precision.

We are going to consider a specific set of GroupSort2 networks $\mathcal G_n$, which will be defined below and whose size will depend on the number of observations $n$. The vector field estimator is then given by any measurable minimizer
\begin{equation}\label{eq:erm}
 \widehat v_n\in\operatorname*{arg\,min}_{g\in\mathcal G_n}\widehat{\mathcal L}_n(g).
\end{equation}
Starting from an independent $X_0\sim N(0,I_q)$, we ultimately solve the ODE
\begin{equation}\label{eq:ode}
 \dot X_t=\widehat v_{n,t}(X_t),\qquad 0\le t\le1 .
\end{equation}
The final distributional estimator is $\widehat P_n=\operatorname{Law}(X_1\mid Y_1,\ldots,Y_n)$. 

In order to evaluate the quality of the estimator $\widehat P_n$ we need to specify a metric. We choose the $2$-Wasserstein distance, defined by
\begin{equation}
    W_2(P,Q)\coloneqq \Big(\inf_{\pi\in\Pi(P,Q)}\int|x-y|^2\,\pi(\dd x,\dd y)\Big)^{\frac{1}{2}},
\end{equation}
where $\Pi(P,Q)$ denotes the couplings of two probability measures $P$ and $Q$. Our goal is to bound $\mathbb{E}_P [W_2(\widehat P_n, P)]$, where we assume that $P$ is of the following form: Fix $\alpha>1$, $0<m<1<M$ and
$H>1$. Let $\varphi_q$ denote the standard Gaussian density, and
let $\mathcal P_\alpha$ be the class of probability measures with
densities
\begin{equation}\label{eq:model}
 p(x)=\frac{f(x)\varphi_q(x)}
 {\int_{\R^q}f(y)\varphi_q(y)\dd y},
 \qquad
 m\le f\le M,\qquad
 \|f\|_{C^\alpha(\R^q)}\le H.
\end{equation}
The bounds on $f$ control the density ratio, while $\alpha$ and
$H$ control its smoothness. For every $P\in\mathcal P_\alpha$,
the population velocity $v^P$ has a continuous version whose flow
transports $N(0,I_q)$ to $P$ and whose spatial Lipschitz constant is bounded (see Lemma~\ref{lem:regularity}).

\subsection{The network class}
In the following, we construct the network class $\mathcal G_n$. Depending on the time $t$, the vector field $v^P$ has different smoothness properties: for smaller $t$, the governing distribution is the Gaussian, which results in a very smooth $v^P$. For larger $t$, the smoothness rather depends on $P$, resulting in a more delicate approximation. Thus, we construct the velocity field separately on a partition of $[0,1]$, using smaller intervals near $t=1$.

Fix an integer $J\ge2$ for the number of intervals, a radius $R\ge1$ determining the region of approximation and a parameter bound
$b_n\ge1$ for the network weights and biases. Set $t_j=1-2^{-j}$ for $0\le j\le J$, put $t_{J+1}=1$, and define
\begin{equation}\label{eq:time-partition}
 I_j\coloneqq [t_{j-1},t_j],\qquad
 \lambda_j \coloneqq 2^{1-j},\qquad 1\le j\le J+1.
\end{equation}
For $j\le J$, we have
$\frac{\lambda_j}{2}\le1-t\le\lambda_j$ on $I_j$, so
$\lambda_j$ describes the distance of the interval from $t=1$.
The last interval is $I_{J+1}=[1-2^{-J},1]$ and has length
$\lambda_{J+1}=2^{-J}$.

On each $I_j$, we use one scalar network for each of the $q$
velocity coordinates. Choose an even width $w_j$, a depth $L_j\ge2$
and an output matrix bound $B_j\ge1$, and define the network construction $ g_{\theta_j} \colon \mathbb{R}^{q+1} \rightarrow \mathbb{R}^q$ with 
\begin{equation}\label{eq:block-parameter-set}
 g_{\theta_j}\coloneqq (g_{\theta_{j,1}},\ldots,g_{\theta_{j,q}})
 \in[\widetilde{\mathcal F}(L_j,w_j,b_n,B_j)]^q.
\end{equation}
 Since the hidden matrices have operator norm at most one,
$B_j$ bounds the Lipschitz constant in the maximum norm.
Each scalar network takes time and position as inputs, so its input
dimension is $q+1$. The architecture of each scalar network is
$k_j=(q+1,w_j,\ldots,w_j,1)$, with $L_j$ hidden layers, and its
number of weights and biases is
\begin{equation}\label{eq:flow-scalar-count}
 s_j=S(L_j, k_j)
 =(q+2)w_j+(L_j-1)w_j(w_j+1)+w_j+1\approx L_j w_j^2.
\end{equation}
To restrict each network to its interval while preserving continuity
in time, define the continuous cutoff function
\[
 \chi_j(t)=\min\left\{
 1,n^8\operatorname{dist}(t,[0,1]\setminus I_j)\right\},
\]
where $\operatorname{dist}(x,A)\coloneqq\inf_{y\in A}|x-y|$ for a nonempty set $A\subset\mathbb R$ and $x\in\mathbb R$.
This function vanishes outside $I_j$ and at its interior partition
endpoints. It equals one whenever the distance to the complement
is at least $n^{-8}$. 
For $z\in\R^q$, set
\begin{equation}\label{eq:gate}
 \Gamma_j(t,z)_i
 =\max\{-2R\chi_j(t),\min\{z_i,2R\chi_j(t)\}\},
 \qquad 1\le i\le q.
\end{equation}
This function clips each coordinate to
$[-2R\chi_j(t),2R\chi_j(t)]$.
We also clip the spatial input to $[-R,R]^q$, so that the network
is evaluated within the region used for approximation. Multiplying
time by $1+R$ puts its range on the same scale as the spatial
coordinates. We then define $\mathcal G_n$ as the class of fields
\begin{equation}\label{eq:flow-network-class}
 g_t(x)=\sum_{j=1}^{J+1}
 \Gamma_j\!\left(t,
 g_{\theta_j}\big((1+R)t,\clip_{[-R,R]^q}x\big)\right),
 \qquad
 g_{\theta_j}\text{ as in \eqref{eq:block-parameter-set}}.
\end{equation}
At each time, at most one summand in
\eqref{eq:flow-network-class} is nonzero. Consequently, every
coordinate of $g$ is bounded by $2R$. Clipping is nonexpansive,
so $g$ is globally Lipschitz in space, uniformly in time, with
Lipschitz constant at most $\max_j B_j$ in the maximum norm.
The construction via \eqref{eq:gate} also makes $g$ continuous in time. Hence the global version of the Picard–Lindelöf theorem guarantees that
\eqref{eq:ode} has a unique solution on $[0,1]$.

Finally, the matrix and parameter bounds define a compact parameter
set. The empirical loss is measurable in the observations and
continuous in the parameters. To see the continuity, note that dominated convergence
applies in \eqref{eq:loss}, since the velocity fields are uniformly
bounded and $\E|Z|^2<\infty$. The measurable maximum theorem
therefore gives a measurable minimizer in \eqref{eq:erm}
\citep[Theorem~18.19]{AliprantisBorder2006}.

\subsection{Minimax optimal rate of convergence}
In order to obtain the $2$-Wasserstein rate of convergence, we need to choose the number of intervals and the
spatial clipping radius: Set
\begin{equation}\label{eq:scales}
 J=\max\left\{2,\left\lceil
 \frac{\log_2 n}{2\alpha+\max\{q,2\}}
 \right\rceil\right\},\qquad
 R=\frac{2M}{m}+6\sqrt{\log(en)}.
\end{equation}
With this choice of $J$, the last interval has length
$h = 2^{-J}\asymp n^{-\frac{1}{2\alpha+\max\{q,2\}}}$.
The choice of $R$ gives $e^{-\frac{R^2}{4}}\le n^{-8}$,
which controls the error in the area where the fitted vector field is clipped instead of being an approximation of $v^P$.

Additionally, we need to fix the number of layers $L_j$, widths $w_j$, the parameter bounds $b_n$ and the Lipschitz bounds $B_j$.
We take sufficiently large even widths satisfying
\begin{equation}\label{eq:flow-joint-range}
 w_j^2\lesssim\lambda_j^{-q},\qquad
 w_j^{1-1/q}\lesssim\lambda_j^{1-\alpha},
 \qquad 1\le j\le J+1
\end{equation}
and choose
\begin{equation}\label{eq:flow-architecture-budget}
 L_j=\left\lceil
 C_L(1+R)^{q+1}[\log(en)]^{q+3}
 \frac{\lambda_j^{-q}}{w_j^2}
 \right\rceil,
\end{equation}
with $C_L$ sufficiently large.
Furthermore, we fix
\begin{equation}\label{eq:block-domains}
 b_n=C_b n,\qquad B_j=C_B,\qquad 1\le j\le J+1,
\end{equation}
with $C_b$ and $C_B$ sufficiently large.
Combining the $q$ coordinates,
every field in the class has spatial Lipschitz constant
at most $\sqrt q\,C_B$, uniformly in time. 

With these choices, we obtain the following result. 

\begin{theorem}\label{thm:main-flow}
For $\alpha>1$, with the choices \eqref{eq:scales}--\eqref{eq:flow-architecture-budget},
the estimator defined by \eqref{eq:erm} and \eqref{eq:ode} satisfies
\begin{equation}\label{eq:flow-fixed-orders}
 \sup_{P\in\mathcal P_\alpha}
 \E_P W_2(\widehat P_n,P)
 \lesssim [\log(en)]^{\frac{3q+11}{4}}
 \begin{cases}
 n^{-\frac12},&q=1,\\
 n^{-\frac12}\log(en),&q=2,\\
 n^{-\frac{\alpha+1}{2\alpha+q}},&q>2.
 \end{cases}
\end{equation}
Conversely, taking the infimum over all estimators
based on $Y_1,\ldots,Y_n$ gives
\begin{equation}\label{eq:lower}
 \inf_{\widetilde P_n}\sup_{P\in\mathcal P_\alpha}
 \E_P W_1(\widetilde P_n,P)
 \gtrsim
 \begin{cases}
 n^{-\frac12},&q\le2,\\
 n^{-\frac{\alpha+1}{2\alpha+q}},&q>2.
 \end{cases}
\end{equation}
\end{theorem}

Since $W_1\le W_2$, \eqref{eq:lower} also bounds the minimax
$W_2$ risk from below. Hence \eqref{eq:flow-fixed-orders} attains
the minimax rate up to logarithmic factors. The lower bound
follows by transferring the construction of
\citet[Theorem~3]{NilesWeedBerthet2022} from the cube to
\eqref{eq:model} (see Appendix~\ref{sec:lower}).

\begin{remark}[Choice of width and depth]\label{rem:flow-width-depth} 
We can use the same fixed width on every interval.
This choice satisfies \eqref{eq:flow-joint-range} for
every $\alpha>1$, but requires more layers on smaller
intervals. Increasing the width on smaller intervals allows us
to reduce the depth. When $\alpha\ge(q+1)/2$, we can
take $w_j\asymp\lambda_j^{-q/2}$. The ratio
$\lambda_j^{-q}/w_j^2$ is then of constant order.
Since $1+R\asymp\sqrt{\log(en)}$,
\eqref{eq:flow-architecture-budget} gives
$L_j\asymp[\log(en)]^{\frac{3q+7}{2}}$
for every $1\le j\le J+1$.
When $1<\alpha<(q+1)/2$, the second condition in
\eqref{eq:flow-joint-range} limits how much the width
can grow. In this case, we can take
$w_j\asymp\lambda_j^{-q(\alpha-1)/(q-1)}$. For all these choices, the total number of fitted
parameters has the same order.
Using $\sum_j\lambda_j^{-q}\asymp h^{-q}$ and
the choice of $h$ in \eqref{eq:scales}, we obtain
\[
 q\sum_{j=1}^{J+1}s_j
 \asymp [\log(en)]^{\frac{3q+7}{2}}h^{-q}
 \asymp
 n^{\frac{q}{2\alpha+\max\{q,2\}}}
 [\log(en)]^{\frac{3q+7}{2}}.
\]
\end{remark}

\section{Nonparametric regression}\label{sec:regression-method}
Let $(X_i,Y_i)_{i=1}^n$ be independent copies of $(X,Y)$, where $X$ has distribution $P_X$ on $[0,1]^q$. Write $P_{X,Y}$ for the joint law. Consider the regression model
\[
 Y_i=f_0(X_i)+\eta_i,\qquad \E[\eta_i\mid X_i]=0,\qquad i=1,\ldots,n .
\]
Fix $\alpha>0$ and $M_Y\ge H\ge1$. Assume that $|Y|\le M_Y$ almost surely, that $\operatorname{Var}(\eta_i)>0$, and that the regression function fulfills $\|f_0\|_{C^\alpha([0,1]^q)}\le H$. Our estimator is any measurable empirical risk minimizer
\begin{equation}\label{eq:regression-erm}
 \widehat f_n\in\operatorname*{arg\,min}_{f\in\mathcal F_n}
           \frac{1}{n}\sum_{i=1}^n(Y_i-f(X_i))^2 ,
\end{equation}
where 
\begin{equation}\label{eq:regression-resolution}
 \epsilon_n=\min\left\{\frac{1}{8},
    \left(\frac{\log(en)}{n}\right)^{\frac{\alpha}{2\alpha+q}}\right\},
 \qquad
 \mathcal F_n=\mathcal F(L_n,w_n,b_{\epsilon_n,w_n},B_{\epsilon_n,w_n},H).
\end{equation}
The bounds \(b_{\epsilon_n,w_n}\) and \(B_{\epsilon_n,w_n}\) are chosen as in Theorem \ref{th:holder-approximation}. If \(\alpha=1\) and \(L_n\gtrsim\epsilon_n^{-q}\), we may alternatively use the constant bounds
\[
B_{\epsilon_n,w_n}=C_B,\qquad
b_{\epsilon_n,w_n}=C_bC_B,
\]
as permitted by the final assertion of that theorem.

The following theorem matches the classical minimax convergence rate for Hölder classes \citep{Stone1982}, up to a logarithmic factor. Its upper bound holds for every covariate distribution \(P_X\) supported on \([0,1]^q\).

\begin{theorem}\label{thm:regression}
Assume \eqref{eq:regression-resolution}  and let \((L_n,w_n)\) satisfy the architecture conditions of Theorem~\ref{th:holder-approximation} at \(\epsilon=\epsilon_n\), and use the corresponding parameter bounds specified above and in Theorem~\ref{th:holder-approximation}. If
\begin{equation}\label{eq:regression-joint-budget}
 w_n^2(L_n+1)\lesssim \epsilon_n^{-\frac{q}{\alpha}},
\end{equation}
then the estimator \eqref{eq:regression-erm} over the class in \eqref{eq:regression-resolution} satisfies
\begin{equation}\label{eq:regression-risk}
 \E\|\widehat f_n-f_0\|_{L^2(P_X)}^2
 \lesssim \left(\frac{\log(en)}{n}\right)^{\frac{2\alpha}{2\alpha+q}} .
\end{equation}
\end{theorem}
The architectural conditions of Theorem \ref{thm:regression} are satisfied by taking a sufficiently large fixed even width \(w_n=w_*\) and
\[
L_n=\left\lceil C\epsilon_n^{-q/\alpha}\right\rceil,
\]
with \(C\) sufficiently large and independent of \(n\). For \(\alpha>1\), the resulting network classes have a common Lipschitz bound independent of \(n\). The same holds at \(\alpha=1\) when using the constant parameter bounds specified above. For \(0<\alpha<1\), the output-row bound grows as \(\epsilon_n\) decreases.

The estimator in Theorem \ref{thm:regression} also applies to sub-Gaussian responses. For this, retain the assumptions on $P_X$ and $f_0$ and suppose, in place of boundedness, that for some $\sigma>0$,
\begin{equation}\label{eq:subgaussian-model}
 \E\big[\exp\{t(Y-f_0(X))\}\mid X\big]\le
 \exp(\frac{\sigma^2t^2}{2}),\qquad t\in\R .
\end{equation}

\begin{corollary}\label{cor:regression-subgaussian}
Under \eqref{eq:subgaussian-model}, set
\begin{equation}\label{eq:regression-subgaussian-resolution}
 \epsilon_n=\min\left\{\frac{1}{8},
 \left(\frac{[\log(en)]^2}{n}\right)^{
       \frac{\alpha}{2\alpha+q}}\right\},
 \qquad
 \mathcal F_n=
 \mathcal F(L_n,w_n,b_{\epsilon_n,w_n},B_{\epsilon_n,w_n},H).
\end{equation}
Let $(w_n,L_n)$ satisfy the hypotheses of
Theorem~\ref{th:holder-approximation} at
$\epsilon=\epsilon_n$ and
$w_n^2(L_n+1)\lesssim\epsilon_n^{-\frac{q}{\alpha}}$.
Then the estimator \eqref{eq:regression-erm} satisfies
\begin{equation}\label{eq:regression-subgaussian-rate}
 \E\|\widehat f_n-f_0\|_{L^2(P_X)}^2
 \lesssim
 \left(\frac{[\log(en)]^2}{n}\right)^{
       \frac{2\alpha}{2\alpha+q}}.
\end{equation}
\end{corollary}

\begin{remark}[Extension to compositional regression functions] Theorem~\ref{thm:regression} and Corollary \ref{cor:regression-subgaussian} suffer from the well-known curse of dimensionality: increasing input dimension $q$ leads to increasingly slow convergence rates. To circumvent this, several results, e.g., \cite{SchmidtHieber2020, KohlerLanger2021, NaglerLanger2026} considered regression models with compositional structure, where
$f_0=h_m\circ\cdots\circ h_0$ and each coordinate of
$h_i$ depends on at most $t_i$ variables and belongs to
a bounded H\"older class of smoothness $\beta_i\ge1$.

Extending our results to these functions requires applying Theorem~\ref{th:holder-approximation}
to each coordinate in each component $h_j$ and clipping intermediate outputs to their
prescribed domains. The additional step is to combine these
networks while preserving the layerwise norm constraints.
Positive homogeneity of GroupSort allows intermediate
scaling factors to be transferred to the final output
layer and in turn, for a fixed composition, the resulting norm budget
can remain independent of the sample size.

With this construction,
a corresponding covering number bound together with the regression arguments would yield,
for suitably chosen network classes,
\[
 \E\|\widehat f_n-f_0\|_{L^2(P_X)}^2
 \lesssim
 \max_{0\le i\le m}
 \left(\frac{\log(en)}{n}\right)^{
       \frac{2\beta_i}{2\beta_i+t_i}}
\]
for bounded responses, with $\log(en)$ replaced by
$[\log(en)]^2$ for conditionally sub-Gaussian errors.
We leave the details of this extension to future work.

\end{remark}

\section{Conclusion}\label{sec:conclusion}

We establish approximation and covering-number bounds for
GroupSort2 networks under  layerwise norm constraints.
For H\"older targets with smoothness $\alpha>1$, the constructions
exploit target regularity while keeping the parameter and
Lipschitz bounds independent of the approximation accuracy.
The covering-number analysis quantifies the statistical complexity of our network class
and establishes optimality of the total parameter count up to
logarithmic factors under the prescribed constraints.

These results yield statistical guarantees for flow matching
and nonparametric regression. Under the Gaussian density-ratio
model, uniform spatial Lipschitz control over the fitted
velocity class keeps the flow stability factor bounded
independently of the sample size. The resulting estimator
attains the minimax Wasserstein-$2$ rate up to logarithmic
factors with polynomially many parameters. For regression, constrained GroupSort2
least-squares estimators attain the optimal minimax convergence rate for H\"older functions up to logarithmic factors.


Two directions are particularly interesting for future research. First, extending
the flow-matching analysis beyond the Gaussian density-ratio
model would clarify the scope of the statistical guarantees.
Second, practical implementations commonly introduce optimization error,
Monte Carlo approximation of the training loss, and numerical
error in solving the flow ODE. Quantifying their combined
effect would help determine the computational accuracy needed
to retain the statistical rates established here. The gap between statistical rates and optimization behavior of the networks remains an open gap in the broader network literature and would require fundamental new proof strategies.

\section*{Acknowledgment} S.L. and L.K.'s work was supported by the Deutsche Forschungsgemeinschaft (DFG) under project number
460867398 (DFG Research unit 5381 Mathematical Statistics in the Information Age).

\clearpage
\appendix

\section{Proof of Theorem~\ref{th:holder-approximation} for \texorpdfstring{$\alpha\geq 1$}{alpha >= 1}}\label{sec:holder-proof}
We first provide the overall proof of Theorem~\ref{th:holder-approximation} and then present the lemmas and their proofs used within the construction.
\begin{proof}[Proof of Theorem~\ref{th:holder-approximation}]
The case $0<\alpha<1$ is shown in Lemma~\ref{lem:fixed-width-cone}. Therefore, we assume smoothness $\alpha\ge1$ in the following. We define first a partition of unity
by cone-type functions on increasingly fine grids.
For each $j\ge0$, we fix a sequence of resolutions $h_j=2^{-j}$
and grids
\[
 \mathcal A_j=(h_j\mathbb Z^q)\cap[0,1]^q
 =
 \left\{\left(\frac{k_1}{2^j},\ldots,\frac{k_q}{2^j}\right):
 k_1,\ldots,k_q\in\{0,1,\ldots,2^j\}\right\}.
\]
We define the hat functions
\[
 \theta_{j,a}(x)=
 \prod_{i=1}^q\left(1-\frac{|x_i-a_i|}{h_j}\right)_+,
 \qquad
 Q_jf(x)=\sum_{a\in\mathcal A_j}\theta_{j,a}(x)P_a(x),
\]
where $P_a(x)=\sum_{|\beta|<\alpha}
\frac{D^\beta f(a)}{\beta!}(x-a)^\beta$
is the Taylor polynomial of $f$ centered at $a$,
of degree at most $\lceil\alpha\rceil-1$.
For $\alpha=1$, we have $P_a=f(a)$.
The functions $\theta_{j,a}$ sum to one and are nonzero
only when $\|x-a\|_\infty<h_j$. The Taylor remainder estimate \citep[Lemma~A.8]{PetersenVoigtlaender2018} then gives
\[
 \|f-Q_jf\|_\infty
 =\left\|\sum_{a\in\mathcal A_j}\theta_{j,a}(f-P_a)\right\|_\infty
 \le \sup_{x\in[0,1]^q}
 \sum_{a\in\mathcal A_j}\theta_{j,a}(x)|f(x)-P_a(x)|
 \lesssim h_j^\alpha.
\]

These partitions also allow us to compare two consecutive
resolutions. Since $h_{j-1}=2h_j$, for $b\in\mathcal A_{j-1}$
we get
\[
\begin{aligned}
 \theta_{j-1,b}(x)
 &=\prod_{i=1}^q
   \left(1-\frac{|x_i-b_i|}{2h_j}\right)_+\\
 &=\prod_{i=1}^q
   \left[
   \sum_{a_i\in(h_j\mathbb Z)\cap[0,1]}
   \left(1-\frac{|a_i-b_i|}{2h_j}\right)_+
   \left(1-\frac{|x_i-a_i|}{h_j}\right)_+
   \right]\\
 &=\sum_{a\in\mathcal A_j}
   \underbrace{\prod_{i=1}^q
   \left(1-\frac{|a_i-b_i|}{2h_j}\right)_+}_{\theta_{j-1,b}(a)}
   \underbrace{\prod_{i=1}^q
   \left(1-\frac{|x_i-a_i|}{h_j}\right)_+}_{\theta_{j,a}(x)}.
\end{aligned}
\]
Substituting the preceding identity into the definition
of $Q_{j-1}f$ and grouping the terms with the same
$\theta_{j,a}$, we obtain
\[
 Q_jf-Q_{j-1}f
 =\sum_{a\in\mathcal A_j}\theta_{j,a}
 \left(
 P_a-\sum_{b\in\mathcal A_{j-1}}
          \theta_{j-1,b}(a)P_b
 \right).
\]
The polynomial in parentheses compares $P_a$ with an average
of the polynomials $P_b$. The weights $\theta_{j-1,b}(a)$
are nonnegative, sum to one, and are nonzero only when
$\|a-b\|_\infty\le2h_j$.
To see why its coefficients are small, we write it
in powers of $\frac{x-a}{h_j}$. The coefficient corresponding
to a multi-index $\beta$ with $|\beta|<\alpha$ is
\[
 c_{j,a,\beta}
 =\frac{h_j^{|\beta|}}{\beta!}
 \sum_{b\in\mathcal A_{j-1}}\theta_{j-1,b}(a)
 \bigl(D^\beta f(a)-D^\beta P_b(a)\bigr).
\]
By Taylor's theorem, when the weight is nonzero,
\[
 |D^\beta f(a)-D^\beta P_b(a)|
 \lesssim h_j^{\alpha-|\beta|}.
\]
Since the weights sum to one, we obtain
$|c_{j,a,\beta}|
\lesssim h_j^{|\beta|}h_j^{\alpha-|\beta|}
=h_j^\alpha$.

We first approximate $Q_0f$ and then each difference
$Q_jf-Q_{j-1}f$. The sum of these approximations approximates
$Q_Jf$ through the telescoping sum
$Q_Jf=Q_0f+\sum_{j=1}^J(Q_jf-Q_{j-1}f)$.
The functions that we need to approximate are the products
of $\theta_{j,a}$ with the monomials centered at $a$:
\[
 \theta_{j,a}(x)
 \left(\frac{x-a}{h_j}\right)^\beta
 =
 \psi_\beta\!\left(\frac{x-a}{h_j}\right),
 \qquad
 \psi_\beta(z)=z^\beta\prod_{i=1}^q(1-|z_i|)_+,
 \qquad |\beta|<\alpha.
\]
Each one vanishes outside the cube $a+[-h_j,h_j]^q$.
With the coefficients obtained above,
\begin{align}
\label{eq:differenceQ}
Q_jf(x)-Q_{j-1}f(x)
 =\sum_{a\in\mathcal A_j}\sum_{|\beta|<\alpha}
 c_{j,a,\beta}\,
 \psi_\beta\!\left(\frac{x-a}{h_j}\right),
 \qquad |c_{j,a,\beta}|\lesssim h_j^\alpha.
\end{align}
For $j=0$, the same representation holds for $Q_0f$ with
$c_{0,a,\beta}=\frac{D^\beta f(a)}{\beta!}$, since $h_0=1$.

Lemma~\ref{lem:localized-taylor-sums} gives a  constant $K_0$ such that, for every error $e_j>0$, every sufficiently large
even width $w$, and every $A_j\gtrsim h_j^{\alpha-1}$,
there exists a network $F_j^\circ$ of width at most $w$ such that
\begin{align}
    \label{eq: taylornetbound}
    \left\|Q_jf-Q_{j-1}f-K_0A_jF_j^\circ\right\|_\infty
 \le e_j.
\end{align}
For $j=0$, replace $Q_jf-Q_{j-1}f$ by $Q_0f$.
All its matrices have absolute row sums at most one,
its biases are bounded, and its depth $L_j$ is at most
\begin{equation}\label{eq:cota-profund}
 C\left\{
 j+1+
 \log\!\left(2+\frac{h_j^\alpha}{e_j}\right)
 \left(
 1+\frac{h_j^{-q}}{w^2}
 +\frac{h_j^{\alpha-1-q}}{A_jw^{1+\frac{1}{q}}}
 \right)
 \right\}.
\end{equation}
We compute $F_0^\circ(x),\ldots,F_J^\circ(x)$ successively,
carrying the input $x$ and the partial weighted average
\[
 a_j(x)=
 \frac{\sum_{l=0}^j A_lF_l^\circ(x)}
      {\sum_{l=0}^j A_l}.
\]
Starting from $a_0=F_0^\circ$, each update is
\begin{equation}\label{eq:dense-level-addition}
 a_j=
 \frac{\sum_{l=0}^{j-1}A_l}{\sum_{l=0}^jA_l}\,a_{j-1}
 +\frac{A_j}{\sum_{l=0}^jA_l}\,F_j^\circ.
\end{equation}
The two coefficients are nonnegative and sum to one.
Each update and the copying of the stored coordinates require
only a fixed number of additional layers, so the total depth
is at most $\sum_{j=0}^J L_j+C(J+1)$.

We choose $J$ as the smallest integer satisfying
$Ch_J^\alpha\le\frac{\epsilon}{2}$, so that
$\|f-Q_Jf\|_\infty\le\frac{\epsilon}{2}$.
For each $0\le j\le J$, we take
$e_j=\frac{\epsilon}{4}\,2^{-(J-j)}$.
Thus, the total error satisfies
$\sum_{j=0}^J e_j\le\frac{\epsilon}{2}$.
Since $e_j$ bounds the error after multiplying the output
by $K_0A_j$, the triangle inequality gives
\[
 \left\|f-\sum_{j=0}^J K_0A_jF_j^\circ\right\|_\infty
 \le \|f-Q_Jf\|_\infty+\sum_{j=0}^J e_j
 \le\epsilon.
\]
This choice also allows us to control the depth:
since $h_j^\alpha\asymp\epsilon\,2^{\alpha(J-j)}$,
we have
$\log(2+\frac{h_j^\alpha}{e_j})\lesssim1+J-j$.

To choose the factors $A_j$, we take
\[
 A_j=
 C h_j^{\alpha-1}
 +
 \frac{\sqrt{(1+J-j)h_j^{\alpha-1-q}}}
 {\displaystyle\sum_{l=0}^J
  \sqrt{(1+J-l)h_l^{\alpha-1-q}}}
 \begin{cases}
 w^{1-\frac1q},&\alpha=1,\\
 1,&\alpha>1,
 \end{cases}
\]
with $C$ sufficiently large.
The first summand ensures the condition
$A_j\gtrsim h_j^{\alpha-1}$.
The fractions in the second summand are positive and sum
to one. Therefore,
\[
 \sum_{j=0}^J A_j
 \le
 \begin{cases}
 C(J+1)+w^{1-\frac1q},&\alpha=1,\\[3pt]
 \displaystyle\frac{C}{1-2^{-(\alpha-1)}}+1,&\alpha>1.
 \end{cases}
\]
This keeps the norm of the final output matrix bounded
when $\alpha>1$. For $\alpha=1$, it gives an output bound
of order $\log(\frac e\epsilon)+w^{1-\frac1q}$.

The second summand also allows us to control the number
of layers. Summing in \eqref{eq:cota-profund}, the part
of the total depth that depends on these factors satisfies
\[
 \frac{1}{w^{1+\frac1q}}
 \sum_{j=0}^J
 \frac{(1+J-j)h_j^{\alpha-1-q}}{A_j}
 \le
 \left(
 \sum_{j=0}^J
 \sqrt{(1+J-j)h_j^{\alpha-1-q}}
 \right)^2
 \begin{cases}
 \frac{1}{w^2},&\alpha=1,\\[3pt]
 \frac{1}{w^{1+\frac1q}},&\alpha>1.
 \end{cases}
\]
The sum on the right satisfies
\[
 \left(
 \sum_{j=0}^J
 \sqrt{(1+J-j)h_j^{\alpha-1-q}}
 \right)^2
 \asymp
 \begin{cases}
 \epsilon^{-\frac{q+1-\alpha}{\alpha}},
     &1\le\alpha<q+1,\\[2pt]
 [\log(\frac e\epsilon)]^3,
     &\alpha=q+1,\\[2pt]
 \log(\frac e\epsilon),
     &\alpha>q+1.
 \end{cases}
\]
Also,
$\sum_{j=0}^J(1+J-j)h_j^{-q}\lesssim h_J^{-q}
\asymp\epsilon^{-\frac q\alpha}$.
The total depth is therefore at most
\begin{equation}
\label{propfundidad-max}
     C\left\{
 \left[\log\!\left(\frac e\epsilon\right)\right]^2
 +\frac{\epsilon^{-\frac q\alpha}}{w^2}
 +
 \begin{cases}
 0,&\alpha=1,\\[2pt]
 \displaystyle
 \frac{\epsilon^{-\frac{q+1-\alpha}{\alpha}}}
      {w^{1+\frac1q}},
     &1<\alpha<q+1,\\[7pt]
 \displaystyle
 \frac{[\log(\frac e\epsilon)]^3}{w^{1+\frac1q}},
     &\alpha=q+1,\\[7pt]
 \displaystyle
 \frac{\log(\frac e\epsilon)}{w^{1+\frac1q}},
     &\alpha>q+1.
 \end{cases}
 \right\}.
\end{equation}
For $\alpha=1$, the contribution involving $A_j$ has been
included in $\frac{\epsilon^{-q}}{w^2}$.
When $\alpha=q+1$, use
\[
 \frac{(J+1)^3}{w^{1+\frac1q}}
 \lesssim (J+1)^2+\frac{h_J^{-q}}{w^2}.
\]
For $w\ge(J+1)^{\frac q{q+1}}$, the first term suffices.
For smaller $w$, the second suffices because $h_J^{-q}=2^{qJ}$.
When $\alpha>q+1$, the last depth term is absorbed by
$[\log(\frac e\epsilon)]^2$. 

The network constructed above outputs $a_J$. Multiplying its
final affine row by $K_0\sum_{j=0}^J A_j$ and clipping to
$[-H,H]$ gives
\[
 g(x)=\clip_{[-H,H]}\!\left(
 K_0\left(\sum_{j=0}^J A_j\right)a_J(x)\right).
\]
Since $\|f\|_\infty\le H$, clipping preserves the error bound
$\|f-g\|_{\infty,[0,1]^q}\le\epsilon$. The output row has
absolute sum at most $K_0\sum_{j=0}^J A_j\le B_{\epsilon,w}$
by \eqref{eq:holder-network-budgets}. Thus
$g\in\mathcal F(L,w,b_{\epsilon,w},B_{\epsilon,w},H)$ and by selecting the parameters in \eqref{propfundidad-max} the result follows. 
\end{proof}

\subsection{Approximation of local Taylor polynomials}
The major building block of Theorem \ref{th:holder-approximation} is a network construction of the local Taylor polynomials multiplied by the locating hat functions that enables to bound \eqref{eq: taylornetbound} in the proof of Theorem \ref{th:holder-approximation} above.
The graphical illustration of the proof for $\beta=0$ and dimension $1$ is in Figures~\ref{fig:a6-ejemplo-coordenadas} and \ref{fig:a6-ejemplo-producto}.

\begin{figure}[htbp]
\centering
\resizebox{1.0\linewidth}{!}{%

\tikzset{every picture/.style={line width=1.05pt}} 

\tikzset{every picture/.style={line width=0.75pt}} 

\begin{tikzpicture}[x=0.75pt,y=0.75pt,yscale=-1,xscale=1]

\draw  [color={rgb, 255:red, 221; green, 244; blue, 202 }  ,draw opacity=1 ][fill={rgb, 255:red, 221; green, 244; blue, 202 }  ,fill opacity=1 ] (91.44,39.67) -- (152,39.67) -- (152,69.33) -- (91.44,69.33) -- cycle ;
\draw  [color={rgb, 255:red, 221; green, 244; blue, 202 }  ,draw opacity=1 ][fill={rgb, 255:red, 221; green, 244; blue, 202 }  ,fill opacity=1 ] (329,39.67) -- (393.56,39.67) -- (393.56,69.33) -- (329,69.33) -- cycle ;
\draw  [draw opacity=0] (91.44,69.33) -- (439,69.33) -- (439,57.67) -- (91.44,57.67) -- cycle ; \draw   (91.44,69.33) -- (91.44,57.67)(111.44,69.33) -- (111.44,57.67)(131.44,69.33) -- (131.44,57.67)(151.44,69.33) -- (151.44,57.67)(171.44,69.33) -- (171.44,57.67)(191.44,69.33) -- (191.44,57.67)(211.44,69.33) -- (211.44,57.67)(231.44,69.33) -- (231.44,57.67)(251.44,69.33) -- (251.44,57.67)(271.44,69.33) -- (271.44,57.67)(291.44,69.33) -- (291.44,57.67)(311.44,69.33) -- (311.44,57.67)(331.44,69.33) -- (331.44,57.67)(351.44,69.33) -- (351.44,57.67)(371.44,69.33) -- (371.44,57.67)(391.44,69.33) -- (391.44,57.67)(411.44,69.33) -- (411.44,57.67)(431.44,69.33) -- (431.44,57.67) ; \draw   (91.44,69.33) -- (439,69.33) ; \draw    ;
\draw  [fill={rgb, 255:red, 74; green, 144; blue, 226 }  ,fill opacity=1 ] (338.11,68.33) .. controls (338.11,66.86) and (339.31,65.67) .. (340.78,65.67) .. controls (342.25,65.67) and (343.44,66.86) .. (343.44,68.33) .. controls (343.44,69.81) and (342.25,71) .. (340.78,71) .. controls (339.31,71) and (338.11,69.81) .. (338.11,68.33) -- cycle ;
\draw   (95,73.67) .. controls (95.09,78.34) and (97.46,80.62) .. (102.13,80.53) -- (111.14,80.36) .. controls (117.8,80.23) and (121.18,82.5) .. (121.27,87.16) .. controls (121.18,82.5) and (124.46,80.1) .. (131.13,79.97)(128.13,80.03) -- (140.14,79.8) .. controls (144.8,79.71) and (147.09,77.34) .. (147,72.67) ;
\draw  [color={rgb, 255:red, 221; green, 244; blue, 202 }  ,draw opacity=1 ][fill={rgb, 255:red, 221; green, 244; blue, 202 }  ,fill opacity=1 ] (616.5,50.22) -- (732.5,50.22) -- (732.5,131.56) -- (616.5,131.56) -- cycle ;
\draw  [draw opacity=0] (534.94,104.85) -- (820.94,104.95) -- (820.94,87.85) -- (534.95,87.75) -- cycle ; \draw   (534.94,104.85) -- (534.95,87.75)(574.94,104.86) -- (574.95,87.76)(614.94,104.87) -- (614.95,87.78)(654.94,104.89) -- (654.95,87.79)(694.94,104.9) -- (694.95,87.8)(734.94,104.92) -- (734.95,87.82)(774.94,104.93) -- (774.95,87.83)(814.94,104.94) -- (814.95,87.84) ; \draw   (534.94,104.85) -- (820.94,104.95) ; \draw    ;
\draw  [color={rgb, 255:red, 130; green, 6; blue, 185 }  ,draw opacity=1 ] (566,104.14) -- (595.74,54.17) -- (756.26,54.17) -- (786,104.14) -- cycle ;
\draw  [color={rgb, 255:red, 221; green, 244; blue, 202 }  ,draw opacity=1 ][fill={rgb, 255:red, 221; green, 244; blue, 202 }  ,fill opacity=1 ] (159,146.67) -- (219,146.67) -- (219,288.44) -- (159,288.44) -- cycle ;
\draw  [color={rgb, 255:red, 221; green, 244; blue, 202 }  ,draw opacity=1 ][fill={rgb, 255:red, 221; green, 244; blue, 202 }  ,fill opacity=1 ] (400,138.67) -- (462.69,138.67) -- (462.69,291.05) -- (400,291.05) -- cycle ;
\draw  [draw opacity=0] (159.46,246.11) -- (540.67,246.11) -- (540.67,240.33) -- (159.46,240.33) -- cycle ; \draw   (159.46,246.11) -- (159.46,240.33)(179.39,246.11) -- (179.39,240.33)(199.32,246.11) -- (199.32,240.33)(219.25,246.11) -- (219.25,240.33)(239.17,246.11) -- (239.17,240.33)(259.1,246.11) -- (259.1,240.33)(279.03,246.11) -- (279.03,240.33)(298.96,246.11) -- (298.96,240.33)(318.89,246.11) -- (318.89,240.33)(338.81,246.11) -- (338.81,240.33)(358.74,246.11) -- (358.74,240.33)(378.67,246.11) -- (378.67,240.33)(398.6,246.11) -- (398.6,240.33)(418.53,246.11) -- (418.53,240.33)(438.45,246.11) -- (438.45,240.33)(458.38,246.11) -- (458.38,240.33)(478.31,246.11) -- (478.31,240.33)(498.24,246.11) -- (498.24,240.33)(518.17,246.11) -- (518.17,240.33)(538.09,246.11) -- (538.09,240.33) ; \draw   (159.46,246.11) -- (540.67,246.11) ; \draw    ;
\draw  [fill={rgb, 255:red, 74; green, 144; blue, 226 }  ,fill opacity=1 ] (405.94,245.81) .. controls (405.94,244.48) and (407.13,243.41) .. (408.6,243.41) .. controls (410.07,243.41) and (411.26,244.48) .. (411.26,245.81) .. controls (411.26,247.13) and (410.07,248.2) .. (408.6,248.2) .. controls (407.13,248.2) and (405.94,247.13) .. (405.94,245.81) -- cycle ;
\draw  [dash pattern={on 0.84pt off 2.51pt}]  (98.54,195.94) -- (583.9,194.35) ;
\draw  [dash pattern={on 0.84pt off 2.51pt}]  (95.2,161.02) -- (577.24,161.82) ;
\draw  [dash pattern={on 0.84pt off 2.51pt}]  (101.55,275.66) -- (575.61,274.66) ;
\draw [color={rgb, 255:red, 74; green, 144; blue, 226 }  ,draw opacity=1 ]   (260,162.67) -- (123,276.67) ;
\draw [color={rgb, 255:red, 74; green, 144; blue, 226 }  ,draw opacity=1 ]   (504.67,162.33) -- (364.67,273.33) ;
\draw [color={rgb, 255:red, 74; green, 144; blue, 226 }  ,draw opacity=1 ]   (364.67,275.33) -- (260,162.67) ;
\draw  [fill={rgb, 255:red, 224; green, 16; blue, 57 }  ,fill opacity=1 ] (97.78,67.89) .. controls (97.78,66.42) and (98.97,65.22) .. (100.44,65.22) .. controls (101.92,65.22) and (103.11,66.42) .. (103.11,67.89) .. controls (103.11,69.36) and (101.92,70.56) .. (100.44,70.56) .. controls (98.97,70.56) and (97.78,69.36) .. (97.78,67.89) -- cycle ;
\draw  [fill={rgb, 255:red, 224; green, 16; blue, 57 }  ,fill opacity=1 ] (167.09,245.71) .. controls (167.09,244.38) and (168.28,243.31) .. (169.75,243.31) .. controls (171.21,243.31) and (172.4,244.38) .. (172.4,245.71) .. controls (172.4,247.03) and (171.21,248.1) .. (169.75,248.1) .. controls (168.28,248.1) and (167.09,247.03) .. (167.09,245.71) -- cycle ;
\draw   (155.67,74.11) .. controls (155.8,78.78) and (158.19,81.05) .. (162.86,80.92) -- (172.03,80.66) .. controls (178.7,80.48) and (182.09,82.72) .. (182.22,87.38) .. controls (182.09,82.72) and (185.36,80.3) .. (192.02,80.11)(189.02,80.19) -- (201.2,79.86) .. controls (205.86,79.73) and (208.13,77.33) .. (208,72.67) ;
\draw   (215.67,73.11) .. controls (215.8,77.78) and (218.19,80.05) .. (222.86,79.92) -- (232.03,79.66) .. controls (238.7,79.48) and (242.09,81.72) .. (242.22,86.38) .. controls (242.09,81.72) and (245.36,79.3) .. (252.02,79.11)(249.02,79.19) -- (261.2,78.86) .. controls (265.86,78.73) and (268.13,76.33) .. (268,71.67) ;
\draw   (275.67,73.11) .. controls (275.8,77.78) and (278.19,80.05) .. (282.86,79.92) -- (292.03,79.66) .. controls (298.7,79.48) and (302.09,81.72) .. (302.22,86.38) .. controls (302.09,81.72) and (305.36,79.3) .. (312.02,79.11)(309.02,79.19) -- (321.2,78.86) .. controls (325.86,78.73) and (328.13,76.33) .. (328,71.67) ;
\draw   (734.5,107.67) .. controls (734.57,110.69) and (736.11,112.16) .. (739.13,112.09) -- (739.13,112.09) .. controls (743.44,112) and (745.64,113.46) .. (745.71,116.48) .. controls (745.64,113.46) and (747.76,111.9) .. (752.07,111.8)(750.13,111.84) -- (752.07,111.8) .. controls (755.09,111.73) and (756.57,110.19) .. (756.5,107.17) ;
\draw  [dash pattern={on 0.84pt off 2.51pt}]  (755,38.67) -- (756.5,110.67) ;
\draw   (761.5,108.67) .. controls (761.57,111.69) and (763.11,113.16) .. (766.13,113.09) -- (766.13,113.09) .. controls (770.44,113) and (772.64,114.46) .. (772.71,117.48) .. controls (772.64,114.46) and (774.76,112.9) .. (779.07,112.8)(777.13,112.84) -- (779.07,112.8) .. controls (782.09,112.73) and (783.57,111.19) .. (783.5,108.17) ;
\draw   (334.67,72.11) .. controls (334.8,76.78) and (337.19,79.05) .. (341.86,78.92) -- (351.03,78.66) .. controls (357.7,78.48) and (361.09,80.72) .. (361.22,85.38) .. controls (361.09,80.72) and (364.36,78.3) .. (371.02,78.11)(368.02,78.19) -- (380.2,77.86) .. controls (384.86,77.73) and (387.13,75.33) .. (387,70.67) ;

\draw (103.44,16.4) node [anchor=north west][inner sep=0.75pt]  [color={rgb, 255:red, 65; green, 117; blue, 5 }  ,opacity=1 ]  {$Q_{0}$};
\draw (356.44,18.73) node [anchor=north west][inner sep=0.75pt]  [color={rgb, 255:red, 65; green, 117; blue, 5 }  ,opacity=1 ]  {$Q_{4}$};
\draw (49.44,91.07) node [anchor=north west][inner sep=0.75pt]    {$\Delta =wh=3h$};
\draw (659.94,135.92) node [anchor=north west][inner sep=0.75pt]  [color={rgb, 255:red, 65; green, 117; blue, 5 }  ,opacity=1 ,rotate=-359.9,xslant=-0.02]  {${\displaystyle Q_{4}}$};
\draw (559.44,30.51) node [anchor=north west][inner sep=0.75pt]  [color={rgb, 255:red, 130; green, 6; blue, 185 }  ,opacity=1 ,rotate=-359.9,xslant=-0.02]  {$\frac{4\chi _{4}}{\Delta }$};
\draw (67.72,258.77) node [anchor=north west][inner sep=0.75pt]  [font=\footnotesize]  {$-\frac{wh}{2}$};
\draw (319.97,197.39) node [anchor=north west][inner sep=0.75pt]  [font=\normalsize,color={rgb, 255:red, 74; green, 144; blue, 226 }  ,opacity=1 ,rotate=-359.9,xslant=-0.02]  {$T_{0}$};
\draw (395.81,177.57) node [anchor=north west][inner sep=0.75pt]    {$\ \frac{}{}$};
\draw (174.44,88.07) node [anchor=north west][inner sep=0.75pt]    {$\Delta $};
\draw (232.44,89.07) node [anchor=north west][inner sep=0.75pt]    {$\Delta $};
\draw (294.44,88.07) node [anchor=north west][inner sep=0.75pt]    {$\Delta $};
\draw (739.44,119.23) node [anchor=north west][inner sep=0.75pt]  [font=\small]  {$\frac{\Delta }{4}$};
\draw (766.44,120.23) node [anchor=north west][inner sep=0.75pt]  [font=\small]  {$\frac{\Delta }{4}$};
\draw (63.72,150.4) node [anchor=north west][inner sep=0.75pt]  [font=\scriptsize]  {$\frac{3wh}{2}$};
\draw (71.72,188.4) node [anchor=north west][inner sep=0.75pt]  [font=\small]  {$\Delta $};
\draw (419.44,109.4) node [anchor=north west][inner sep=0.75pt]  [color={rgb, 255:red, 65; green, 117; blue, 5 }  ,opacity=1 ]  {$Q_{4}$};
\draw (183.44,122.4) node [anchor=north west][inner sep=0.75pt]  [color={rgb, 255:red, 65; green, 117; blue, 5 }  ,opacity=1 ]  {$Q_{0}$};
\draw (353.44,87.07) node [anchor=north west][inner sep=0.75pt]    {$\Delta $};

\end{tikzpicture}
}
\caption[The common coordinate in the grid intervals]{%
We want to approximate $F(x)=\sum_{a\in\mathcal A}c_{a,0}\sigma_a(x)$,
where $\sigma_a(x)=(1-|x-a|/h)_+$ is the triangle centered at $a$.
Starting with an interval of length $\Delta=wh$ (with $w=3$ in this
example), we create $N\approx w$ intervals by shifting it as shown
in the figure, forming a set $\mathcal{C}$ of intervals.
The function $T_0$ serves to detect the relative position in each
interval, and $4\chi_4/\Delta$ serves to locate the interval within
the set $\mathcal{C}$.}
\label{fig:a6-ejemplo-coordenadas}

\medskip
\resizebox{1.0\linewidth}{!}{%
\tikzset{every picture/.style={line width=1.05pt}} 

\begin{tikzpicture}[x=0.75pt,y=0.75pt,yscale=-1,xscale=1]

\draw  [color={rgb, 255:red, 221; green, 244; blue, 202 }  ,draw opacity=1 ][fill={rgb, 255:red, 221; green, 244; blue, 202 }  ,fill opacity=1 ] (63.11,35) -- (123.67,35) -- (123.67,118.33) -- (63.11,118.33) -- cycle ;
\draw  [draw opacity=0] (63.11,118.33) -- (139.67,118.33) -- (139.67,103) -- (63.11,103) -- cycle ; \draw   (63.11,118.33) -- (63.11,103)(83.11,118.33) -- (83.11,103)(103.11,118.33) -- (103.11,103)(123.11,118.33) -- (123.11,103) ; \draw   (63.11,118.33) -- (139.67,118.33) ; \draw    ;
\draw  [dash pattern={on 0.84pt off 2.51pt}]  (58.67,47) -- (145.67,47) ;
\draw [color={rgb, 255:red, 224; green, 16; blue, 57 }  ,draw opacity=1 ][line width=1.5]    (86.67,47) -- (63.11,118.33) ;
\draw [color={rgb, 255:red, 224; green, 16; blue, 57 }  ,draw opacity=1 ][line width=1.5]    (86.67,47) -- (103.11,118.33) ;
\draw [color={rgb, 255:red, 224; green, 16; blue, 57 }  ,draw opacity=1 ][line width=1.5]    (235.33,118.11) -- (103.11,118.33) ;
\draw [color={rgb, 255:red, 224; green, 16; blue, 57 }  ,draw opacity=1 ][line width=1.5]    (63.11,118.33) -- (13.56,118.56) ;
\draw [color={rgb, 255:red, 74; green, 144; blue, 226 }  ,draw opacity=1 ][line width=1.5]    (181.67,66.67) .. controls (221.07,37.12) and (277.28,36.89) .. (314.64,69.59) ;
\draw [shift={(316.33,71.11)}, rotate = 222.51] [color={rgb, 255:red, 74; green, 144; blue, 226 }  ,draw opacity=1 ][line width=1.5]    (14.21,-4.28) .. controls (9.04,-1.82) and (4.3,-0.39) .. (0,0) .. controls (4.3,0.39) and (9.04,1.82) .. (14.21,4.28)   ;
\draw  [color={rgb, 255:red, 221; green, 244; blue, 202 }  ,draw opacity=1 ][fill={rgb, 255:red, 221; green, 244; blue, 202 }  ,fill opacity=1 ] (338.11,31.44) -- (398.67,31.44) -- (398.67,117) -- (338.11,117) -- cycle ;
\draw  [color={rgb, 255:red, 221; green, 244; blue, 202 }  ,draw opacity=1 ][fill={rgb, 255:red, 221; green, 244; blue, 202 }  ,fill opacity=1 ] (575.67,29.44) -- (640.22,29.44) -- (640.22,117) -- (575.67,117) -- cycle ;
\draw  [draw opacity=0] (338.11,117) -- (685.67,117) -- (685.67,105.33) -- (338.11,105.33) -- cycle ; \draw   (338.11,117) -- (338.11,105.33)(358.11,117) -- (358.11,105.33)(378.11,117) -- (378.11,105.33)(398.11,117) -- (398.11,105.33)(418.11,117) -- (418.11,105.33)(438.11,117) -- (438.11,105.33)(458.11,117) -- (458.11,105.33)(478.11,117) -- (478.11,105.33)(498.11,117) -- (498.11,105.33)(518.11,117) -- (518.11,105.33)(538.11,117) -- (538.11,105.33)(558.11,117) -- (558.11,105.33)(578.11,117) -- (578.11,105.33)(598.11,117) -- (598.11,105.33)(618.11,117) -- (618.11,105.33)(638.11,117) -- (638.11,105.33)(658.11,117) -- (658.11,105.33)(678.11,117) -- (678.11,105.33) ; \draw   (338.11,117) -- (685.67,117) ; \draw    ;
\draw [color={rgb, 255:red, 224; green, 16; blue, 57 }  ,draw opacity=1 ][line width=1.5]    (361.67,45.67) -- (338.11,117) ;
\draw [color={rgb, 255:red, 224; green, 16; blue, 57 }  ,draw opacity=1 ][line width=1.5]    (361.67,45.67) -- (378.11,117) ;
\draw [color={rgb, 255:red, 224; green, 16; blue, 57 }  ,draw opacity=1 ][line width=1.5]    (478.11,117) -- (378.11,117) ;
\draw [color={rgb, 255:red, 224; green, 16; blue, 57 }  ,draw opacity=1 ][line width=1.5]    (578.11,117) -- (518.11,117) ;
\draw [color={rgb, 255:red, 224; green, 16; blue, 57 }  ,draw opacity=1 ][line width=1.5]    (501.67,45.67) -- (478.11,117) ;
\draw [color={rgb, 255:red, 224; green, 16; blue, 57 }  ,draw opacity=1 ][line width=1.5]    (601.67,45.67) -- (578.11,117) ;
\draw [color={rgb, 255:red, 224; green, 16; blue, 57 }  ,draw opacity=1 ][line width=1.5]    (601.67,45.67) -- (618.11,117) ;
\draw [color={rgb, 255:red, 224; green, 16; blue, 57 }  ,draw opacity=1 ][line width=1.5]    (501.67,45.67) -- (518.11,117) ;
\draw  [color={rgb, 255:red, 221; green, 244; blue, 202 }  ,draw opacity=1 ][fill={rgb, 255:red, 221; green, 244; blue, 202 }  ,fill opacity=1 ] (89.44,188.78) -- (150,188.78) -- (150,274.33) -- (89.44,274.33) -- cycle ;
\draw  [color={rgb, 255:red, 221; green, 244; blue, 202 }  ,draw opacity=1 ][fill={rgb, 255:red, 221; green, 244; blue, 202 }  ,fill opacity=1 ] (327,186.78) -- (391.56,186.78) -- (391.56,274.33) -- (327,274.33) -- cycle ;
\draw  [draw opacity=0] (89.44,274.33) -- (437,274.33) -- (437,262.67) -- (89.44,262.67) -- cycle ; \draw   (89.44,274.33) -- (89.44,262.67)(109.44,274.33) -- (109.44,262.67)(129.44,274.33) -- (129.44,262.67)(149.44,274.33) -- (149.44,262.67)(169.44,274.33) -- (169.44,262.67)(189.44,274.33) -- (189.44,262.67)(209.44,274.33) -- (209.44,262.67)(229.44,274.33) -- (229.44,262.67)(249.44,274.33) -- (249.44,262.67)(269.44,274.33) -- (269.44,262.67)(289.44,274.33) -- (289.44,262.67)(309.44,274.33) -- (309.44,262.67)(329.44,274.33) -- (329.44,262.67)(349.44,274.33) -- (349.44,262.67)(369.44,274.33) -- (369.44,262.67)(389.44,274.33) -- (389.44,262.67)(409.44,274.33) -- (409.44,262.67)(429.44,274.33) -- (429.44,262.67) ; \draw   (89.44,274.33) -- (437,274.33) ; \draw    ;
\draw [color={rgb, 255:red, 224; green, 16; blue, 57 }  ,draw opacity=1 ][line width=1.5]    (329.44,274.33) -- (89.44,274.33) ;
\draw [color={rgb, 255:red, 224; green, 16; blue, 57 }  ,draw opacity=1 ][line width=1.5]    (353,203) -- (329.44,274.33) ;
\draw [color={rgb, 255:red, 224; green, 16; blue, 57 }  ,draw opacity=1 ][line width=1.5]    (353,203) -- (369.44,274.33) ;
\draw [color={rgb, 255:red, 130; green, 6; blue, 185 }  ,draw opacity=1 ][line width=1.5]    (504.67,135) .. controls (508.61,154.7) and (474.71,198.66) .. (419.22,207.62) ;
\draw [shift={(416.67,208)}, rotate = 352.01] [color={rgb, 255:red, 130; green, 6; blue, 185 }  ,draw opacity=1 ][line width=1.5]    (14.21,-4.28) .. controls (9.04,-1.82) and (4.3,-0.39) .. (0,0) .. controls (4.3,0.39) and (9.04,1.82) .. (14.21,4.28)   ;
\draw [color={rgb, 255:red, 224; green, 16; blue, 57 }  ,draw opacity=1 ][line width=1.5]    (678.11,117) -- (618.11,117) ;
\draw [color={rgb, 255:red, 224; green, 16; blue, 57 }  ,draw opacity=1 ][line width=1.5]    (429.44,274.33) -- (369.44,274.33) ;

\draw (85.11,16.4) node [anchor=north west][inner sep=0.75pt]  [color={rgb, 255:red, 65; green, 117; blue, 5 }  ,opacity=1 ]  {$Q_{0}$};
\draw (57,121.04) node [anchor=north west][inner sep=0.75pt]  [font=\small]  {$b_{0}$};
\draw (78,122.04) node [anchor=north west][inner sep=0.75pt]  [font=\small]  {$b_{1}$};
\draw (95.39,121.73) node [anchor=north west][inner sep=0.75pt]  [font=\small]  {$b_{2}$};
\draw (41.67,40.04) node [anchor=north west][inner sep=0.75pt]  [font=\small]  {$1$};
\draw (227.97,18.39) node [anchor=north west][inner sep=0.75pt]  [font=\normalsize,color={rgb, 255:red, 74; green, 144; blue, 226 }  ,opacity=1 ,rotate=-359.9,xslant=-0.02]  {$T_{0}$};
\draw (348.78,7.4) node [anchor=north west][inner sep=0.75pt]  [color={rgb, 255:red, 65; green, 117; blue, 5 }  ,opacity=1 ]  {$Q_{0}$};
\draw (594.44,5.07) node [anchor=north west][inner sep=0.75pt]  [color={rgb, 255:red, 65; green, 117; blue, 5 }  ,opacity=1 ]  {$Q_{4}$};
\draw (133.97,82.73) node [anchor=north west][inner sep=0.75pt]  [font=\normalsize,color={rgb, 255:red, 224; green, 16; blue, 57 }  ,opacity=1 ,rotate=-359.9,xslant=-0.02]  {$\sigma _{b_{1}}$};
\draw (417.97,34.73) node [anchor=north west][inner sep=0.75pt]  [font=\normalsize,color={rgb, 255:red, 224; green, 16; blue, 57 }  ,opacity=1 ,rotate=-359.9,xslant=-0.02]  {$\sigma _{b_{1}} \circ T_{0}$};
\draw (100.11,164.73) node [anchor=north west][inner sep=0.75pt]  [color={rgb, 255:red, 65; green, 117; blue, 5 }  ,opacity=1 ]  {$Q_{0}$};
\draw (345.78,162.4) node [anchor=north west][inner sep=0.75pt]  [color={rgb, 255:red, 65; green, 117; blue, 5 }  ,opacity=1 ]  {$Q_{4}$};
\draw (498.44,160.51) node [anchor=north west][inner sep=0.75pt]  [color={rgb, 255:red, 130; green, 6; blue, 185 }  ,opacity=1 ,rotate=-359.9,xslant=-0.02]  {$\times \frac{4\chi _{4}}{\Delta }$};
\draw (167.44,208.51) node [anchor=north west][inner sep=0.75pt]  [color={rgb, 255:red, 224; green, 16; blue, 57 }  ,opacity=1 ,rotate=-359.9,xslant=-0.02]  {$\frac{4\chi _{4}}{\Delta } \sigma _{b_{1}} \circ T_{0} =\sigma _{4\Delta +b_{1}}$};

\end{tikzpicture}

}
\caption[A single product for several coefficients]{%
In the intervals $Q_0$ and $Q_4$ of the group $\mathcal C$ in
Figure~\ref{fig:a6-ejemplo-coordenadas}, the first mark $b_1=h$
corresponds to $h$ and $13h$, respectively. In the figure,
$\beta=0$. Since $4\chi_0/\Delta$ selects the first triangle
and $4\chi_4/\Delta$ selects the last one, a single product gives
each its coefficient:
$c_{h,0}\sigma_h+c_{13h,0}\sigma_{13h}
 =
 \left(c_{h,0}\frac{4\chi_0}{\Delta}
      +c_{13h,0}\frac{4\chi_4}{\Delta}\right)
 (\sigma_h\circ T_0).$
In general, the mark at distance $b$ from the beginning of
$Q_\nu$ has coordinate $a=\nu\Delta+b$ and hence
coefficient $c_{\nu\Delta+b,\beta}$. For fixed $\beta$, the
same affine layer simultaneously computes, for
$b=0,h,\ldots,\Delta-h$, the normalized sums
$U_{\mathcal C,b,\beta}(x)
 =\sum_{Q_\nu\in\mathcal C}
   \frac{4c_{\nu\Delta+b,\beta}}{A\Delta}\,\chi_\nu(x).$
The following layers approximate the product of each
$U_{\mathcal C,b,\beta}(x)$ and the corresponding local function,
$v_{b,\beta}(T_0(x))$. For this purpose, one copy of the multiplier
is used for each $b$: all copies receive their two inputs in the
same layer and compute their products in parallel. Therefore,
approximating the $\Delta/h$ products requires as many hidden
layers as approximating a single one. In dimension one, a full
group contains $\#\mathcal C\asymp w$ intervals with
$\Delta/h\asymp w$ marks each. Thus, these layers process
on the order of $w^2$ coefficients per group, and
$O(1+h^{-1}/w^2)$ successive groups suffice to cover the
$h^{-1}+1$ grid points.}
\label{fig:a6-ejemplo-producto}
\end{figure}

\begin{lemma}\label{lem:localized-taylor-sums}
Fix integers $q\ge1$ and $m\ge0$. There exist constants
$C,K_0\ge1$ and an integer $w_0\ge2$, depending only on $q,m$,
with the following property.
Let $j\ge0$ be an integer, set $h=2^{-j}$ and
$\mathcal A=(h\mathbb Z^q)\cap[0,1]^q$, and define
\[
 \psi_\beta(z)=z^\beta\prod_{i=1}^q(1-|z_i|)_+,
 \qquad
 F(x)=\sum_{a\in\mathcal A}\sum_{|\beta|\le m}
 c_{a,\beta}\psi_\beta\!\left(\frac{x-a}{h}\right).
\]
Suppose $D,\tau>0$, $|c_{a,\beta}|\le D$ and
$A\ge C\frac{D}{h}$. Then, for every even $w\ge w_0$,
there is a GroupSort2 network $F^\circ$ of width at most $w$
and depth at most
\begin{equation}\label{eq:localized-taylor-depth}
 C\left\{j+1+\log\!\left(2+\frac{D}{\tau}\right)
 \left(1+\frac{h^{-q}}{w^2}
 +\frac{Dh^{-q-1}}{A w^{1+\frac{1}{q}}}\right)\right\}
\end{equation}
such that
\[
 \|F-K_0AF^\circ\|_{\infty,[0,1]^q}\le\tau.
\]
Every weight matrix of $F^\circ$, including the output matrix,
has absolute row sums at most one. All biases have absolute
value at most $C$, and the output bias is zero.
\end{lemma}

\begin{proof}
We choose $\eta=c\min\{1,\tau/D\}$, with $c>0$ sufficiently
small. To construct $\psi_\beta$, we multiply the $q$ factors
$(1-|z_i|)_+$ and, for each $i$, $\beta_i$ copies of
$\clip_{[-1,1]}(z_i)$. All the factors take values in $[-1,1]$,
and there are at most $q+m$ of them.
Applying Lemma~\ref{lemma:multiplication} successively, with
the normalization used in the proof of Lemma~\ref{lemma:monomial},
we obtain a network $p_\beta$ such that
$\|p_\beta-\psi_\beta/C_0\|_\infty\le\eta$, where $C_0\ge1$
is fixed. Its width is bounded, its depth is
$O(\log(e/\eta))$, and each row of its matrices satisfies
$\sum_i|w_i|\le1$. Moreover, $p_\beta$ vanishes outside
$(-1,1)^q$ because the approximate product is zero whenever
one of its factors is zero.

For $a\in h\mathbb Z^q$, we define
$v_{a,\beta}(x)=\clip_{[-h,h]}(h p_\beta((x-a)/h))$.
By the homogeneity of GroupSort, we obtain this network by
replacing the first affine map $Wz+b$ with $W(x-a)+hb$,
multiplying the subsequent biases by $h$, and clipping the
output to $[-h,h]$. Thus,
\begin{equation}\label{eq:taylor-local-features}
 \left\|v_{a,\beta}
 -\frac{h}{C_0}\psi_\beta\!\left(\frac{\cdot-a}{h}\right)
 \right\|_\infty\le h\eta,\qquad
 |v_{a,\beta}|\le h,\qquad
 v_{a,\beta}(x)=0\quad\text{if }\|x-a\|_\infty\ge h.
\end{equation}

For simplicity and without loss of generality, we may assume
that $w=2^{J/q}$ for some integer $J\geq 2q$.
We define $\Delta=h w^{1/q}$.
If $\Delta>1$, the grid $\mathcal A$ has
$(1+h^{-1})^q\lesssim w$ points. We can therefore compute all
the outputs $c_{a,\beta}v_{a,\beta}/(4Ah)$ in parallel.
Multiplication by $c_{a,\beta}/(4Ah)$ respects the matrix bound
because $|c_{a,\beta}|\le D$ and $A\ge CD/h$.

Suppose in what follows that $\Delta\le1$. For $\nu\in\mathbb Z^q$,
we define the cubes
\[
 Q_\nu=\prod_{i=1}^q[\Delta\nu_i,\Delta(\nu_i+1))
\]
which form a coarser grid than $\mathcal A$.
Each $a\in\mathcal A$ belongs to a unique $Q_\nu$ and can be
written uniquely as $a=\Delta\nu+b$, with
$b\in\mathcal B:=\{0,h,\ldots,(w^{1/q}-1)h\}^q$.
We set $c_{a,\beta}=0$ for $a\notin\mathcal A$.

We want to use the same networks $v_{b,\beta}$ in several
cubes, changing only their coefficients. To do so, we associate
with each $Q_\nu$ the function
\[
 \chi_\nu(x)=\min\left\{\frac{\Delta}{4},
 \left(\frac{\Delta}{2}
       -\operatorname{dist}_\infty(x,Q_\nu)\right)_+\right\},
 \quad
 \operatorname{dist}_\infty(x,Q_\nu)
 =\max_{1\le i\le q}
   \max\{\Delta\nu_i-x_i,\,0,\,x_i-\Delta(\nu_i+1)\}.
\]
This formula uses only maxima and minima of affine functions
whose rows have norm at most one. Therefore, $\chi_\nu$ requires
a fixed number of layers. Moreover, we have
\[
 \frac{4\chi_\nu(x)}{\Delta}
 =
 \begin{cases}
 1, & \operatorname{dist}_\infty(x,Q_\nu)\le \Delta/4,\\
 2-4\operatorname{dist}_\infty(x,Q_\nu)/\Delta,
    & \Delta/4<\operatorname{dist}_\infty(x,Q_\nu)<\Delta/2,\\
 0, & \operatorname{dist}_\infty(x,Q_\nu)\ge \Delta/2.
 \end{cases}
\]
Thus, $\chi_\nu(x)$ serves to select the coefficients of the
grid points that belong to $Q_\nu$.
If $a\in Q_\nu$ and $\|x-a\|_\infty\le h$, then
$\operatorname{dist}_\infty(x,Q_\nu)\le h\le\Delta/4$.
Hence $\chi_\nu=\Delta/4$ throughout the support of
$v_{a,\beta}$. For each $\sigma\in\{0,1,2,3\}^q$, we define
\[
 \mathcal Q_\sigma
 =\{Q_\nu:Q_\nu\cap\mathcal A\ne\varnothing,\
                  \nu-\sigma\in4\mathbb Z^q\},
\]
which represents the family of cubes obtained by translating
$Q_\sigma$ by integer multiples of $4\Delta$ in each coordinate.
Two distinct cubes in $\mathcal Q_\sigma$ are separated by
at least $3\Delta$ in some coordinate. Since $\chi_\nu$
vanishes at distance greater than or equal to $\Delta/2$ from
$Q_\nu$, for each $x$ there is at most one
$Q_\nu\in\mathcal Q_\sigma$ with $\chi_\nu(x)>0$.

For this family of cubes, we define
\[
 T_{\sigma,i}(x)=
 \operatorname{dist}\!\left(
 x_i-(\sigma_i-\tfrac12)\Delta,4\Delta\mathbb Z\right)
 -\frac{\Delta}{2},\qquad 1\le i\le q,
\]
which represents the position of $x_i$ relative to the left
face of the cube $Q_\nu\in\mathcal Q_\sigma$ when
$\chi_\nu(x)>0$.
If $Q_\nu\in\mathcal Q_\sigma$ and $\chi_\nu(x)>0$, then
$x_i-\nu_i\Delta+\Delta/2\in(0,2\Delta)$. Since
$\nu_i-\sigma_i\in4\mathbb Z$, the distance in the formula
is $x_i-\nu_i\Delta+\Delta/2$. Therefore,
$T_\sigma(x)=x-\Delta\nu$ whenever $\chi_\nu(x)>0$.
For $a=\Delta\nu+b$, the construction gives
$v_{b,\beta}(x-\Delta\nu)=v_{a,\beta}(x)$. Using the fact
that $\chi_\nu=\Delta/4$ on the support of $v_{a,\beta}$,
we obtain
\begin{equation}\label{eq:localized-cube-translation}
 \frac{4\chi_\nu(x)}{\Delta}\,v_{b,\beta}(T_\sigma(x))
 =v_{\Delta\nu+b,\beta}(x),
 \qquad x\in[0,1]^q.
\end{equation}

To compute $T_{\sigma,i}$, we start with
$z=x_i-(\sigma_i-\tfrac12)\Delta+4\in[0,8]$ and apply
$z\mapsto\min\{z,D'-z\}$ for $D'=8,4,\ldots,4\Delta$,
dividing $D'$ by two at each step. Since $4$ is a multiple
of $4\Delta$, we obtain the distance in the preceding formula.
We subtract $\Delta/2$ to obtain $T_{\sigma,i}$. This requires
$O(j+1)$ layers. For each $\sigma$, we compute the outputs
$v_{b,\beta}(T_\sigma(x))$ only once.

We now form subgroups with $N=\min\{w,\lfloor A\Delta/(4D)\rfloor\}$
cubes from each class $\mathcal Q_\sigma$.
The assumption $A\ge CD/h$ and the equality $\Delta=h w^{1/q}$ give
$A\Delta/(4D)\ge2$ for $C$ sufficiently large.
For each $\sigma$, we enumerate the cubes in $\mathcal Q_\sigma$
in an arbitrary order as $Q_{\nu_1},\ldots,Q_{\nu_{n_\sigma}}$,
where $n_\sigma=\#\mathcal Q_\sigma$. We take the first $N$
cubes, then the next $N$, and so on, obtaining
\[
 \mathcal C_{\sigma,k}
 =\{Q_{\nu_i}:(k-1)N<i\le\min\{kN,n_\sigma\}\},
 \qquad 1\le k\le\lceil n_\sigma/N\rceil.
\]
For each of these groups $\mathcal C=\mathcal C_{\sigma,k}$,
each $b\in\mathcal B$, and each $|\beta|\le m$, we compute
\[
 U_{\mathcal C,b,\beta}(x)
 =\sum_{Q_\nu\in\mathcal C}
   \frac{4c_{\Delta\nu+b,\beta}}{A\Delta}\chi_\nu(x).
\]
For a fixed position $b$ and a fixed index $\beta$,
$U_{\mathcal C,b,\beta}(x)$ selects the coefficient
$c_{\Delta\nu+b,\beta}$ of the cube in $\mathcal C$ to which
$x$ belongs. The choice of the number of cubes per group
ensures that the row computing this sum satisfies
\begin{equation}\label{eq:dense-table-row}
 \sum_{Q_\nu\in\mathcal C}
 \frac{4|c_{\Delta\nu+b,\beta}|}{A\Delta}
 \le\frac{4ND}{A\Delta}\le1.
\end{equation}
Moreover, $|U_{\mathcal C,b,\beta}(x)|\le D/A\le h$,
since at most one $\chi_\nu(x)$ is positive and its value
is at most $\Delta/4$.

We apply Lemma~\ref{lemma:multiplication} to the two inputs
$U_{\mathcal C,b,\beta}(x)$ and $v_{b,\beta}(T_\sigma(x))$,
both in $[-h,h]$. We multiply the biases of the multiplier
by $h$ and divide its output by four. The resulting network
has matrices of norm at most one and computes
\[
 Z_{\mathcal C,b,\beta}(x)
 =\frac h4\operatorname{Mult}_\eta\!\left(
 \frac{U_{\mathcal C,b,\beta}(x)}h,\,
 \frac{v_{b,\beta}(T_\sigma(x))}h\right).
\]
By \eqref{eq:localized-cube-translation}, the exact product
of the two inputs, divided by $4h$, is
$\sum_{Q_\nu\in\mathcal C}
 c_{\Delta\nu+b,\beta}v_{\Delta\nu+b,\beta}(x)/(4Ah)$.
The error in approximating the product is at most
$\eta|U_{\mathcal C,b,\beta}(x)|/4\le D\eta/(4A)$
by \eqref{eq:relative-multiplication}. Using also
\eqref{eq:taylor-local-features}, we obtain
\[
 \bigg|Z_{\mathcal C,b,\beta}(x)
 -\underbrace{\sum_{Q_\nu\in\mathcal C}
   \frac{c_{\Delta\nu+b,\beta}}{4C_0A}
   \psi_\beta\!\left(\frac{x-\Delta\nu-b}{h}\right)}_{=:P_{\mathcal C,b,\beta}(x)}
 \bigg|
 \le\frac{D\eta}{2A}.
\]
In this last sum, only the cube for which $\chi_\nu(x)>0$
can contribute. To sum the outputs, we fix
$\sigma\in\{0,1,2,3\}^q$, $\gamma\in\{0,1\}^q$, and
$|\beta|\le m$. We consider the points
$\mathcal B_\gamma=\{b\in\mathcal B:b/h-\gamma\in2\mathbb Z^q\}$.
We will prove that, for $x\in[0,1]^q$,
\[
 \left|
 \sum_{\mathcal C}\sum_{b\in\mathcal B_\gamma}
 \bigl(Z_{\mathcal C,b,\beta}(x)-P_{\mathcal C,b,\beta}(x)\bigr)
 \right|
 \le \frac{D\eta}{2A},
\]
where the first sum runs over the groups $\mathcal C_{\sigma,k}$
of $\mathcal Q_\sigma$.
Fix $x\in[0,1]^q$. For $Z_{\mathcal C,b,\beta}(x)$ to be nonzero,
both inputs of ${\rm Mult}_\eta$ must be nonzero. By the definition
of $U_{\mathcal C,b,\beta}$ and by \eqref{eq:taylor-local-features},
this requires
\[
 \chi_\nu(x)>0\ \text{for some }Q_\nu\in\mathcal C
 \quad\text{and}\quad
 \|T_\sigma(x)-b\|_\infty<h.
\]
The same two conditions are necessary for
$P_{\mathcal C,b,\beta}(x)$ to be nonzero. If one of its
summands is nonzero, then $\|x-\Delta\nu-b\|_\infty<h$.
Since $\Delta\nu+b\in Q_\nu$ and $h\le\Delta/4$, we have
$\chi_\nu(x)=\Delta/4$ and hence $T_\sigma(x)=x-\Delta\nu$.
The first condition determines at most one group $\mathcal C$,
since at most one cube in $\mathcal Q_\sigma$ satisfies
$\chi_\nu(x)>0$ and each cube belongs to a unique group.
The second determines at most one $b\in\mathcal B_\gamma$,
since two distinct points of $\mathcal B_\gamma$ are at least
$2h$ apart in the maximum norm. Thus, for the fixed
$\sigma,\gamma,\beta$, both functions vanish except possibly
for the same single pair $(\mathcal C,b)$, which proves the claim.

For fixed $\sigma,\gamma,\beta$, we want to compute at each
point $x$ the sum
\[
 \sum_{k=1}^{\lceil n_\sigma/N\rceil}
 \sum_{b\in\mathcal B_\gamma}
 Z_{\mathcal C_{\sigma,k},b,\beta}(x).
\]
The first sum runs over the groups of cubes in $\mathcal Q_\sigma$,
and the second runs over the positions $b$ in $\mathcal B_\gamma$.
Since at most one pair $(k,b)$ gives a nonzero output, this
sum is equal to
\[
 \max_{k,\,b\in\mathcal B_\gamma}
 \max\{0,Z_{\mathcal C_{\sigma,k},b,\beta}(x)\}
 +
 \min_{k,\,b\in\mathcal B_\gamma}
 \min\{0,Z_{\mathcal C_{\sigma,k},b,\beta}(x)\},
\]
where $k$ runs over the same groups as in the sum.
If $\mathcal Q_\sigma$ is empty, we take both values to be zero.
The network computes these two quantities as it goes through
the groups. Fix $\sigma,\gamma,\beta$ and $x$. The network goes
through the groups $\mathcal C_{\sigma,1},\ldots,\mathcal C_{\sigma,K}$,
where $K=\lceil n_\sigma/N\rceil$. For each position
$b\in\mathcal B_\gamma$, it stores two numbers, the largest and
the smallest among zero and the outputs already computed at
that position.

Before processing the first group, both numbers are zero,
$p_{0,b}=n_{0,b}=0$. When the network computes the output
of group $k$, it compares it with the two stored numbers
and obtains
\[
 p_{k,b}=\max\{p_{k-1,b},Z_{\mathcal C_{\sigma,k},b,\beta}(x)\},
 \qquad
 n_{k,b}=\min\{n_{k-1,b},Z_{\mathcal C_{\sigma,k},b,\beta}(x)\}.
\]
These comparisons are made in parallel for all positions
$b\in\mathcal B_\gamma$.

After the last group, $p_{K,b}$ and $n_{K,b}$ are the maximum
and the minimum for position $b$, including zero. It remains
to compare the positions with one another. By pairwise
comparisons, the network computes
\[
 p_{\sigma,\gamma,\beta}(x)=\max_{b\in\mathcal B_\gamma}p_{K,b},
 \qquad
 n_{\sigma,\gamma,\beta}(x)=\min_{b\in\mathcal B_\gamma}n_{K,b}.
\]
GroupSort2 performs each comparison by sorting its two inputs.
Since at most one pair $(k,b)$ gives a nonzero output,
\[
 p_{\sigma,\gamma,\beta}(x)+n_{\sigma,\gamma,\beta}(x)
 =
 \sum_{k=1}^{K}\sum_{b\in\mathcal B_\gamma}
 Z_{\mathcal C_{\sigma,k},b,\beta}(x).
\]

There are $M=4^q2^q\binom{q+m}{m}$ choices of
$(\sigma,\gamma,\beta)$. To construct $F^\circ$, we carry
to the last layer the two values $p_{\sigma,\gamma,\beta}(x)$
and $n_{\sigma,\gamma,\beta}(x)$ computed for each
$\sigma\in\{0,1,2,3\}^q$, $\gamma\in\{0,1\}^q$, and
$\beta\in\mathbb N_0^q$ with $|\beta|\le m$. There are
$M=4^q2^q\binom{q+m}{m}$ triples of these indices.
We assign the weight $1/(2M)$ to each of the $2M$ inputs
and set the bias to zero. Thus,
\[
 F^\circ(x)=\frac{1}{2M}
 \sum_{\sigma,\gamma,\beta}
 \bigl(p_{\sigma,\gamma,\beta}(x)
      +n_{\sigma,\gamma,\beta}(x)\bigr).
\]
The sum of the absolute values of the weights in this row is one.

To identify the function approximated by this network, we replace
each output $Z_{\mathcal C,b,\beta}$ with the function
$P_{\mathcal C,b,\beta}$ that it approximates. Each $a\in\mathcal A$
belongs to a unique cube $Q_\nu$ and determines
$b=a-\Delta\nu\in\mathcal B$. The index $\nu$ determines the
class $\sigma$ of its cube, and $b/h$ determines the class
$\gamma$ of $b$. Moreover, $Q_\nu$ belongs to a unique group
$\mathcal C_{\sigma,k}$. Therefore, for each $\beta$, the
summand corresponding to $a$ appears exactly once, giving
\[
 \sum_{\sigma,\gamma,\beta}
 \sum_{k=1}^{\lceil n_\sigma/N\rceil}
 \sum_{b\in\mathcal B_\gamma}P_{\mathcal C_{\sigma,k},b,\beta}(x)
 =
 \frac{1}{4C_0A}\sum_{a\in\mathcal A}\sum_{|\beta|\le m}
 c_{a,\beta}\psi_\beta\!\left(\frac{x-a}{h}\right)
 =\frac{F(x)}{4C_0A}.
\]
Since $p_{\sigma,\gamma,\beta}+n_{\sigma,\gamma,\beta}$ is
the sum of the outputs $Z_{\mathcal C_{\sigma,k},b,\beta}$
for that triple, the preceding error bound gives
\[
 \left\|2M F^\circ-\frac{F}{4C_0A}\right\|_{\infty,[0,1]^q}
 \le M\frac{D\eta}{2A}.
\]
Multiplying by $4C_0A$ and taking $K_0=8C_0M$, we obtain
$\|K_0AF^\circ-F\|_{\infty,[0,1]^q}\le 2C_0MD\eta\le\tau$.

If $\Delta>1$, we group the points $a\in\mathcal A$ directly
according to $a/h-\gamma\in2\mathbb Z^q$. For each
$\gamma,\beta$ and each $x$, at most one of these points satisfies
$\|x-a\|_\infty<h$. We therefore compute the maximum with zero
and the minimum with zero of the outputs
$c_{a,\beta}v_{a,\beta}/(4Ah)$ in each group. We add zeros
to complete the $2M$ inputs of the average. By
\eqref{eq:taylor-local-features}, the error for each group is
at most $D\eta/(4A)$, so we obtain the same final bound with
the same $K_0$.

To count the layers when $\Delta\le1$, recall that
$n_\sigma=\#\mathcal Q_\sigma$. There are at most
$1+\Delta^{-1}$ cube indices per coordinate. Therefore,
$\sum_\sigma n_\sigma\le(1+\Delta^{-1})^q\lesssim h^{-q}/w$.
Since each group contains $N$ cubes except possibly the last
group of each $\mathcal Q_\sigma$, the total number of groups
satisfies
\[
 \sum_\sigma\left\lceil\frac{n_\sigma}{N}\right\rceil
 \le 4^q+\frac{1}{N}\sum_\sigma n_\sigma
 \lesssim
 1+\frac{h^{-q}}{w^2}
   +\frac{Dh^{-q-1}}{Aw^{1+1/q}},
\]
where we have used $N\asymp\min\{w,Ah w^{1/q}/D\}$.

For each $\sigma$, we compute $T_\sigma$ and the outputs
$v_{b,\beta}(T_\sigma(x))$ only once, in
$O(j+1+\log(e/\eta))$ layers. In each group, we compute the
functions $\chi_\nu$ in parallel, then the sums
$U_{\mathcal C,b,\beta}$ and their approximate products
$Z_{\mathcal C,b,\beta}$, and finally update the maxima and
minima. This requires $O(\log(e/\eta))$ layers per group.
The final comparisons over $b$ are made in pairs in parallel
and require $O(\log(2+w))=O(j+1)$ layers, since
$w^{1/q}\le h^{-1}=2^j$. If $\Delta>1$, the comparisons are made
over the $(1+h^{-1})^q$ points of $\mathcal A$ and also
require $O(j+1)$ layers. Substituting
$\log(e/\eta)\lesssim\log(2+D/\tau)$, we obtain
\eqref{eq:localized-taylor-depth}.

In each group, we compute at most $N\le w$ functions
$\chi_\nu$ and at most $w\binom{q+m}{m}$ products in
parallel. We carry the inputs, the outputs
$v_{b,\beta}(T_\sigma(x))$, and the maxima and minima already
computed as pairs of equal coordinates. The width is therefore
at most $Cw\le w$ by the choice of $C_*$. The bound
\eqref{eq:dense-table-row}, the normalization of the multiplier,
and the choice of the final row ensure that the sum of the
absolute values in each row is at most one. The biases are
bounded independently of $j,w,D,A,\tau$ because the translations
used are bounded and the biases of the multiplication networks
are multiplied by $h\le1$.
\end{proof}

\subsection{Approximation of monomials}
In the following, we construct GroupSort2 networks that approximate the monomials appearing in \eqref{eq:differenceQ}. We further present multiplication networks that allow us to multiply the hat functions with the local Taylor polynomials within the proof of Lemma \ref{lem:localized-taylor-sums}.

The definition of $\GS$ gives
\begin{equation}\label{eq:analytic-proof-sorting}
 \max\{u,v\}=(1,0)\GS(u,v),\qquad
 \min\{u,v\}=(0,1)\GS(u,v),\qquad
 \GS(u,u)=(u,u).
\end{equation}
GroupSort is also positively homogeneous, so
\begin{equation}\label{eq:analytic-sorting-homogeneity}
\GS(tv)=t\GS(v),\qquad t\ge0.
\end{equation}
We store scalar outputs as equal pairs (see \eqref{eq:analytic-proof-sorting}), before the final affine
map. Additional layers can copy these pairs, and unused
coordinates can be set to zero pairs. This gives any larger
prescribed depth and even width without changing the function
or the matrix bounds.
We use these identities to implement the square interpolant of \citet[Proposition~2]{Yarotsky2017} under our matrix constraints \eqref{eq:network-constraints}.

\begin{lemma}\label{lemma:squared}
For $0<\epsilon<1$, there is a network
$g_\epsilon\in\widetilde{\mathcal F}(L,4,2,2)$ with
$L\lesssim\log\!\left(\frac{e}{\epsilon}\right)$ such that
\[
 \sup_{x\in[0,1]}|g_\epsilon(x)-x^2|\le\epsilon.
\]
All hidden biases have absolute value at most one.
\end{lemma}

\begin{proof}
Let $m\ge1$ and set $a_i=2-2^{-i}$. For $u\in[0,1]$, define
\[
 v_0(u)=q_0(u)=u,\qquad
 v_i(u)=\min\{v_{i-1}(u),2^{1-i}-v_{i-1}(u)\},\qquad
 q_i(u)=q_{i-1}(u)-2^{-i}v_i(u),
\]
where $q_m$ is the dyadic linear interpolant of the square. Hence
$ \|q_m-(\cdot)^2\|_{\infty,[0,1]}\le 4^{-m-1},$ see \citet[Proposition~2]{Yarotsky2017}.
It remains to implement this recursion under the network constraints
\eqref{eq:network-constraints}. We store the vector
$z_i(u)=\bigl(\frac{q_i(u)}{a_i},\frac{q_i(u)}{a_i},v_i(u),v_i(u)\bigr)^\top.$
The first hidden layer computes $z_0=(u,u,u,u)^\top$.
Suppose that $z_{i-1}$ is available, and write
$s=\frac{q_{i-1}}{a_{i-1}}$ and $v=v_{i-1}$.
One affine map followed by GroupSort2 computes
$ \bigl(\max\{v,2^{1-i}-v\},\ v_i,\ s,\ s\bigr)^\top.$
Its matrix has absolute row sums equal to one and its only
nonzero bias is $2^{1-i}\le1$. A second affine map followed
by GroupSort2 computes
\[
 z_i=\left(
 \frac{a_{i-1}s-2^{-i}v_i}{a_i},
 \frac{a_{i-1}s-2^{-i}v_i}{a_i},v_i,v_i\right)^\top.
\]
The two rows forming its first coordinate pair have absolute
row sum $\frac{a_{i-1}+2^{-i}}{a_i}=1$, and the other rows copy
$v_i$. Thus every hidden matrix has
$\ell^\infty$ operator norm one. The output map selects the
first coordinate and multiplies it by $a_m<2$.
All biases are bounded by one.

Choose the smallest $m\ge1$ such that $4^{-m-1}\le\epsilon$.
The construction has width four and $2m+1$ hidden layers,
which gives the asserted depth.
\end{proof}

The identity
\begin{align*}
    uv=\left(\frac{|u+v|}{2}\right)^2-\left(\frac{|u-v|}{2}\right)^2
\end{align*}
for $(u,v)\in[-1,1]^2$ allows us to use the squaring networks of Lemma~\ref{lemma:squared} to approximate multiplication.

\begin{lemma}\label{lemma:multiplication}
For $0<\epsilon<1$, there is a network
$\operatorname{Mult}_\epsilon\in\widetilde{\mathcal F}(L,8,4,4)$
with $L\lesssim\log\!\left(\frac{e}{\epsilon}\right)$ such that
\begin{equation}\label{eq:relative-multiplication}
 |\operatorname{Mult}_\epsilon(u,v)-uv|
 \le\epsilon\min\{|u|,|v|\},\qquad (u,v)\in[-1,1]^2.
\end{equation}
Moreover, $\operatorname{Mult}_\epsilon(\R^2)\subseteq[-1,1]$ and
$\operatorname{Mult}_\epsilon(u,0)=\operatorname{Mult}_\epsilon(0,v)=0$
for all real $u,v$. All hidden biases have absolute value at most one.
\end{lemma}

\begin{proof}
Let $q_m$ be the dyadic piecewise affine interpolant of $x^2$ used in Lemma~\ref{lemma:squared}. Its error $e_m(x)=q_m(x)-x^2$ has $\Lip(e_m|_{[0,1]})\le2^{-m}$, since on each grid interval its slope differs from $2x$ by at most $2^{-m}$.
Choose $2^{-m}\le\epsilon$ and define
\[
 \operatorname{Mult}_\epsilon(u,v)
 =\clip_{[-1,1]}\left\{q_m\left(\frac{|u+v|}{2}\right)
                  -q_m\left(\frac{|u-v|}{2}\right)\right\}.
\]
We note that 
\begin{align*}
    \left|\frac{|u+v|}{2} - \frac{|u-v|}{2}\right| \leq \min\{|u|,|v|\}.
\end{align*}
The Lipschitz bound on $e_m$ and nonexpansiveness of clipping prove \eqref{eq:relative-multiplication}. Equal arguments give the exact-zero assertions. Write $q=q_m$.

To represent this function, first compute
$\frac{|u+v|}{2}$ and $\frac{|u-v|}{2}$ using GroupSort2, and
apply two copies of $q$ in parallel. The hidden
layers have width at most $8$ and matrix norms
at most one. If $W$ is the output row of $q$,
the row $\frac{(W,-W)}{4}$ has norm at most one by
Lemma~\ref{lemma:squared}. The output biases
of the two copies cancel. This gives
\[
 s=\frac{1}{4}\left(
 q\left(\frac{|u+v|}{2}\right)
 -q\left(\frac{|u-v|}{2}\right)\right).
\]
Two additional GroupSort2 layers compute
$\clip_{[-\frac{1}{4},\frac{1}{4}]}(s)$ using maxima
and minima with the constants $\pm\frac{1}{4}$.
The output layer multiplies the result by $4$.
Thus all hidden matrix norms are at most one,
the output matrix norm is $4$, and all biases
are bounded by one. The construction adds
a bounded number of hidden layers to two parallel copies of $q$,
giving width eight and depth
$L\lesssim\log\!\left(\frac{e}{\epsilon}\right)$.
\end{proof}
Applying the multiplication networks of Lemma \ref{lemma:multiplication} pairwise allows us to approximate monomials successively.

\begin{lemma}\label{lemma:monomial}
Let $\beta\in\mathbb N_0^q$ and $0<\epsilon<1$. For $j=|\beta|\ge 0$ and
$r=\#\{i:\beta_i>0\}$, there is a network
$p_\beta\in\widetilde{\mathcal F}(L,2r+8,4^{j-1},4^{j-1})$ with
\[
 L\lesssim j\log\!\left(\frac{ej}{\epsilon}\right),\qquad
 \sup_{u\in[-1,1]^q}|p_\beta(u)-u^\beta|\le\epsilon.
\]
It satisfies $p_\beta(\R^q)\subseteq[-1,1]$, and all hidden
biases have absolute value at most one.
\end{lemma}

\begin{proof}
For $|\beta|\le1$, the claim follows from constant and
coordinate networks. Hence, assume $j=|\beta|\ge2$. Choose indices
$\iota_1,\ldots,\iota_j\in\{1,\ldots,q\}$ such that
each index $i$ occurs exactly $\beta_i$ times. Then
\[
 u^\beta=\prod_{i=1}^q u_i^{\beta_i}
        =\prod_{\ell=1}^j u_{\iota_\ell}.
\]
Set $\delta=\frac{\epsilon}{j-1}$ and define
\begin{equation}\label{eq:analytic-monomial-recursion}
 p_1(u)=u_{\iota_1},\qquad
 p_\ell(u)=\operatorname{Mult}_\delta
       (p_{\ell-1}(u),u_{\iota_\ell}),
 \quad 2\le\ell\le j.
\end{equation}
By Lemma~\ref{lemma:multiplication}, $p_\ell$ takes
values in $[-1,1]$ for $\ell\ge2$, and
\[
 \sup_{u\in[-1,1]^q}
 \left|p_\ell(u)-\prod_{h=1}^\ell u_{\iota_h}\right|
 \le(\ell-1)\delta.
\]
Thus $p_j$ approximates $u^\beta$ within $\epsilon$.
The coordinates $u_i$ used in later multiplications are
copied through the preceding layers as pairs $(u_i,u_i)$.
Since $\GS(u_i,u_i)=(u_i,u_i)$, each copy is
preserved exactly. This requires at most $2r$ additional
neurons per layer, with matrix norm at most one. The resulting
width is at most $2r+8$, and the depth is at most
$Cj\log\!\left(\frac{ej}{\epsilon}\right)$.

It remains to control the matrix norms. Each multiplication
network has hidden matrix norms at most one and output
matrix norm at most four. By positive homogeneity of GroupSort2,
rescale the $\ell$th multiplication so that its hidden
activations and biases are divided by $4^{\ell-1}$,
and its output by $4^\ell$, for $1\le\ell\le j-1$.
Rescale the carried coordinates by the same factors.
Then every hidden matrix has norm at most one, while
a final affine map of norm $4^{j-1}$ recovers $p_j$.
All biases remain bounded by one. This proves that the network lies in $\widetilde{\mathcal F}(L,2r+8,4^{j-1},4^{j-1}) $.
\end{proof}

\subsection{Proof of Theorem \ref{th:holder-approximation} for \texorpdfstring{$0<\alpha \leq 1$}{0 < alpha <= 1}}
\label{sec:holder-cones}
To treat the case $0<\alpha\le1$, we approximate the function
from its closest function values on a grid, using maxima and minima.
The following lemma shows that this construction achieves
error $\epsilon$ with fixed width and depth
$O(\epsilon^{-q/\alpha})$.
In particular, for $\alpha=1$ it allows the output row norm
to remain bounded independently of $\epsilon$. This completes the proof of Theorem~\ref{th:holder-approximation}.
\begin{lemma}\label{lem:fixed-width-cone}
For $0<\alpha\le1$, a fixed sufficiently large even width and depth $O(\epsilon^{-q/\alpha})$ suffice for \eqref{eq:holder-approximation} with $B=C(1+\epsilon^{(\alpha-1)/\alpha})$ and $b=C_b B$. 
\end{lemma}
\begin{proof}
Let $0<\alpha\le1$ and $0<\epsilon\le1$. Choose the smallest
integer $j\ge0$ such that $2^{-j}\le\epsilon^{1/\alpha}$, and
set $h=2^{-j}$. Then
$\frac12\epsilon^{1/\alpha}<h\le\epsilon^{1/\alpha}$.
Dividing each coordinate interval $[0,1]$ into $2^j$ equal
intervals gives the grid
\[
 \mathcal A
 =\{(k_1h,\ldots,k_qh):k_1,\ldots,k_q\in\{0,\ldots,2^j\}\}.
\]
It has $N=(2^j+1)^q=(1+h^{-1})^q$ points. For every
$x\in[0,1]^q$, rounding each coordinate to the nearest multiple
of $h$ gives a point $a_x\in\mathcal A$ with
$\|x-a_x\|_\infty\le h/2$.

Take $K\ge1$ such that
$|f(x)-f(y)|\le K\|x-y\|_\infty^\alpha$ on $[0,1]^q$.
Set $B=Kh^{\alpha-1}$, so that $Bh=Kh^\alpha$.
For each $a\in\mathcal A$, define
\[
 G_a(x)=f(a)+B\|x-a\|_\infty,
 \qquad
 G(x)=\min_{a\in\mathcal A}G_a(x).
\]
We first bound $G(x)$ from above. Using the nearest grid point $a_x \in \arg\min_{a \in \mathcal{A}}\|x-a\|_\infty$ and the H\"older bound on $f$, we obtain
\[
 G(x)\le G_{a_x}(x)
 \le f(x)+K\|x-a_x\|_\infty^\alpha+B\|x-a_x\|_\infty
 \le f(x)+2Kh^\alpha.
\]

For the lower bound, we must bound $G_a(x)$ from below for
every $a\in\mathcal A$. Write $r=\|x-a\|_\infty$.
If $r\le h$, then $f(a)\ge f(x)-Kr^\alpha\ge f(x)-Kh^\alpha$.
Adding $Br\ge0$ gives $G_a(x)\ge f(x)-Kh^\alpha$.
If $r\ge h$, then $1-\alpha\ge0$ gives
\[
 Br=Kh^{\alpha-1}r
   =Kr^\alpha\left(\frac rh\right)^{1-\alpha}
   \ge Kr^\alpha.
\]
Since $f(a)\ge f(x)-Kr^\alpha$, we have $G_a(x)\ge f(x)$
in this case. Taking the minimum over $a$ gives
$G(x)\ge f(x)-Kh^\alpha$. Finally, put
$g=\clip_{[-H,H]}\circ G$. Since $f(x)\in[-H,H]$, clipping
can only decrease the error. Hence
\[
 \|g-f\|_{\infty,[0,1]^q}\le2Kh^\alpha\le2K\epsilon.
\]

We now construct the network. We compute $G/B$ in the hidden
layers and multiply by $B$ in the output row. For each grid
point $a$,
\[
 \frac{G_a(x)}{B}
 =\max_{\substack{1\le i\le q\\ \sigma\in\{-1,1\}}}
   \left\{\frac{f(a)}{B}+\sigma(x_i-a_i)\right\}.
\]
Each affine function in this maximum has one nonzero weight,
equal to $1$ or $-1$. Its bias has absolute value at most
$H/B+1\le H+1$, since $B\ge1$. Sorting two values with
GroupSort2 gives their maximum and minimum. Pairwise comparisons
therefore compute the maximum of these $2q$ affine functions
in a number of layers depending only on $q$.

To compute the minimum over the $N$ grid points, choose
$p=\min\{N,\lfloor w/C_q\rfloor\}$, where $C_q$ is large
enough to allow the computations below within width $w$.
Divide the grid into $p$ groups of at most $\lceil N/p\rceil$
points each. In each group, start with the value $G_a(x)/B$
at its first point. At every subsequent point, compute
$G_a(x)/B$ and take its minimum with the value already stored.
Carry each input coordinate and each stored value as a pair
of equal coordinates, which GroupSort2 preserves.

The $p$ groups are computed in parallel. Each new grid point
in a group costs at most $C_q$ layers, so after at most
$C_q\lceil N/p\rceil$ layers we have the minimum for each
group. Comparing these $p$ minima in pairs gives $G(x)/B$
in at most $\lceil\log_2 p\rceil$ further comparison rounds.
Since $N\lesssim h^{-q}$, the total depth is at most
\[
 C_q\left(1+\frac{h^{-q}}{w}+\log(2+w)\right).
\]
All hidden rows have absolute sum at most one. The final row
multiplies $G(x)/B$ by $B$, and output clipping gives $g$.
For fixed sufficiently large even $w$, the depth is therefore
$O(\epsilon^{-q/\alpha})$. The output row has norm
$B=Kh^{\alpha-1}$, which equals $K$ for every $\epsilon$
when $\alpha=1$.
\end{proof}

\section{Covering numbers}
\subsection{The output matrix bound for Lipschitz functions}
\label{sec:bounded-endpoint}
The following result shows that, for $q\ge3$ and
$L\asymp[\log(ew)]^2$ and a bounded layerwise constraint $B$, an approximation rate $(w^2L)^{-\frac1q}$ cannot be derived uniformly over the $C^1$
ball, even with unrestricted biases. The proof uses a different
entropy bound with  polynomial dependence on $B$. 
\begin{proposition}\label{prop:bounded-endpoint}
Let $q\ge3$, $H,B\ge1$ be fixed and let
$L\asymp[\log(ew)]^2$.
As $w$ tends to infinity through even integers,
\[
 (w^2L)^{1/q}
 \sup_{\|f\|_{C^1([0,1]^q)}\le H}
 \inf_{g\in\bigcup_{b\ge1}\mathcal F(L,w,b,B,H)}
 \|f-g\|_{\infty,[0,1]^q}\longrightarrow\infty.
\]
\end{proposition}
\begin{proof}
To simplify the notation, let $\mathcal F_*=\bigcup_{b\ge1}\mathcal F(L,w,b,B,H)$ for this proof.

\medskip
\noindent\textit{Step 1: Independence of $b$.} For $g\in\mathcal F_*$,
\[
 g(x)=\clip_{[-H,H]}\!\left(
 W_{L+1}\operatorname{GS}_2\!\left(
 W_L\operatorname{GS}_2\!\left(
 \cdots\operatorname{GS}_2(W_1x+b_1)\cdots
 \right)+b_L\right)+b_{L+1}\right),
\]
with $\max_i\sum_k|(W_\ell)_{ik}|\le1$ for
$1\le\ell\le L$ and $\|W_{L+1}\|_1\le B$.
We define $y_0(x)=z_0(x)=x$,
and recursively $y_\ell=\operatorname{GS}_2(W_\ell y_{\ell-1}+b_\ell)$
and $z_\ell=y_\ell-y_\ell(0)$.
By induction and the Lipschitz property of GroupSort2,
$z_\ell(0)=0$,
$\|z_\ell(x)-z_\ell(x')\|_\infty\le\|x-x'\|_\infty$
and $\|z_\ell\|_\infty\le1$.
For $c_\ell=W_\ell y_{\ell-1}(0)+b_\ell$,
\[
 z_\ell=
 \operatorname{GS}_2(W_\ell z_{\ell-1}+c_\ell)
 -\operatorname{GS}_2(c_\ell).
\]
Interchanging the rows of each pair if necessary,
we assume $c_{\ell,2r-1}\ge c_{\ell,2r}$.
Let $v_1=(W_\ell z_{\ell-1})_{2r-1}$,
$v_2=(W_\ell z_{\ell-1})_{2r}$ and
$d_{\ell,r}=c_{\ell,2r-1}-c_{\ell,2r}$.
With the maximum in the first coordinate of the pair,
\begin{equation}\label{eq:centered-layer-pairs}
 z_{\ell,2r-1}=\max\{v_1,v_2-d_{\ell,r}\},
 \qquad
 z_{\ell,2r}=\min\{v_1+d_{\ell,r},v_2\}.
\end{equation}
Since $|v_i|\le1$, the condition $d_{\ell,r}\ge2$ implies
$v_2-d_{\ell,r}\le v_1$ and $v_1+d_{\ell,r}\ge v_2$.
By \eqref{eq:centered-layer-pairs}, replacing $d_{\ell,r}$
by $\min\{d_{\ell,r},2\}$ preserves $z_\ell$.

Likewise, $g=\clip_{[-H,H]}(W_{L+1}z_L+c)$,
where $c=W_{L+1}y_L(0)+b_{L+1}$.
The inequality $|W_{L+1}z_L|\le B$ implies
\[
 \clip_{[-H,H]}(W_{L+1}z_L+c)
 =
 \clip_{[-H,H]}\!\left(
 W_{L+1}z_L+\clip_{[-H-B,H+B]}(c)\right).
\]
Every function in $\mathcal F_*$ therefore admits
this description with $d_{\ell,r}\in[0,2]$
and $c\in[-H-B,H+B]$, independently of $b$.

\medskip
\noindent\textit{Step 2: Covering the rows.}
Let $z=(z_1,\ldots,z_p)$ with $z(0)=0$ and
$\|z(x)-z(x')\|_\infty\le\|x-x'\|_\infty$.
For $\mathcal V_z=\{a\cdot z:\|a\|_1\le1\}$, we will prove
\begin{equation}\label{eq:centered-row-cover}
 \log\mathcal N(t,\mathcal V_z,\|\cdot\|_\infty)
 \lesssim t^{-\frac{2q}{q+2}}
 \left[\log\!\left(\frac{e(p+1)}{t}\right)\right]^2,
 \qquad 0<t<1.
\end{equation}

Let $\mathcal D_{t^{\frac{2}{q+2}}},\mathcal D_t
\subset[0,1]^q$ be such that
$\#\mathcal D_\rho\lesssim\rho^{-q}$ and
$\sup_{x\in[0,1]^q}\min_{y\in\mathcal D_\rho}
\|x-y\|_\infty\le\rho$ for
$\rho\in\{t^{\frac{2}{q+2}},t\}$.
For each $x\in\mathcal D_t$, we choose
$y(x)\in\mathcal D_{t^{\frac{2}{q+2}}}$ with
$\|x-y(x)\|_\infty\le t^{\frac{2}{q+2}}$.
It follows that
$\|z(x)-z(y(x))\|_\infty\le t^{\frac{2}{q+2}}$.
We apply \citet[Corollary~5]{Zhang2002Covering} with
bound one for the $\ell^1$ norm of the coefficients,
bound $t^{\frac{2}{q+2}}$ for the maximum norm
of the inputs and accuracy $t$.
The class
\[
 \left\{
 \bigl(a\cdot(z(x)-z(y(x)))\bigr)_{x\in\mathcal D_t}:
 \|a\|_1\le1
 \right\}
\]
admits a cover of radius $t$ and cardinality $N$
in the maximum norm, with
\begin{equation}\label{eq:centered-increment-cover}
\begin{aligned}
 \log N
 &\lesssim t^{-\frac{2q}{q+2}}\log(e(p+1))
       \log\!\left(2+t^{\frac{2}{q+2}-1}\#\mathcal D_t\right)
 \lesssim t^{-\frac{2q}{q+2}}\log(e(p+1))
       \log\!\left(\frac{e}{t}\right).
\end{aligned}
\end{equation}

For $v\in\mathcal V_z$, we choose a center of this
cover and round
$(v(y))_{y\in\mathcal D_{t^{\frac{2}{q+2}}}}$
with error at most $t$.
Since $\|v\|_\infty\le1$, the logarithm of the number
of choices is at most, up to a constant factor,
\[
 \#\mathcal D_{t^{\frac{2}{q+2}}}
 \log\!\left(\frac{e}{t}\right)+\log N.
\]
If $v,v'\in\mathcal V_z$ have the same rounded
values and the same center, then
\[
 \max_{y\in\mathcal D_{t^{\frac{2}{q+2}}}}
 |v(y)-v'(y)|\le2t,
 \qquad
 \max_{x\in\mathcal D_t}
 |(v-v')(x)-(v-v')(y(x))|\le2t.
\]
Therefore, $\max_{x\in\mathcal D_t}|v(x)-v'(x)|\le4t$
and, by Lipschitz continuity, $\|v-v'\|_\infty\le6t$.
Choosing a representative of each nonempty class
and applying \eqref{eq:centered-increment-cover}, we obtain
\[
 \log\mathcal N(6t,\mathcal V_z,\|\cdot\|_\infty)
 \lesssim t^{-\frac{2q}{q+2}}\log\!\left(\frac{e}{t}\right)
 +t^{-\frac{2q}{q+2}}\log(e(p+1))
       \log\!\left(\frac{e}{t}\right).
\]
Replacing $t$ by $\frac{t}{6}$ gives
\eqref{eq:centered-row-cover}.

\medskip
\noindent\textit{Step 3: Covering $\mathcal F_*$.}
We proceed layer by layer as in the proof of
Lemma~\ref{lem:network-entropy}, now using
\eqref{eq:centered-row-cover}.
Let $\widetilde z_0=z_0$.
For each $\widetilde z_{\ell-1}$, we choose the rows
of $\widetilde W_\ell$ so that
\[
 \|(W_\ell-\widetilde W_\ell)
       \widetilde z_{\ell-1}\|_\infty\le6t,
 \qquad
 \max_i\sum_k|(\widetilde W_\ell)_{ik}|\le1.
\]
We choose $\widetilde d_{\ell,r}\in[0,2]$
on a fixed grid of spacing at most $t$, with
$|d_{\ell,r}-\widetilde d_{\ell,r}|\le t$,
and define $\widetilde z_\ell$ through
\eqref{eq:centered-layer-pairs}.

By induction, $\widetilde z_\ell(0)=0$ and
$\|\widetilde z_\ell(x)-\widetilde z_\ell(x')\|_\infty
\le\|x-x'\|_\infty$.
For $E_\ell=\|z_\ell-\widetilde z_\ell\|_\infty$,
the matrix constraints and
\eqref{eq:centered-layer-pairs} imply
\begin{equation}\label{eq:centered-layer-error}
\begin{aligned}
 \|W_\ell z_{\ell-1}
       -\widetilde W_\ell\widetilde z_{\ell-1}\|_\infty
 &\le E_{\ell-1}+6t,\quad {\rm and}\quad
 E_\ell\le E_{\ell-1}+7t.
\end{aligned}
\end{equation}
Since $E_0=0$, it follows that $E_L\le7Lt$.
We write $W_{L+1}=Ba$, with $\|a\|_1\le1$.
Applying \eqref{eq:centered-row-cover} to
$a\cdot\widetilde z_L$ and approximating $c$ within
$[-H-B,H+B]$ with error at most $t$,
\eqref{eq:centered-layer-error} gives
$\|g-\widetilde g\|_\infty\le B(7L+6)t+t$.
For $t=\frac{u}{32B(L+1)}$, this bound is smaller than
$\frac{u}{2}$.

There are at most $w(L+1)$ rows, each with at most
$\max\{q,w\}$ inputs.
By \eqref{eq:centered-row-cover} and the rounding
of $d_{\ell,r}$ and $c$, the number $M$
of functions $\widetilde g$ satisfies
\[
 \log M
 \lesssim w(L+1)t^{-\frac{2q}{q+2}}
 \left[\log\!\left(\frac{e(q+w)}{t}\right)\right]^2.
\]
We choose a function from $\mathcal F_*$
in each ball of radius $\frac{u}{2}$ that intersects the class.
The triangle inequality and the choice of $t$ give
\begin{equation}\label{eq:bounded-output-entropy}
 \log\mathcal N(u,\mathcal F_*,\|\cdot\|_\infty)
 \lesssim w(L+1)^{1+\frac{2q}{q+2}}u^{-\frac{2q}{q+2}}
 \left[\log\!\left(\frac{ew(L+1)}{u}\right)\right]^2.
\end{equation}

\medskip
\noindent\textit{Step 4: Lower bound on the error.}
Let $\mathcal H=\{f:\|f\|_{C^1([0,1]^q)}\le H\}$ and
$E_{w,L}=\sup_{f\in\mathcal H}
\inf_{g\in\mathcal F_*}\|f-g\|_\infty$.
If $E_{w,L}\le u$, a cover of $\mathcal F_*$
of radius $u$ provides a cover of $\mathcal H$
of radius $3u$, with centers in the ambient space.
The entropy lower bound from
Section~\ref{sec:entropy} implies
$u^{-q}\lesssim
\log\mathcal N(u,\mathcal F_*,\|\cdot\|_\infty)$.

For each $m>0$, define $u_m=m(w^2L)^{-\frac{1}{q}}$.
For $L\asymp[\log(ew)]^2$, the bound
\eqref{eq:bounded-output-entropy} and $E_{w,L}\le u_m$ imply
\[
 1\lesssim
 w^{-\frac{q-2}{q+2}}
 (L+1)^{\frac{2(q+1)}{q+2}}[\log(ew)]^2
 \longrightarrow0,\qquad {\rm for}\  q>2.
\]
Therefore, $E_{w,L}>m(w^2L)^{-\frac{1}{q}}$ for
sufficiently large $w$. Since $m>0$ is arbitrary,
$(w^2L)^{\frac{1}{q}}E_{w,L}\longrightarrow\infty$.
\end{proof}

\subsection{Proof of the covering number bound}\label{sec:entropy-proofs}

\begin{proof}[Proof of Lemma~\ref{lem:network-entropy}]
Let $\mathcal F_k$ denote the subclass of $\widetilde{\mathcal F}(L,w,b,B)$ with fixed architecture $(L,k)$.
Let $g_\theta,g_{\theta'}$ have the same architecture and let their corresponding parameters differ by at most $\eta$.
For $x\in[-R,R]^q$, the hidden outputs satisfy
$\|z_\ell\|_\infty,\|z'_\ell\|_\infty\le R+\ell b$.
Each affine row has at most $q+w$ weight entries and one bias. Nonexpansiveness of GroupSort therefore gives
\[
 \|z_{\ell+1}-z'_{\ell+1}\|_\infty
 \le\|z_\ell-z'_\ell\|_\infty
       +\eta(q+w+1)(1+R+\ell b).
\]
Iterating and using the output row norm bound yields
\begin{equation}\label{eq:parameter-sensitivity}
 \|g_\theta-g_{\theta'}\|_{\infty,[-R,R]^q}
 \le\eta(q+w+1)(L+1)(1+B)(1+R+(L+1)b).
\end{equation}
For $L=0$ the same estimate follows directly from the affine map.

Mimicking \citet[Lemma~5]{SchmidtHieber2020}, we take $\eta$ to make the right side of \eqref{eq:parameter-sensitivity} equal to $u$, and replace every parameter $v$ by
$Q_\eta(v)=\eta\operatorname{sgn}(v)\lfloor\frac{|v|}{\eta}\rfloor$. Since $|Q_\eta(v)|\le|v|$, all entrywise bounds and row-sum constraints are preserved. Each of the $s=S(L,k)$ parameter positions has at most $1+\frac{2b}{\eta}$ possible values. These  new networks belong to the class and form a $u$-cover, so
$\mathcal N(u,\mathcal F_k,\|\cdot\|_{\infty,[-R,R]^q})\le(1+\frac{2b}{\eta})^s$.
Taking logarithms proves the claim for fixed $k$.  Any network with $L$ hidden layers of even widths at most $w$
can be represented with width exactly $w$ in every hidden layer:
add pairs of zero coordinates and set all their incoming and
outgoing weights to zero. These additional pairs remain zero
under GroupSort2 and do not change the network output or the
matrix bounds. Thus the covering bound applies to the whole
class with hidden widths at most $w$, using the parameter
count $S(L,w)$.

Clipping the outputs of both a network and its approximant
does not increase their uniform distance, since
$|\clip_{[-H,H]}(a)-\clip_{[-H,H]}(b)|\le|a-b|$.
Likewise, restricting the domain to $[0,1]^q$ cannot increase
the uniform error. Both operations therefore preserve the
covering radius.
\end{proof}

\section{Least-squares bounds}\label{sec:statistics}
We establish risk bounds for  least-squares estimators and apply them to
flow matching on each time interval.
Let $D_1,\ldots,D_n$ be independent copies of $D$, and let
$(X,T)$ be defined jointly with $D$, with $T\in\R^q$ and
$\E|T|^2<\infty$. For a class $\mathcal H$ of
measurable functions taking values in $\R^q$, define
\[
 \ell_g(D)=\E[|g(X)-T|^2\mid D],\qquad
 m(X)=\E[T\mid X],\qquad
 r(g)=\E|g(X)-m(X)|^2.
\]
Write $P_n\ell_g=\frac{1}{n}\sum_{i=1}^n\ell_g(D_i)$
and $\|g\|_\infty=\sup_x|g(x)|$.

\begin{lemma}\label{lem:bounded-regression}
Suppose that $|T|\le K$ almost surely and
$\sup_{g\in\mathcal H}\|g\|_\infty\le M_0$, and set
$A=M_0+K$. Let $u>0$ satisfy
$N=\mathcal N(u,\mathcal H,\|\cdot\|_\infty)<\infty$.
For every $x>0$, with probability at least $1-e^{-x}$,
\begin{equation}\label{eq:bounded-comparison}
 \frac{1}{2}r(g)-a_n(x)
 \le P_n(\ell_g-\ell_m)
 \le\frac{3}{2}r(g)+a_n(x),
 \qquad g\in\mathcal H,
\end{equation}
where $a_n(x)=11A^2\frac{\log(2N)+x}{n}+5Au$.
Consequently, any measurable 
$\widehat g\in\arg\min_{g\in\mathcal H}P_n\ell_g$ satisfies
\begin{equation}\label{eq:bounded-oracle}
 \E r(\widehat g)
 \le3\inf_{g\in\mathcal H}r(g)+4a_n(1).
\end{equation}
\end{lemma}

\begin{proof}
Put $Q_g=|g(X)-T|^2-|m(X)-T|^2$ and
$\Delta_g(D)=\E[Q_g\mid D]=\ell_g(D)-\ell_m(D)$.
Since $\E[T-m(X)\mid X]=0$, we have
$\E\Delta_g=r(g)$. Moreover, $|m(X)|\le K$ and
\[
 Q_g=\langle g(X)-m(X),g(X)+m(X)-2T\rangle
\]
give $|Q_g|\le3A|g(X)-m(X)|\le3A^2$.
Conditional Jensen's inequality therefore yields
\[
 \E[\Delta_g^2]\le\E[Q_g^2]\le9A^2r(g),
 \qquad
 |\Delta_g-r(g)|\le4A^2.
\]
Bernstein's inequality
\citep[Theorem~2.10]{BoucheronLugosiMassart2013}
gives, for every $t>0$,
\[
 \Prb\left(
 |P_n\Delta_g-r(g)|
 >
 \sqrt{\frac{18A^2r(g)t}{n}}
 +\frac{4A^2t}{3n}
 \right)\le2e^{-t}.
\]
Since $\sqrt{18ab}\le\frac{a}{2}+9b$, the threshold is
at most $\frac{1}{2}r(g)+\frac{11A^2t}{n}$.
Choose a uniform $u$-cover $g_1,\ldots,g_N\in\mathcal H$.
Taking $t=\log(2N)+x$ and a union bound gives this
comparison simultaneously at all its centers.

For any $g\in\mathcal H$, choose $g_j$ with
$\|g-g_j\|_\infty\le u$. The difference of squares gives
\[
 |\ell_g-\ell_{g_j}|\le2Au,
 \qquad |r(g)-r(g_j)|\le2Au.
\]
Thus, on the same event,
\[
\begin{aligned}
 |P_n\Delta_g-r(g)|
 &\le |P_n\Delta_{g_j}-r(g_j)|+4Au\\
 &\le\frac{1}{2}r(g)
      +11A^2\frac{\log(2N)+x}{n}+5Au,
\end{aligned}
\]
which proves \eqref{eq:bounded-comparison}.

Since $\widehat g$ minimizes $P_n\ell_g$, applying
\eqref{eq:bounded-comparison} to $\widehat g$ and any
$g\in\mathcal H$ gives
\[
 \frac{1}{2}r(\widehat g)-a_n(x)
 \le P_n(\ell_{\widehat g}-\ell_m)
 \le P_n(\ell_g-\ell_m)
 \le\frac{3}{2}r(g)+a_n(x).
\]
Hence, with probability at least $1-e^{-x}$,
$r(\widehat g)\le3\inf_g r(g)+4a_n(x)$.
Set $W=r(\widehat g)-3\inf_g r(g)$.
Since $a_n(x)=a_n(0)+\frac{11A^2x}{n}$, integration yields
\[
 \E W
 \le4a_n(0)+
 \int_0^\infty
 \Prb\bigl(W>4a_n(0)+t\bigr)\dd t
 \le4a_n(0)+\frac{44A^2}{n}
 =4a_n(1).
\]
This proves \eqref{eq:bounded-oracle}.
\end{proof}

We next apply Lemma~\ref{lem:bounded-regression} to the
flow-matching loss on a time interval. The estimator minimizes
the original loss; the response is truncated only in the proof.

\begin{lemma}\label{prop:interval-oracle}
Suppose $\Prb(|Y|>r)\le C_0e^{-c_0r^2}$ for $r\ge0$,
and let $I\subset[0,1]$ be an interval of positive length.
Let $\mathcal H$ be a class of fields on $I\times\R^q$
with $\sup_{g\in\mathcal H}\|g\|_\infty
\lesssim\sqrt{\log(en)}$ and
\begin{equation}\label{eq:block-entropy-hyp}
 \log\mathcal N(u,\mathcal H,\|\cdot\|_\infty)
 \lesssim s\log\!\left(\frac{en}{u}\right),
 \qquad 0<u<1,
\end{equation}
for some $s\ge1$.
If $\widehat g$ is a measurable minimizer of \eqref{eq:loss}
restricted to $I$ over $\mathcal H$, then, for every fixed
$a>0$,
\begin{equation}\label{eq:interval-oracle}
\begin{aligned}
 \E\int_I
 \|\widehat g_t-v_t^P\|_{L^2(P_t)}^2\dd t
 &\le3\inf_{g\in\mathcal H}\int_I
 \|g_t-v_t^P\|_{L^2(P_t)}^2\dd t+C_a|I|
\left(\frac{s[\log(en)]^2}{n}+n^{-a}\right).
\end{aligned}
\end{equation}
The bound holds uniformly over $I$.
\end{lemma}

\begin{proof}
Let $\tau$ be uniform on $I$, independently of $Y,Z$,
and take $D=Y$, $X=(\tau,U_\tau)$ and $T=Y-Z$.
Then $m(X)=v_\tau^P(U_\tau)$ and
\begin{equation}\label{eq:interval-excess-identity}
 r(g)=\frac{1}{|I|}\int_I
       \|g_t-v_t^P\|_{L^2(P_t)}^2\dd t.
\end{equation}
Moreover, $P_n\ell_g$ is the empirical loss on $I$
divided by $|I|$, so it is minimized by $\widehat g$.

Set $K=C_a\sqrt{\log(en)}$ and let $T^K$ be the Euclidean
projection of $T$ onto $\{z:|z|\le K\}$. Define
\[
 \ell_g^K(Y)=\E[|g(X)-T^K|^2\mid Y],\qquad
 m_K(X)=\E[T^K\mid X],\qquad
 r_K(g)=\E|g(X)-m_K(X)|^2.
\]
With $M_0=\sup_{g\in\mathcal H}\|g\|_\infty$, the
difference of squares gives
$|\ell_g(Y)-\ell_g^K(Y)|\le D_K(Y)$ uniformly in $g$, where
\[
 D_K(Y)=2\E\!\left[
 (M_0^2+|T|^2)\ind_{\{|T|>K\}}\mid Y\right].
\]
The tail assumptions imply $\Prb(|T|>r)\le Ce^{-cr^2}$.
Also, conditional Jensen's inequality gives
\[
 \|m-m_K\|_{L^2(P_X)}
 \le\bigl(\E|T-T^K|^2\bigr)^{\frac{1}{2}},
\]
and therefore
\[
 \sup_{g\in\mathcal H}|r(g)-r_K(g)|
 \le
 \left(2M_0+2(\E|T|^2)^{\frac{1}{2}}\right)
 \bigl(\E|T-T^K|^2\bigr)^{\frac{1}{2}}.
\]
Choosing $C_a$ sufficiently large ensures
\[
 \E D_K(Y)+
 \sup_{g\in\mathcal H}|r(g)-r_K(g)|
 \lesssim n^{-a-2}.
\]

Since $\widehat g$ minimizes the original empirical loss,
for every $g\in\mathcal H$,
\[
 P_n\ell_{\widehat g}^K
 \le P_n\ell_{\widehat g}+P_nD_K
 \le P_n\ell_g+P_nD_K
 \le P_n\ell_g^K+2P_nD_K.
\]
Apply \eqref{eq:bounded-comparison} to $T^K$, with
$A=M_0+K$. With probability at least $1-e^{-x}$,
simultaneously for every $g\in\mathcal H$,
\[
\begin{aligned}
 \frac{1}{2}r_K(\widehat g)-a_n(x)
 &\le P_n(\ell_{\widehat g}^K-\ell_{m_K}^K)\\
 &\le P_n(\ell_g^K-\ell_{m_K}^K)+2P_nD_K\le\frac{3}{2}r_K(g)+a_n(x)+2P_nD_K.
\end{aligned}
\]
Consequently,
\[
 \Prb\left(
 r_K(\widehat g)-3\inf_{g\in\mathcal H}r_K(g)
 -4P_nD_K>4a_n(x)
 \right)\le e^{-x}.
\]
Integrating as in the proof of
Lemma~\ref{lem:bounded-regression} and using
$\E P_nD_K=\E D_K(Y)$ gives
\[
 \E r_K(\widehat g)
 \le3\inf_{g\in\mathcal H}r_K(g)
   +4a_n(1)+4\E D_K(Y).
\]

Choose $u=(en)^{-a-4}$. By \eqref{eq:block-entropy-hyp},
$\log(2N)\lesssim s\log(en)$, while
$M_0+K\lesssim\sqrt{\log(en)}$. Hence
\[
 a_n(1)\lesssim
 \frac{s[\log(en)]^2}{n}+n^{-a-2}.
\]
Replacing $r_K$ by $r$ and multiplying by $|I|$
in \eqref{eq:interval-excess-identity} proves
\eqref{eq:interval-oracle}.
\end{proof}
\section{Proofs for flow matching}\label{sec:flow-proofs}

\subsection{The population velocity}\label{sec:geometry}
First, we are going to derive a representation of the population vector field $v^P$ that eases the access to smoothness properties in case of $P \in \mathcal{P}_{\alpha}$.
We apply the change from linear to trigonometric interpolation used in \citet[Appendix~D]{Chen_2025_ICCV}. Related parametrizations appear in \citet{AlbergoVandenEijnden2023,SalimansHo2022}. 
For $0\le t\le1$ we set
\begin{equation}\label{eq:coordinates}
 r_t=\sqrt{t^2+(1-t)^2}.
\end{equation}
The map $(Y,Z)\mapsto\bigl(\frac{U_t}{r_t},\frac{(1-t)Y-tZ}{r_t}\bigr)$ is orthogonal. Under \eqref{eq:model}, the joint density of the transformed vectors at $(x,z)$ is
\begin{equation}\label{eq:joint}
 \varphi_q(x)\varphi_q(z)\,
 \frac{f\bigl(\frac{tx+(1-t)z}{r_t}\bigr)}{\int f\varphi_q}.
\end{equation}
Integrating out $z$ therefore gives
\begin{equation}\label{eq:rho}
 \frac{m}{M}r_t^{-q}\varphi_q\!\left(\frac{x}{r_t}\right)
 \le \frac{\mathrm dP_t}{\mathrm dx}(x)
 \le \frac{M}{m}r_t^{-q}\varphi_q\!\left(\frac{x}{r_t}\right).
\end{equation}
For $G\sim N(0,I_q)$, define the vector
\begin{equation}\label{eq:holder-gaussian-ratio}
 \mathcal V(a,z)=\frac{\E[G f(z+aG)]}{\E[f(z+aG)]},
 \qquad {\rm for}\qquad (a,z)\in\R\times\R^q .
\end{equation}
Since $r_t^{\prime}=\frac{2 t-1}{r_t}$ and consequently $Y-Z=\frac{r_t'}{r_t}U_t+\frac{(1-t)Y-tZ}{r_t^2}$, conditioning in \eqref{eq:joint} on $U_t = x$ gives
\begin{equation}\label{eq:population-field}
 v_t^P(x)=\frac{r_t'}{r_t}x+
   \frac{1}{r_t}\mathcal V\!\left(\frac{1-t}{r_t},\frac{tx}{r_t^2}\right).
\end{equation}
Here $\frac{r_t'}{r_t}=\frac{2t-1}{r_t^2}$ has absolute value at most one, $r_t^{-1}\le\sqrt2$, and $|\mathcal V_i(a,z)|\le\frac{M}{m}\E|G_i|$. Thus each component of $v_t^P$ is bounded by $2R$ on $[-R,R]^q$ for $R$ as in \eqref{eq:scales}.

\subsection{Joint H\"older regularity}\label{sec:velocity-holder}

\begin{lemma}\label{lem:velocity-holder}
Under \eqref{eq:model}, $\|\mathcal V\|_{C^\alpha(\R^{q+1})}\lesssim1$ and $\mathcal V(0,z)=0$. For $R\ge2$, the function
\begin{equation}\label{eq:holder-physical-remainder}
 (u,x)\longmapsto v_t^P(x)-\frac{r_t'}{r_t}x,\qquad
 t=\frac{u}{1+R},\quad (u,x)\in[0,1+R]\times[-R,R]^q,
\end{equation}
has $C^\alpha$ norm bounded uniformly in $R$ and $P$ on the displayed rectangle.
\end{lemma}

\begin{proof}
First, we show that $\mathcal{V}$ inherits Hölder regularity from $f$. To this end, consider the numerator and denominator in
\eqref{eq:holder-gaussian-ratio} as functions of $(a,z)$.
Differentiating their integrands with respect to $z$ differentiates
$f$, while each derivative with respect to $a$ also introduces
a coordinate of $G$. Since the derivatives up to order  $\lfloor \alpha \rfloor$ of $f$ are bounded and
$\E|G|^p<\infty$ for every $p>0$, these derivatives can be taken under the
expectations.   If $\alpha$ is not an integer, put
$k=\lfloor\alpha\rfloor$ and $\beta=\alpha-k$.
The derivatives of order $k$ of the numerator and denominator
are expectations of derivatives of $f$ multiplied by polynomials
in $G$. For the numerator, with $\ell+|\gamma|=k$, we obtain
\begin{equation*}
    \partial_a^\ell \partial_z^\gamma \big( \mathbb{E}[G_i f(z+a G)] \big)=\mathbb{E}\big[G_i(G^{\top} \nabla)^\ell \partial^\gamma f(z+a G)\big].
\end{equation*}
The derivative of the denominator is constructed accordingly. Since each derivative of $f$ of order $k$ is
$\beta$-H\"older continuous and
\[
 |z-z'+(a-a')G|
 \le(1+|G|)(|z-z'|+|a-a'|),
\]
their increments are bounded by a constant times
$(|z-z'|+|a-a'|)^\beta$ after taking expectations.
For integer $\alpha$, boundedness of the derivatives of $f$
through order $\alpha$ gives the required derivative bounds
directly. Thus both the numerator and denominator have bounded
$C^\alpha$ norms. Since the denominator satisfies
$\E f(z+aG)\ge m$, the quotient also has a bounded
$C^\alpha$ norm. At $a=0$, its numerator vanishes because
$f(z)$ is deterministic and $\E G=0$.

Our purpose is now to bound the $C^\alpha$ norm of the $\mathcal V$ multiplied with $\frac{1}{r_t}$, as in \eqref{eq:population-field},
independently of $R$. To this end, we fix
$a(u)=\frac{1-t}{r_t}$ and $z(u,x)=\frac{tx}{r_t^2}$,
with the rescaled time $u=(1+R)t$.
Since $r_t^2=t^2+(1-t)^2\ge1/2$, the functions
$r_t^{-1}$, $(1-t)/r_t$ and $t/r_t^2$ have bounded
derivatives of each fixed order on $[0,1]$.
A fortiori,
\[
 \left|\partial_u^k z_i(u,x)\right|
 =
 \frac{|x_i|}{(1+R)^k}
 \left|\frac{\dd^k}{\dd t^k}\left(\frac{t}{r_t^2}\right)\right|
 \le C_k\frac{R}{(1+R)^k}\le C_k,
 \qquad k\ge1.
\]
Since $\|x\|_\infty\le R$, we obtain
\[
 \left|\partial_u^k\partial_{x_i}z_i(u,x)\right|
 =
 \frac{1}{(1+R)^k}
 \left|\frac{\dd^k}{\dd t^k}\left(\frac{t}{r_t^2}\right)\right|
 \le C_k,\qquad k\ge0.
\]
The other spatial derivatives of $z_i$ are zero.
The derivatives of $a(u)$ and $r_{u/(1+R)}^{-1}$
are also bounded uniformly in $R$.

We have proved that
$\|\mathcal V\|_{C^\alpha(\R^{q+1})}\lesssim1$.
We apply this bound and the bounds on the derivatives of $a$,
$z$ and $r_t^{-1}$ to the chain and product rules.
If $\alpha$ is not an integer, the derivatives of order
$\lfloor\alpha\rfloor$ of $\mathcal V$ are
H\"older with exponent $\alpha-\lfloor\alpha\rfloor$.
Composition with the Lipschitz functions $a$ and $z$
preserves this exponent. We thus obtain
\[
 \left\|
 (u,x)\longmapsto
 \frac{1}{r_t}\mathcal V(a(u),z(u,x))
 \right\|_{C^\alpha([0,1+R]\times[-R,R]^q)}
 \lesssim1.
\]

\end{proof}

\begin{lemma}\label{lem:regularity}
For $P\in\mathcal P_\alpha$ the field $v^P$ in \eqref{eq:population-field} has linear growth and a unique flow with time-$t$ law $P_t$ when started at $N(0,I_q)$. In particular its terminal law is $P$.
Uniformly over the model,
\begin{equation}\label{eq:w-lip}
 \sup_{0\le t\le1}\Lip_x(v_t^P)\lesssim1.
\end{equation}
\end{lemma}
\begin{proof}
Lemma~\ref{lem:velocity-holder} gives bounded $\mathcal V$ and $D_z\mathcal V$. Since $2^{-1/2}\le r_t\le1$, formula~\eqref{eq:population-field} gives both \eqref{eq:w-lip} and
\begin{equation}\label{eq:regularity-pop-growth}
 |v_t^P(x)|\lesssim 1+|x|,\qquad 0\le t\le1.
\end{equation}
As $t\uparrow1$, we have $r_t\to1$ and
$r_t'/r_t=(2t-1)/r_t^2\to1$, so the first term in
\eqref{eq:population-field} tends to $x$. In the second term,
$(1-t)/r_t\to0$ and $tx/r_t^2\to x$. By continuity of
$\mathcal V$, this term tends to $\mathcal V(0,x)=0$, since
$\E[Gf(x)]=f(x)\E G=0$. Thus $v_t^P(x)\to x$.

The linear growth bound and the integrability of the spatial
Lipschitz constant give a unique global flow on $[0,1]$
\citep[Remark~2.4]{AmbrosioCrippa2014}.
To identify its time marginals, take
$\psi\in C_c^\infty(\R^q)$. Since $\nabla\psi$ is bounded
and $\E|Y-Z|<\infty$, differentiation under the expectation
is justified. The definition of $v_t^P$ then gives
\[
 \frac{\dd}{\dd t}\int\psi\,\dd P_t
 =\E[\nabla\psi(U_t)\cdot(Y-Z)]
 =\int\nabla\psi(x)\cdot v_t^P(x)\,P_t(\dd x).
\]
Thus $(P_t)$ solves the continuity equation with initial law
$P_0=N(0,I_q)$. The laws of the flow started from $N(0,I_q)$
solve the same equation. Uniqueness of nonnegative
measure-valued solutions identifies these laws with $P_t$
\citep[Theorem~3.1]{AmbrosioCrippa2014}.
In particular, the terminal law is $P_1=P$.
\end{proof}

\subsection{Approximation of the velocity on each interval}
\label{sec:flow-verification}

The following lemma provides regularity estimates for the velocity depending on the time interval $I_j, j \in \{1,...,J+1\}$.

\begin{lemma}\label{lem:jets}
Let $\alpha>1$ and $P\in\mathcal P_\alpha$. For $1\le j\le J+1$,
set $\Omega_j=((1+R)I_j)\times[-R,R]^q$. For $1\le i\le q$,
write $V_i(u,x)=v^P_{t,i}(x)$ with $t=u/(1+R)$.
On each $\Omega_j$, the function $V_i$ is Lipschitz in time
and space jointly:
\[
 |V_i(u,x)-V_i(u',x')|
 \le C\max\{|u-u'|,\|x-x'\|_\infty\},
 \qquad (u,x),(u',x')\in\Omega_j.
\]
Moreover, for every rectangle $Q\subset\Omega_j$ with side lengths
at most one and every $(u_0,x_0)\in Q$,
\[
 \bigl\|V_i-V_i(u_0,x_0)\bigr\|_{C^\alpha(Q)}\le C.
\]
On the intervals $I_1,\ldots,I_J$, derivatives of all orders
exist. For every multi-index $\nu\in\mathbb N_0^{q+1}$ of order
$k=|\nu|\ge\lceil\alpha\rceil$,
\begin{equation}\label{eq:nonterminal-analytic-derivatives}
 \sup_{(u,x)\in\Omega_j}
 \frac{|D^\nu_{u,x}V_i(u,x)|}{\nu!}
 \le C^{k+1}\lambda_j^{\alpha-k},
 \qquad 1\le j\le J.
\end{equation}
The same $C$ may be chosen for all $P,i,j,R$ and $k$.
\end{lemma}

\begin{proof}
By \eqref{eq:coordinates}, \eqref{eq:holder-gaussian-ratio} and \eqref{eq:population-field},
\begin{equation}\label{eq:jets-physical-formula}
 V_i(u,x)=\frac{2t-1}{r_t^2}x_i
 +\frac{1}{r_t}\mathcal V_i\!\left(
       \frac{1-t}{r_t},\frac{tx}{r_t^2}\right),
 \qquad t=\frac{u}{1+R}.
\end{equation}
Lemma~\ref{lem:velocity-holder} gives the uniform $C^\alpha$ bound
for the second term on every interval. 
In the first term, each derivative of the function $t \mapsto \frac{2t-1}{r_t^2}$ is bounded, hence each derivative of the first term in $u$ gives a factor
$(1+R)^{-1}$, and a derivative in $x_i$ removes $x_i$.
Its derivatives of every fixed positive order are therefore uniformly
bounded. This gives the Lipschitz bound and, on rectangles of side
lengths at most one, the $C^\alpha$ bound after subtracting the
value at one point.

For the higher order derivatives bound, write $H_af(z)=\E[f(z+aG)]$ for $0<a\le1$.
Fix $z$ and let $P_z$ be the Taylor polynomial of $f$ at $z$
of degree $\lceil\alpha\rceil-1$, so that
$(a',z')\mapsto\E[P_z(z'+a'G)]$ is a polynomial of the same
or lower degree. {\it Ergo}, its derivatives of order $k\ge\lceil\alpha\rceil$
vanish. We may therefore subtract $P_z(y)$ when integrating
the corresponding Gaussian kernel derivative against $f(y)$.
Using $|f(y)-P_z(y)|\lesssim|y-z|^\alpha$ gives 
\begin{equation}\label{eq:jets-gaussian-derivatives}
 \sup_{z\in\R^q}
 \frac{|D^\nu_{a,z}H_af(z)|}{\nu!}
 \le C^{k+1}a^{\alpha-k},
 \qquad k=|\nu|\ge\lceil\alpha\rceil .
\end{equation}
At lower orders, differentiating $f$ under the expectation
gives a uniform bound.

Differentiation of the Gaussian kernel gives
$\mathcal V_i(a,z)=a\partial_{z_i}H_af(z)/H_af(z)$.
By the product rule and \eqref{eq:jets-gaussian-derivatives},
the numerator satisfies the same derivative bounds as $H_af$.
Since $H_af\ge m$, differentiating
$H_af\,\mathcal V_i=a\partial_{z_i}H_af$ gives 
\[
 \frac{|D^\nu_{a,z}\mathcal V_i(a,z)|}{\nu!}
 \le C^{k+1}a^{-(k-\alpha)_+},\qquad k=|\nu|\ge0.
\]
Here each product of derivatives of orders $l$ and $k-l$ uses
$(l-\alpha)_++(k-l-\alpha)_+\le(k-\alpha)_+$.

We next bound the derivatives with respect to $u$ and $x$.
Since $r_t^2=t^2+(1-t)^2\ge1/2$ on $[0,1]$, the binomial
series for $1/r_t$ and $1/r_t^2$ gives derivatives of order $l$
bounded by $C^{l+1}l!$. For example, for $l\ge1$,
\[
 \frac{1}{l!}
 \left|\partial_u^l\left(\frac{tx_i}{r_t^2}\right)\right|
 \le C^{l+1}\frac{|x_i|}{(1+R)^l}\le C^{l+1}.
\]
A derivative in $x_i$ removes the factor $x_i$, and all second
spatial derivatives are zero. Thus the Taylor
coefficients of both arguments of $\mathcal V_i$ and of
$(2t-1)x_i/r_t^2$ are bounded by $C^{l+1}$ at order $l$.

On $I_j$, $j\le J$, we have $a=(1-t)/r_t\asymp\lambda_j$.
In \eqref{eq:jets-physical-formula}, the chain rule at order $k$
uses derivatives of $\mathcal V_i$ of orders $s\le k$, and
$a^{-(s-\alpha)_+}\le a^{-(k-\alpha)_+}$. After division by
$\nu!$, products of the Taylor coefficients of the arguments
contribute at most an exponential factor in $k$: for each $s$,
there are $\binom{k-1}{s-1}$ ways to split total degree $k$ into
$s$ positive degrees, and at most $C^k$ choices of coordinates.
Applying the product rule to $1/r_t$ gives
\eqref{eq:nonterminal-analytic-derivatives} for $k\ge\lceil\alpha\rceil$.
\end{proof}
We now use the regularity of Lemma~\ref{lem:jets} to approximate
each coordinate of the population velocity on each time interval
by a GroupSort2 network with the architecture of
Theorem~\ref{th:holder-approximation}.
The resulting networks take advantage of the high regularity
before the terminal interval and have uniform approximation
error of order $n^{-8}$ before the last interval and
$h^\alpha$ on the last interval.
\begin{proposition}\label{prop:blocks}
Let $\alpha>1$ and choose $J,h,R,b_n$ and $B_j=C_B$ as in
\eqref{eq:scales} and \eqref{eq:block-domains}. For each
$1\le j\le J+1$, let $w_j$ be a sufficiently large even width
satisfying \eqref{eq:flow-joint-range}, and let $L_j$ be given
by \eqref{eq:flow-architecture-budget}. For every
$P\in\mathcal P_\alpha$ there exist scalar networks
$g_{\theta_{j,i}}\in\widetilde{\mathcal F}(L_j,w_j,b_n,C_B)$,
$1\le i\le q$, $1\le j\le J+1$, such that, writing
$g_{\theta_j}=(g_{\theta_{j,1}},\ldots,g_{\theta_{j,q}})$,
\begin{equation}\label{eq:block-approximation}
 \sup_{t\in I_j,\,x\in[-R,R]^q}
 \left\|\clip_{[-2R,2R]^q}
 \bigl(g_{\theta_j}((1+R)t,x)\bigr)-v_t^P(x)\right\|_\infty
 \lesssim
 \begin{cases}
 n^{-8},&1\le j\le J,\\
 h^\alpha,&j=J+1.
 \end{cases}
\end{equation}
Each scalar network can be taken to have width $w_j$ in every
hidden layer. Its architecture then has exactly
$s_j=(q+2)w_j+(L_j-1)w_j(w_j+1)+w_j+1$ weights and biases.
\end{proposition} 
\begin{proof}
We fix $P\in\mathcal P_\alpha$, an interval $I_j$ and a
coordinate $i$. We write $\lambda=\lambda_j$, $u=(1+R)t$ and
$\delta=(en)^{-10}$. Let $w=w_j-c_q$, where $c_q$ is an even
integer depending only on $q$, chosen below. Since $w_j$ is
sufficiently large, we may assume $w\ge w_j/2$.
We choose a sufficiently
large fixed integer $\ell_0\ge3$ and set
$r=2^{-\ell_0}\lambda=2^{-j'}$, with $j'=j-1+\ell_0$.

We cover $[-R,R]^q$ by cubes of side length one contained in it
and divide $(1+R)I_j$ into
$\lceil(1+R)|I_j|/r\rceil$ equal intervals.
Since $|I_j|\le\lambda$, there are at most
$\lceil 2^{\ell_0}(1+R)\rceil$ time intervals and hence at most
$C(1+R)^{q+1}$ rectangles. On each rectangle
$Q=[u_0,u_1]\times(x_0+[0,1]^q)$, we define
\[
 V(u,y)=V_i(u,x_0+y)-V_i(u_0,x_0),
 \qquad y\in[0,1]^q.
\]
Its $C^\alpha$ norm is bounded by
Lemma~\ref{lem:jets}, and $u_1-u_0\le r$.

We set $m=\lceil\alpha\rceil-1$. We take $K=m$ if $j=J+1$
and $K=\max\{m,\lceil C_1\log(en)\rceil\}$ if $j\le J$.
For $h_l=2^{-l}$, we use the grids
$\mathcal A_l=(h_l\mathbb Z^q)\cap[0,1]^q$ and the functions
$\theta_{l,a}$ from the proof of
Theorem~\ref{th:holder-approximation}.
We consider
\[
 F(u,y)=\sum_{a\in\mathcal A_{j'}}\theta_{j',a}(y)
 \sum_{b+|\beta|\le K}
 \frac{\partial_u^bD_y^\beta V(u_0,a)}{b!\beta!}
 (u-u_0)^b(y-a)^\beta,
\]
where $b\in\mathbb N_0$ and $\beta\in\mathbb N_0^q$.
The functions $\theta_{j',a}$ sum to one and vanish when
$\|y-a\|_\infty>r$. Taylor's formula therefore gives
$\|V-F\|_\infty\lesssim r^\alpha$ if $j=J+1$.
For $j\le J$, \eqref{eq:nonterminal-analytic-derivatives}
gives
\[
 \|V-F\|_\infty
 \le C\lambda^\alpha(C_2r/\lambda)^{K+1}\le\delta,
\]
choosing first $\ell_0$ and then $C_1$ sufficiently large.

To represent $F$, we reuse the decomposition between
consecutive grids from the proof of
Theorem~\ref{th:holder-approximation}.
For $0\le b\le m$, let
$f_b(y)=\partial_u^bV(u_0,y)/b!$.
We denote by $Q_lf_b$ the approximation defined in that
proof using Taylor polynomials of degree $m-b$ and weights
$\theta_{l,a}$, and set $Q_{-1}f_b=0$. Then
\[
 (u-u_0)^bQ_{j'}f_b
 =
 \sum_{l=0}^{j'}
 \left(\frac{u-u_0}{h_l}\right)^b
 h_l^b(Q_lf_b-Q_{l-1}f_b).
\]
Since $f_b$ has bounded $C^{\alpha-b}$ norm, the calculation
of the coefficients in that proof gives
$O(h_l^{\alpha-b})$ for the spatial difference, expressed
in the functions $\psi_\beta((y-a)/h_l)$.
The factor $h_l^b$ therefore gives coefficients of size
$O(h_l^\alpha)$.

The terms of degree $k>m$ appear only on the grid
of spacing $r$ when $j\le J$. Their rescaled coefficients are
\[
 d_a=\frac{r^k}{b!\beta!}\partial_u^bD_y^\beta V(u_0,a),
 \qquad k=b+|\beta|.
\]
By \eqref{eq:nonterminal-analytic-derivatives},
$|d_a|\le C^{k+1}r^k\lambda^{\alpha-k}
\le C'r^\alpha64^{-k}$, increasing $\ell_0$ if necessary.
For $k\le m$, the factor $64^{-k}$ is included by increasing
the constant, since $m$ is fixed.
Thus, it suffices to represent expressions of the form
\begin{equation}\label{eq:flow-local-polynomial-term}
 \sum_{a\in\mathcal A_l}d_a\theta_{l,a}(y)
 \left(\frac{u-u_0}{s}\right)^b
 \left(\frac{y-a}{s}\right)^\beta,
 \qquad
 s=h_l,\quad k=b+|\beta|,\quad
 |d_a|\le Cs^\alpha64^{-k}.
\end{equation}

We now repeat the construction of the proof
of Lemma~\ref{lem:localized-taylor-sums}.
For $v\in\{0,1,2,3\}^q$, we set
\[
 E_v=\{a\in\mathcal A_l:
          a_i/s-v_i\in4\mathbb Z\text{ for every }i\},
 \qquad
 G(y)=\sum_{a\in E_v}d_a\theta_{l,a}(y).
\]
Within $E_v$, at most one function $\theta_{l,a}$ is
positive at each point. Therefore,
$|G|\le D:=Cs^\alpha64^{-k}$.
The functions
\[
 Y_i(y)=\operatorname{dist}(y_i-(v_i-1)s,4s\mathbb Z)-s
\]
satisfy $|Y_i|\le s$ and $Y_i=y_i-a_i$ where
$\theta_{l,a}(y)>0$ for $a\in E_v$.
Here $\operatorname{dist}(z,4s\mathbb Z)
=\min_{p\in\mathbb Z}|z-4sp|$.
The construction using minima in that proof represents
$Y=(Y_1,\ldots,Y_q)$ with matrices of norm at most one
and $O(1+l)$ layers. For $s=1$, it suffices to use
$Y_i=|y_i-(v_i-1)|-1$.
We define $M(z_0,\ldots,z_q)=z_0^b\prod_{i=1}^q z_i^{\beta_i}$.
The sum of the terms of \eqref{eq:flow-local-polynomial-term}
corresponding to the points $a\in E_v$ satisfies
\[
 \sum_{a\in E_v}d_a\theta_{l,a}(y)
 \left(\frac{u-u_0}{s}\right)^b
 \left(\frac{y-a}{s}\right)^\beta
 =
 G(y)M\!\left(\frac{u-u_0}{s},\frac{Y(y)}s\right).
\]
We take $A=C_3\,8^{-k}(s^{\alpha-1}+(j'+1)^{-1})$.
We apply Lemma~\ref{lem:localized-taylor-sums} with degree
zero to $G$, setting the coefficients outside $E_v$ to zero
and taking error $A\delta$.
Its condition $A\gtrsim D/s$ holds.
We obtain a network that approximates $G/(K_0A)$ with error
$C\delta$ and has all its matrices of norm at most one.
For $C_3$ sufficiently large, $|G/(K_0A)|\le s$.
We clip the output to $[-s,s]$ and call the resulting
network $S$. Clipping adds two hidden layers of width two
and preserves the error and the matrix constraints.
By \eqref{eq:localized-taylor-depth}, $S$ has width
at most $w$ and depth bounded by
\begin{equation}\label{eq:profundidad-flow-prueba}
    C\left\{
 l+1+\log\!\left(2+\frac{D}{A\delta}\right)
 \left(
 1+\frac{s^{-q}}{w^2}
 +\frac{Ds^{-q-1}}{Aw^{1+1/q}}
 \right)
 \right\}. 
\end{equation}

We want to approximate
$G(y)M((u-u_0,Y(y))/s)$.
By Lemma~\ref{lemma:monomial}, there exists a network $p$
that approximates $M$ on $[-1,1]^{q+1}$ with error $\delta$,
satisfies $|p|\le1$ and has an output row of norm at
most $D_k=4^{(k-1)_+}$.
Multiplying its affine constants by $s$ and dividing
its output by $D_k$, we obtain a network that computes
$s\,p(z/s)/D_k$ with all its matrices of norm at most
one, by \eqref{eq:analytic-sorting-homogeneity}.

We apply Lemma~\ref{lemma:multiplication} to the values
$S(y)$ and $s\,p((u-u_0,Y(y))/s)/D_k$, both in $[-s,s]$.
We also multiply the affine constants of this network
by $s$ and divide its output by four.

Let $L_S$ be the depth of $S$. Computing $Y$ requires
$C(1+l)$ layers; computing the monomial and the product requires,
respectively, $C(k+1)\log(e(k+1)/\delta)$ and
$C\log(e/\delta)$ layers. We keep the input and the
necessary results in a fixed number of coordinates.
The resulting network $Z$ has all its matrices of norm
at most one and satisfies
\begin{equation}\label{eq:bound-Z}
    \begin{aligned}
 Z(u,y)
 &=\frac{s}{4}\operatorname{Mult}_\delta\!\left(
       \frac{S(y)}s,
       \frac{p((u-u_0,Y(y))/s)}{D_k}\right),\\
 \text{width}(Z)
 &\le w+C_q,\\
 \text{depth}(Z)
 &\le L_S+C\left\{1+l+
       (k+1)\log\!\left(\frac{e(k+1)}{\delta}\right)\right\},\\
 \left\|4K_0AD_kZ-
       G(y)M\!\left(\frac{u-u_0}{s},\frac{Y(y)}s\right)
       \right\|_\infty
 &\lesssim 4^kA\delta,\\
 4K_0AD_k
 &\asymp 2^{-k}\left(s^{\alpha-1}+\frac1{j'+1}\right).
\end{aligned}
\end{equation}
The error bound follows from
\eqref{eq:relative-multiplication} and the approximation
errors of $S$ and $p$.

For each scale $s=2^{-l}$, each set $E_v$ and each monomial
of degree $k=b+|\beta|$ occurring at that scale, we use one
network $Z$ for the sum over $a\in E_v$ in
\eqref{eq:flow-local-polynomial-term}. Thus there are at most
$4^q\binom{k+q}{q}$ such networks per scale and degree.
Each network $Z$ is multiplied by $4K_0AD_k$ as in \eqref{eq:bound-Z}, where
$A=C_3\,8^{-k}(s^{\alpha-1}+(j'+1)^{-1})$ and
$D_k=4^{(k-1)_+}$.
Since $s=2^{-l}$, this coefficient satisfies
\[
 4K_0AD_k
 \le C2^{-k}\left(2^{-l(\alpha-1)}+\frac1{j'+1}\right).
\]
We denote by $Z_1,\ldots,Z_N$ the constructed networks and by
$\rho_\nu$ the coefficient $4K_0AD_k$ corresponding to
$Z_\nu$. Summing the preceding bound over $l$ and $k$ gives
$\sum_{\nu=1}^N\rho_\nu\lesssim1$, since $\alpha>1$ and
$\sum_{k\ge0}\binom{k+q}{q}2^{-k}<\infty$.

We apply \eqref{eq:dense-level-addition}, using $Z_\nu$
in place of $F_j^\circ$ and $\rho_\nu$ as the weight.
Each update has nonnegative coefficients that sum to
one, so the hidden matrices have norm at most one.
Multiplying the final output by
$\sum_{\nu=1}^N\rho_\nu$, we obtain a network
$\widetilde F=\sum_{\nu=1}^N\rho_\nu Z_\nu$
with bounded output row. Each term of $F$ is approximated
with error at most $C\rho_\nu\delta$.
Adding these errors, we obtain
\[
 \|F-\widetilde F\|_{\infty,[u_0,u_1]\times[0,1]^q}
 \le C\delta\sum_{\nu=1}^N\rho_\nu
 \lesssim\delta.
\]
By \eqref{eq:profundidad-flow-prueba},
\eqref{eq:bound-Z} and \eqref{eq:flow-joint-range},
each network $Z_\nu$ has depth
$O(r^{-q}w^{-2}[\log(en)]^2)$.
Here we use $D/A\lesssim(j'+1)s^\alpha$ and
$j'+1+K\lesssim\log(en)$.
There are $O(j'+1)$ networks corresponding to degrees
at most $m$ and $O((K+1)^{q+1})$ corresponding to
higher degrees. Adding their depths and the layers of
the weighted averages, we obtain
\[
 \text{depth}(\widetilde F)
 \lesssim
 \bigl(j'+1+(K+1)^{q+1}\bigr)
 \frac{r^{-q}}{w^2}[\log(en)]^2
 \lesssim \frac{r^{-q}}{w^2}[\log(en)]^{q+3}.
\]

Since $\widetilde F(u,y)$ approximates
$V_i(u,x_0+y)-V_i(u_0,x_0)$ on $[u_0,u_1]\times[0,1]^q$,
the network
\[
 g_Q(u,x)=V_i(u_0,x_0)+\widetilde F(u,x-x_0)
\]
approximates $V_i$ on $Q$ with the same uniform error.
The translation adds $-W_0(0,x_0)$ to the first-layer bias.
By \eqref{eq:network-constraints},
$\|W_0(0,x_0)\|_\infty\le\|x_0\|_\infty\le R$.
The output bias changes by $V_i(u_0,x_0)$, whose absolute value
is at most $2R$ by \eqref{eq:population-field} and
\eqref{eq:scales}.
We apply Lemma~\ref{lem:join-rectangles} to the family
$\{g_Q\}_Q$, using the Lipschitz bound of Lemma~\ref{lem:jets}.
The depth of the resulting network is
at most
\[
 \sum_Q\bigl(\text{depth}(g_Q)+C_q\bigr).
\]
Since there are at most $C(1+R)^{q+1}$
rectangles, the total depth is at most
\[
 C(1+R)^{q+1}
 \left(1+\frac{r^{-q}}{w_j^2}[\log(en)]^{q+3}\right)
 \le
 C'(1+R)^{q+1}
 \frac{\lambda_j^{-q}}{w_j^2}[\log(en)]^{q+3},
\]
where we use $w\ge w_j/2$, $r=2^{-\ell_0}\lambda_j$ and
$w_j^2\lesssim\lambda_j^{-q}$.
Taking $C_L\ge C'$ in
\eqref{eq:flow-architecture-budget}, this depth
is at most $L_j$.
By \eqref{eq:bound-Z}, each $Z_\nu$ has width at most $w+C_q$.
We argue as in \eqref{eq:dense-level-addition}, with weights
$\rho_\nu$ in place of $A_j$ and networks $Z_\nu$ in place of
$F_j^\circ$. 
On each rectangle $Q$, the network $g_Q$ approximates
$V_i(u,x)=v^P_{u/(1+R),i}(x)$. By Lemma~\ref{lem:join-rectangles},
there exists a single network $g_{\theta_{j,i}}$ such that
\[
  \sup_{t\in I_j,\,x\in[-R,R]^q}
  \bigl|g_{\theta_{j,i}}((1+R)t,x)-v^P_{t,i}(x)\bigr|
  \le \max_Q \sup_{(u,x)\in Q}|g_Q(u,x)-V_i(u,x)|.
\]
Its width is at most of order $\max_Q \operatorname{width}(g_Q)$. The output row has bounded $\ell_1$-norm, and all
biases are $O(1+R)$, as allowed by \eqref{eq:block-domains}.

The error is $O(\delta)$ if $j\le J$ and
$O(r^\alpha+\delta)$ if $j=J+1$.
Since $\delta\le n^{-8}$, and on the last interval
$r\asymp h$ and $\delta\lesssim h^\alpha$, we obtain
\eqref{eq:block-approximation}. Clipping preserves
these errors because $|v^P_{t,i}(x)|\le2R$.
Finally, \eqref{eq:flow-scalar-count} gives the number
of parameters.
\end{proof}

For $\mathcal G_{n,j}=\{g|_{I_j\times\R^q}:g\in\mathcal G_n\}$,
let $\widehat{\mathcal L}_{n,j}$ be \eqref{eq:loss} restricted to $I_j$.
The separate parameters of the blocks give
\begin{equation}\label{eq:empirical-product}
 \begin{aligned}
 \widehat{\mathcal L}_n(g)&=\sum_j\widehat{\mathcal L}_{n,j}(g|_{I_j\times\R^q}),\qquad
 \min_{g\in\mathcal G_n}\widehat{\mathcal L}_n(g)
=\sum_j\min_{g\in\mathcal G_{n,j}}\widehat{\mathcal L}_{n,j}(g).
 \end{aligned}
\end{equation}
Hence the restriction of every minimizer in \eqref{eq:erm} minimizes the corresponding interval loss.

We have already controlled the approximation error of the  class $\mathcal G_n$. 
To analyze its complexity,
we use the following result, which also provides
a uniform Lipschitz bound for the networks. 

\begin{proposition}\label{prop:actual-class}
Every $g\in\mathcal G_n$ satisfies $\sup_{t,x}|g_t(x)|\le2\sqrt q\,R$ and
\begin{equation}\label{eq:actual-envelope}
 \sup_{0\le t\le1}\Lip_x(g_t)\le\sqrt q\,C_B .
\end{equation}
For $0<u<1$,
\begin{equation}\label{eq:actual-entropy}
 \log N(u,\mathcal G_{n,j},\|\cdot\|_\infty)
 \lesssim s_j\log\!\left(\frac{en}{u}\right),
\end{equation}
and, uniformly over $P\in\mathcal P_\alpha$,
\begin{equation}\label{eq:actual-approx}
 \inf_{g\in\mathcal G_{n,j}}\int_{I_j}\|g_t-v_t^P\|_{L^2(P_t)}^2\dd t
 \lesssim |I_j|\{h^{2\alpha}\ind_{\{j=J+1\}}+n^{-8}\}+R^2n^{-8}.
\end{equation}
\end{proposition}
\begin{proof}
Recall that 
$  \chi_j(t)=\min\{1,n^8\operatorname{dist}(t,[0,1]\setminus I_j)\}$
vanishes outside $I_j$ and at every endpoint shared with another
interval. Hence, for $t\in I_j$, all terms with index different
from $j$ vanish in \eqref{eq:flow-network-class}, and
\[
  g_t(x)=\Gamma_j\!\left(t,
  g_{\theta_j}\bigl((1+R)t,\clip_{[-R,R]^q}x\bigr)\right).
\]
For fixed $t$, the map $\Gamma_j(t,\cdot)$ clips each coordinate
to the fixed interval $[-2R\chi_j(t),2R\chi_j(t)]$. Thus
$|g_t(x)|\le 2\sqrt q\,R$. Moreover, both the spatial input
clipping and the scalar output clipping are $1$-Lipschitz in the
maximum norm. Since $B_j=C_B$ in \eqref{eq:block-domains},
\eqref{eq:network-constraints} gives
\[
  |g_{t,i}(x)-g_{t,i}(x')|
  \le C_B\|x-x'\|_\infty,\qquad {\rm for}\qquad 1\le i\le q, 
\]
which proves
\eqref{eq:actual-envelope}. To derive \eqref{eq:actual-entropy}, fix $0<u<1$. The inputs
$((1+R)t,\clip_{[-R,R]^q}x)$ lie in $[-(1+R),1+R]^{q+1}$.
We apply Lemma~\ref{lem:network-entropy} on this cube to each network in \eqref{eq:block-parameter-set}, with covering
radius $u/\sqrt q$. For each $t$, the clipping map
$\Gamma_j(t,\cdot)$ in \eqref{eq:gate} preserves the covering. Taking the product of the $q$ covers
therefore gives a cover of radius $u$ in the uniform Euclidean
norm, with
\[
 \log N(u,\mathcal G_{n,j},\|\cdot\|_\infty)
 \lesssim q\,w_j^2L_j
 \log\!\left(
   \frac{e(q+1+w_j)(L_j+1)b_n C_B(1+R)\sqrt q}{u}
 \right).
\]
By \eqref{eq:scales}, $1+R\lesssim\sqrt{\log(en)}$, while
\eqref{eq:flow-joint-range}--\eqref{eq:block-domains} imply that
$w_j$, $L_j$ and $b_n$ grow at most polynomially in $n$,
uniformly in $j$. Hence the logarithm is bounded by
$C\log(en/u)$. Since $q$ is fixed and $w_j^2L_j\lesssim s_j$
by \eqref{eq:flow-scalar-count}, we prove the claim. 
For \eqref{eq:actual-approx}, choose the networks from
Proposition~\ref{prop:blocks} and set
\[
  E_j=\{t\in I_j:\chi_j(t)<1\}.
\]
By the definition of $\chi_j$, this set is contained in the
intervals of length $n^{-8}$ adjacent to the shared endpoints of
$I_j$. Hence $|E_j|\le 2n^{-8}$. For $t\in I_j\setminus E_j$,
we have $\chi_j(t)=1$, so \eqref{eq:gate} gives exactly the
output clipping used in \eqref{eq:block-approximation}. On
$[-R,R]^q$ the spatial input clipping is also the identity.
Consequently,
\[
  \int_{I_j\setminus E_j}\int_{[-R,R]^q}
  |g_t(x)-v_t^P(x)|^2\,P_t(\dd x)\,\dd t
  \lesssim |I_j|\bigl(n^{-16}
            +h^{2\alpha}\ind_{\{j=J+1\}}\bigr).
\]

It remains to bound the error on $E_j$ and outside $[-R,R]^q$.
By \eqref{eq:coordinates} and \eqref{eq:rho},
\[
  P_t(\dd x)\le \frac{M}{m}r_t^{-q}
  \varphi_q(x/r_t)\,\dd x,\qquad 2^{-1/2}\le r_t\le1.
\]
In particular, the second moments of $P_t$ are uniformly bounded.
Since $|g_t(x)|\le 2\sqrt q\,R$ and
$|v_t^P(x)|\lesssim 1+|x|$ by \eqref{eq:population-field}, we get
\[
  \int_{E_j}\|g_t-v_t^P\|_{L^2(P_t)}^2\,\dd t
  \lesssim R^2|E_j|\lesssim R^2n^{-8}.
\]
Finally, outside $[-R,R]^q$ we have $|x|>R$, so the squared error
is bounded by $C(1+|x|^2)$. Integrating against the Gaussian
bound above gives, uniformly in $t$,
\[
  \int_{\mathbb R^q\setminus[-R,R]^q}
  |g_t(x)-v_t^P(x)|^2\,P_t(\dd x)
  \lesssim e^{-R^2/4}\le n^{-8},
\]
where the last inequality follows from \eqref{eq:scales}.
Integrating this bound over $I_j\setminus E_j$ and adding the
three contributions proves \eqref{eq:actual-approx}.
\end{proof}
\subsection{Upper bound}\label{sec:upper}

\begin{proof}[Proof of \eqref{eq:flow-fixed-orders}]
Fix $P\in\mathcal P_\alpha$. By \eqref{eq:model} and Markov's inequality,
\[
 \Prb_P(|Y|>r)
{ \le \frac{M}{m}\Prb(|Z|>r)
 \le \frac{M}{m}e^{-\frac{r^2}{4}}\E e^{\frac{|Z|^2}{4}}
 =\frac{M}{m}2^{\frac{q}{2}}e^{-\frac{r^2}{4}},}
 \qquad Z\sim N(0,I_q),\quad r\ge0.
\]
By \eqref{eq:empirical-product}, $\widehat v_n$ minimizes
the empirical loss on each interval. Applying
Lemma~\ref{prop:interval-oracle} with
\eqref{eq:actual-entropy} and \eqref{eq:actual-approx} gives
\begin{equation}\label{eq:risk-block}
\begin{aligned}
 \E_P\int_{I_j}
 \|\widehat v_{n,t}-v_t^P\|_{L^2(P_t)}^2\dd t
 &\lesssim |I_j|\left(
 h^{2\alpha}\ind_{\{j=J+1\}}
 +\frac{s_j[\log(en)]^2}{n}+n^{-8}\right)
 +R^2n^{-8}.
\end{aligned}
\end{equation}

The population flow has law $P_t$ at time $t$
by Lemma~\ref{lem:regularity}.  Applying Gronwall's
inequality, as in the
proof of Lemma~1 of \citet{Stephanovitch2026}, gives
\begin{equation}\label{eq:stability}
 W_2(\widehat P_n,P)
 \le\int_0^1
 \exp\!\left(\int_t^1\Lip_x(\widehat v_{n,s})\dd s\right)
 \|\widehat v_{n,t}-v_t^P\|_{L^2(P_t)}\dd t.
\end{equation}
By \eqref{eq:actual-envelope}, the exponential is bounded
uniformly over the fitted fields. Applying Cauchy--Schwarz
on each interval and then taking expectations yields
\[
 \E_P W_2(\widehat P_n,P)
 \lesssim\sum_{j=1}^{J+1}|I_j|^{\frac{1}{2}}
 \left(
 \E_P\int_{I_j}
 \|\widehat v_{n,t}-v_t^P\|_{L^2(P_t)}^2\dd t
 \right)^{\frac{1}{2}}.
\]
The last interval has length $h$, while
$|I_j|\le\lambda_j$ and
$\sum_j|I_j|^{\frac{1}{2}}\lesssim1$.
Substituting \eqref{eq:risk-block}, using that the last interval length is $h$ and
$s_j\asymp w_j^2L_j$ therefore gives
\[
 \E_P W_2(\widehat P_n,P)
 \lesssim h^{\alpha+1}
 +\frac{\log(en)}{\sqrt n}
   \sum_{j=1}^{J+1}\lambda_jw_j\sqrt{L_j}
 +Rn^{-4}.
\]
By \eqref{eq:scales},
\[
 h^{\alpha+1}\asymp
 \begin{cases}
 n^{-\frac{1}{2}},&q\le2,\\
 n^{-\frac{\alpha+1}{2\alpha+q}},&q>2,
 \end{cases}
 \qquad {\rm and}\qquad  Rn^{-4}\lesssim h^{\alpha+1}.
\]
By \eqref{eq:flow-architecture-budget}, the summand $\lambda_jw_j\sqrt{L_j}$ is bounded by $C[\log(en)]^{(3q+7)/4}\lambda_j^{1-q/2}$. The geometric sum satisfies
\begin{equation}\label{eq:upper-geometric-sum}
 \sum_{j=1}^{J+1}\lambda_j^{1-q/2}\asymp
 \begin{cases}
 1,&q=1,\\
 J+1,&q=2,\\
 h^{1-q/2},&q>2 .
 \end{cases}
\end{equation}
Since $J+1\lesssim\log(en)$ and $n^{-1/2}h^{1-q/2}\asymp h^{\alpha+1}$ for $q>2$, substitution gives \eqref{eq:flow-fixed-orders}.

\end{proof}
\subsection{Minimax lower bound}\label{sec:lower}

\begin{proof}[Proof of the lower bound in
Theorem~\ref{thm:main-flow}]
Let $\Phi$ be the standard normal cdf and set
$G(x)=(\Phi(x_1),\ldots,\Phi(x_q))$.
Then $G_\#N(0,I_q)=U_q$ and
$\Lip(G)\le(2\pi)^{-1/2}$.

For $\eta>0$, let $\mathcal Q$ be the class of
probabilities on $[0,1]^q$ with densities $\varrho$
satisfying $\|\varrho-1\|_{C^\alpha}\le\eta$.
The lower-bound construction in
\citet[Proof of Theorem~3]{NilesWeedBerthet2022}
applies to this class for every fixed $\eta>0$.
Consequently,
\begin{equation}\label{eq:lower-nwb}
 \inf_{\widetilde Q_n}\sup_{Q\in\mathcal Q}
 \E_Q W_1(\widetilde Q_n,Q)
 \ge c n^{-\min\left\{\frac{1}{2},\frac{\alpha+1}{2\alpha+q}\right\}}.
\end{equation}

For $Q\in\mathcal Q$, define
$P_Q=(G^{-1})_\#Q$. Its density is
$p_Q(x)=\varphi_q(x)\varrho(G(x))$.
Since $G$ has bounded derivatives of every order,
choosing $\eta$ sufficiently small ensures
$P_Q\in\mathcal P_\alpha$ for all $Q\in\mathcal Q$.

Now $G_\#P_Q=Q$, and the Lipschitz property of $G$
gives
$ W_1(P_Q,P_{Q'})
 \ge\sqrt{2\pi}\,W_1(Q,Q').$ 
Moreover, applying $G^{-1}$ to independent
observations from $Q$ gives independent observations
from $P_Q$. Applying $G$ to any estimator of $P_Q$
therefore gives an estimator of $Q$.
Consequently,
\[
 \inf_{\widetilde P_n}\sup_{P\in\mathcal P_\alpha}
 \E_P W_1(\widetilde P_n,P)
 \ge\sqrt{2\pi}
 \inf_{\widetilde Q_n}\sup_{Q\in\mathcal Q}
 \E_Q W_1(\widetilde Q_n,Q).
\]
The assertion follows from \eqref{eq:lower-nwb}.
\end{proof}

\subsection{Combining approximations on rectangles}
The following lemma combines approximations on overlapping
rectangles, preserving their uniform error and a bounded
output row norm. We use it in the proof of
Proposition~\ref{prop:blocks} to combine the local
approximations in time and space into a single network
for each time interval.
\begin{lemma}\label{lem:join-rectangles}
Let $\Omega\subset\R^q$ be a rectangle and let
$\mathcal Q$ be a finite family of nonempty closed
rectangles such that $\Omega=\bigcup_{Q\in\mathcal Q}Q$.
Let $V:\Omega\to\R$ satisfy $\Lip_\infty(V)\le K$.
For each $Q\in\mathcal Q$, let $g_Q$ be a GroupSort2 network
of width $w_Q$ and depth $L_Q$, with hidden matrices
of norm at most one and output row of norm at most $B$,
such that $\|g_Q-V\|_{\infty,Q}\le\delta$.
We define
\[
 g(y)=\min_{Q\in\mathcal Q}
 \left\{g_Q(\clip_Q y)
       +K\|y-\clip_Q y\|_\infty\right\},
 \qquad y\in\R^q.
\]
Then $\|g-V\|_{\infty,\Omega}\le\delta$.
Moreover, $g$ admits a GroupSort2 representation
of width $w$ and depth $L$ satisfying
\[
 w\le\max_{Q\in\mathcal Q}w_Q+C_q,
 \qquad
 L\le\sum_{Q\in\mathcal Q}(L_Q+C_q).
\]
The hidden matrices of this representation have norm
at most one and its output row has norm at most $1+B+K$.
\end{lemma}
\begin{proof}
Lipschitz continuity gives
\[
 g_Q(\clip_Q y)+K\|y-\clip_Q y\|_\infty
 \ge V(y)-\delta,\qquad y\in\Omega.
\]
A rectangle containing $y$ provides the upper bound
$V(y)+\delta$. To compute the minimum, we keep $y$ and the minimum
computed so far as pairs of equal coordinates.
Clipping each coordinate uses maxima and minima.
The distance is the maximum of zero and the affine
functions that measure how far each coordinate exceeds
the endpoints of its interval. We divide each expression entering the minimum by
$1+B+K$ before comparing it with the minimum already computed.
Its row then has absolute sum at most one.
The final row restores the common factor.
The additional affine constants are bounded by
the maximum of one and the absolute values of the
rectangle endpoints.
\end{proof}

\section{Proofs for regression}\label{sec:regression-proofs}
\begin{proof}[Proof of Theorem~\ref{thm:regression}]
For every bounded predictor $g$, conditioning on $X$ gives
\[
 \E(g(X)-Y)^2-\E(f_0(X)-Y)^2
 =\|g-f_0\|_{L^2(P_X)}^2.
\]
Apply Lemma~\ref{lem:bounded-regression} with observed
variable $(X,Y)$, predictor input $X$, response $Y$,
predictor bound $H$ and response bound $M_Y$.
With $u=(en)^{-2}$, exact empirical minimization gives
\[
 \E\|\widehat f_n-f_0\|_{L^2(P_X)}^2
 \le3\inf_{g\in\mathcal F_n}
       \|g-f_0\|_{L^2(P_X)}^2
 +C\frac{
 \log\!\left(2\mathcal N(u,\mathcal F_n,
              \|\cdot\|_{\infty,[0,1]^q})\right)+1}{n}
 +Cn^{-2}.
\]
Theorem~\ref{th:holder-approximation} bounds the infimum
by $C\epsilon_n^2$.
By \eqref{eq:holder-network-budgets},
Lemma~\ref{lem:network-entropy} and
\eqref{eq:regression-joint-budget},
\[
 \log\!\left(2\mathcal N(u,\mathcal F_n,
              \|\cdot\|_{\infty,[0,1]^q})\right)+1
 \lesssim S(L_n, w_n)\log(en)
 \lesssim\epsilon_n^{-\frac{q}{\alpha}}\log(en).
\]
Consequently,
\begin{equation}\label{eq:regression-balance}
 \E\|\widehat f_n-f_0\|_{L^2(P_X)}^2
 \le C\epsilon_n^2
 +\frac{C\epsilon_n^{-\frac{q}{\alpha}}\log(en)}{n}
 +Cn^{-2}.
\end{equation}
For all sufficiently large $n$,
\eqref{eq:regression-resolution} gives
$\epsilon_n^{2+\frac{q}{\alpha}}=\frac{\log(en)}{n}$.
The second term in \eqref{eq:regression-balance} is
therefore at most $C\epsilon_n^2$, and $n^{-2}$ is smaller.
This proves \eqref{eq:regression-risk}. 
\end{proof}

\begin{proof}[Proof of Corollary~\ref{cor:regression-subgaussian}]
By \eqref{eq:subgaussian-model} and $\|f_0\|_\infty\le H$, we have $ \Prb(|Y|>r)\le2e^{-\frac{r^2}{8\sigma^2}},$ for  $r\ge2H.$ 
Integrating this bound and applying the union bound, we
can choose $K_n=C\sqrt{\log(en)}\ge2H$ such that
\begin{equation}\label{eq:regression-tail-truncation}
 \Prb\left(\max_{1\le i\le n}|Y_i|>K_n\right)\le n^{-3},
 \qquad
 \E\bigl[(K_n+|Y|)^2\ind_{\{|Y|>K_n\}}\bigr]\le n^{-4}.
\end{equation}

Set $E_n=\bigcap_{i=1}^n\{|Y_i|\le K_n\}$ and
$P^K=P_{X,Y}(\,\cdot\mid |Y|\le K_n)$.
Since the observations are independent and each condition
in $E_n$ concerns only one observation,
\begin{equation}\label{eq:regression-conditional-sample}
 \operatorname{Law}\bigl((X_i,Y_i)_{i=1}^n\mid E_n\bigr)
 =(P^K)^{\otimes n}.
\end{equation}
Thus, conditionally on $E_n$, the observations are iid
with law $P^K$ and responses bounded by $K_n$.
Write $P_X^K$ for its covariate marginal and
$f_K(x)=\E_{P^K}[Y\mid X=x]$ for its regression function,
choosing a version bounded by $K_n$.

Every function in the class $\mathcal F_n$ is  bounded by
$H$, and $\widehat f_n$ still minimizes the empirical
squared loss under this conditional law.
We may therefore apply
Lemma~\ref{lem:bounded-regression} under
$\Prb(\,\cdot\mid E_n)$, with regression function $f_K$,
predictor bound $H$ and response bound $K_n$.
Taking $u=(en)^{-5}$ and writing
$N=\mathcal N(u,\mathcal F_n,\|\cdot\|_{\infty,[0,1]^q})$
gives
\begin{equation}\label{eq:regression-conditional-oracle}
\begin{aligned}
 \E\left[
 \|\widehat f_n-f_K\|_{L^2(P_X^K)}^2
 \,\middle|\,E_n\right]
 &\le3\inf_{g\in\mathcal F_n}
       \|g-f_K\|_{L^2(P_X^K)}^2\\
 &\quad+C(H+K_n)^2\frac{\log(2N)+1}{n}
       +C(H+K_n)u.
\end{aligned}
\end{equation}

We next compare the conditional risk with the original
risk. For a measurable function $a$, write
$\mathcal R(a)=\E(a(X)-Y)^2$ and
$\mathcal R_K(a)=\E_{P^K}(a(X)-Y)^2$.
For $\|a\|_\infty\le K_n$, splitting the expectation
according to $|Y|\le K_n$ gives
\[
\begin{aligned}
 |\mathcal R(a)-\mathcal R_K(a)|
 &\le
 \E\bigl[(K_n+|Y|)^2\ind_{\{|Y|>K_n\}}\bigr]
 +4K_n^2\Prb(|Y|>K_n)\le2n^{-4}.
\end{aligned}
\]
Both terms are controlled by
\eqref{eq:regression-tail-truncation}, since
$(K_n+|Y|)^2\ge4K_n^2$ when $|Y|>K_n$.

The functions $f_0$ and $f_K$ minimize their respective
risks over all measurable functions bounded by $K_n$.
Taking infima in the preceding comparison gives
$|\mathcal R(f_0)-\mathcal R_K(f_K)|\le2n^{-4}$.
Subtracting these minimum risks yields
\begin{equation}\label{eq:regression-excess-comparison}
 \left|
 \|g-f_0\|_{L^2(P_X)}^2
 -
 \|g-f_K\|_{L^2(P_X^K)}^2
 \right|
 \le4n^{-4},
 \qquad  \|g\|_\infty\le H.
\end{equation}
This bound holds for every function in $\mathcal F_n$,
so it also holds for each realization of $\widehat f_n$. By \eqref{eq:holder-approximation} and
\eqref{eq:regression-excess-comparison},
\[
 \inf_{g\in\mathcal F_n}
 \|g-f_K\|_{L^2(P_X^K)}^2
 \le C\epsilon_n^2+4n^{-4}.
\]
Lemma~\ref{lem:network-entropy}, the parameter bounds
\eqref{eq:holder-network-budgets} and the assumed
bound on $w_n^2(L_n+1)$ give
$ \log(2N)
 \lesssim\epsilon_n^{-\frac{q}{\alpha}}\log(en).$ 
Using $(H+K_n)^2\lesssim\log(en)$ in
\eqref{eq:regression-conditional-oracle}, and then applying
\eqref{eq:regression-excess-comparison}, we obtain
\[
 \E\left[
 \|\widehat f_n-f_0\|_{L^2(P_X)}^2
 \,\middle|\,E_n\right]
 \le C\epsilon_n^2
 +\frac{C\epsilon_n^{-\frac{q}{\alpha}}
               [\log(en)]^2}{n}
 +Cn^{-4}.
\]
On $E_n^c$, the error is bounded by $4H^2$.
Multiplying the conditional bound by $\Prb(E_n)$ and
using \eqref{eq:regression-tail-truncation} therefore gives
\begin{equation}\label{eq:regression-subgaussian-balance}
 \E\|\widehat f_n-f_0\|_{L^2(P_X)}^2
 \le C\epsilon_n^2
 +\frac{C\epsilon_n^{-\frac{q}{\alpha}}
               [\log(en)]^2}{n}
 +Cn^{-3}.
\end{equation}
For all sufficiently large $n$,
\eqref{eq:regression-subgaussian-resolution} gives
$\epsilon_n^{2+\frac{q}{\alpha}}
=\frac{[\log(en)]^2}{n}$. 
This proves \eqref{eq:regression-subgaussian-rate}. 
\end{proof}

\bibliographystyle{plainnat}
\bibliography{biblio}

\end{document}